\documentclass[10pt]{article}

\usepackage{amsmath}
\usepackage{amssymb}
\usepackage{amsthm}
\usepackage{graphicx}

\newtheorem{theorem}{Theorem}

\newcommand{\eq}{{\sim}^z}
\newcommand{\txt}{\textrm}
          \newtheorem{defin}{Definition}[section]
          
          \newtheorem{teo}{Theorem}[section]
          
          \newtheorem{con}{Conjecture}
          \newtheorem{cond}{Condition}
          \newtheorem{prop}[teo]{Proposition}
          \newtheorem{lem}{Lemma}[section]
          \newtheorem{rmk}[teo]{Remark}
          \newtheorem{cor}{Corollary}[section]
          
          \newcommand{\beq}{\begin{equation}}
          \newcommand{\eeq}{\end{equation}}
          \newcommand{\beqn}{\begin{eqnarray}}
          \newcommand{\beqnn}{\begin{eqnarray*}}
          \newcommand{\eeqn}{\end{eqnarray}}
          \newcommand{\eeqnn}{\end{eqnarray*}}
          \newcommand{\bprop}{\begin{prop}}
          \newcommand{\eprop}{\end{prop}}
          \newcommand{\bteo}{\begin{teo}}
          \newcommand{\bcor}{\begin{cor}}
          \newcommand{\ecor}{\end{cor}}
          \newcommand{\bcon}{\begin{con}}
          \newcommand{\econ}{\end{con}}
          \newcommand{\bcond}{\begin{cond}}
          \newcommand{\econd}{\end{cond}}
          \newcommand{\eteo}{\end{teo}}
          \newcommand{\brm}{\begin{rmk}}
          \newcommand{\erm}{\end{rmk}}
          \newcommand{\blem}{\begin{lem}}
          \newcommand{\elem}{\end{lem}}
          \newcommand{\ben}{\begin{enumerate}}
          \newcommand{\een}{\end{enumerate}}
          \newcommand{\bei}{\begin{itemize}}
          \newcommand{\eei}{\end{itemize}}
          \newcommand{\bdf}{\begin{defin}}
          \newcommand{\edf}{\end{defin}}
\newcommand{\ro}{{_{[1]}r_0}}
\newcommand{\rz}{{_{[1]}r_z}}
\newcommand{\brz}{{_{[1]}\bar{r}_z}}
\newcommand{\bro}{{_{[1]}\bar{r}_0}}
\newcommand{\nro}{{_{[n]}r_0}}
\newcommand{\Bo}{{_{[1]}B_0}}
\newcommand{\lo}{{_{[1]}l_0}}
\newcommand{\Bz}{{_{[1]}B_0}}
\newcommand{\Bzz}{{_{[1]}B_z}}

          \newcommand{\fr}{\frac}
          \renewcommand{\r}{{\mathbb R}}
          \newcommand{\br}{\bar{\mathbb R}}
          \newcommand{\Z}{{\mathbb Z}}

          \newcommand{\R}{{\mathbb R}}
          
          \newcommand{\E}{{\mathbb E}}
          
          \renewcommand{\P}{{\mathbb P}}
          \newcommand{\N}{{\mathbb N}}

	\newcommand{\M}{{\cal M}}
          \newcommand{\W}{{\cal W}}

          \newcommand{\h}{{\cal H}}
          \newcommand{\f}{{\cal F}}

\newcommand{\Mun}{\mathcal{M}(\tau_1)}
\newcommand{\Mde}{\mathcal{M}(\tau_2)}
\newcommand{\Lo}{\mathcal{L}}
\newcommand{\Bsw}{\left[\sim^z_\e B_{switch}\right] }
\newcommand{\Bol}{\left[\sim^z_\e B \right] }
\newcommand{\nB}{\nparallel^z_\e}

          \newcommand{\g}{\gamma}

          \newcommand{\D}{\Delta}
          \newcommand{\e}{\epsilon}

          \newcommand{\s}{\sigma}
          \renewcommand{\o}{\Pi}

\newcommand{\dd}{\cdot}
\newcommand{\btt}{\begin{theorem}}
\newcommand{\ett}{\end{theorem}}
\newcommand{\Net}{\mathcal{N}}
          
\newcommand{\daw}{\downarrow}
\newcommand{\uaw}{\uparrow}
\newcommand{\raw}{\rightarrow}

\newcommand{\Wl}{\mathcal{W}_l}
\newcommand{\Wr}{\mathcal{W}_r}
\newcommand{\be}{\begin{equation}}
\newcommand{\ee}{\end{equation}}

          \newcommand\sqr{\vcenter{
          \hrule height.1mm
          \hbox{\vrule width.1mm height2.2mm\kern2.18mm\vrule width.1mm}
          \hrule height.1mm}}        % This is a slimmer sqr.

\begin{document}
\title{Marking (1,2) Points of the Brownian Web and Applications}
\author{
{\bf C.~M.~Newman}
{\small Courant Inst.~of Mathematical Sciences, NYU, New York, NY 10012}
\and
{\bf K.~Ravishankar}
{\small  Dept.~of Mathematics, SUNY College at New Paltz, New Paltz, NY 12561}
\and
{\bf E.~Schertzer}
{\small Dept.~of Mathematics, Columbia University, New York, NY 10027}
}

%\date{\small\file\today}
\maketitle

\begin{abstract}
The Brownian web (BW), which developed from the work of Arratia and then
T\'{o}th and Werner, is a random collection of paths (with specified starting
points) in one plus one dimensional space-time that arises as the scaling
limit of the discrete web (DW) of coalescing simple random walks. Two recently
introduced extensions of the BW, the Brownian net (BN) constructed by Sun and
Swart, and the dynamical Brownian web (DyBW) proposed by Howitt and Warren,
are (or should be) scaling limits of corresponding discrete extensions of the
DW --- the discrete net (DN) and the dynamical discrete web (DyDW). These
discrete extensions have a natural geometric structure in which the underlying
Bernoulli left {\it or\/} right ``arrow'' structure of the DW is extended by means of
branching (i.e., allowing left {\it and\/} right simultaneously) to construct the 
DN or by means of switching (i.e., from left to right and vice-versa) to
construct the DyDW. In this paper we show that there is a similar structure 
in the continuum where arrow direction is replaced by the left or right  
parity of the (1,2) space-time points of the BW (points with one incoming 
path from the past and two outgoing paths to the future, only one of which is 
a continuation of the incoming path). We then provide a complete construction 
of the DyBW and an alternate construction of the BN to that of Sun and Swart 
by proving that the switching or branching can be implemented by a Poissonian 
marking of the (1,2) points.

Le r\'{e}seau Brownien (BW) construit \`{a} partir des travaux de Arratia,
de T\`{o}th et de Werner est une collection
al\'{e}atoire de chemins
(avec des points de depart determin\'{e}s) 
dans un espace deux-dimensionnel (une dimension en
temps et une autre en espace), 
qui est la limite d'\'{e}chelle d'un r\'{e}seau discret (DW) 
de marches al\'{e}atoires coalescentes.
R\'{e}cemment, deux extensions du BW ont \'{e}t\'{e} introduites: 
le filet Brownien (BN), construit par Sun et Swart, et le r\'{e}seau Brownien
dynamique (DyBW), propos\'{e} par Howitt et Warren. Ces deux 
objets sont (ou devraient \^{e}tre) 
la limite d'\'{e}chelle de deux extensions naturelles du r\'{e}seau 
discret---le filet discret (DN)
et le r\'{e}seau dynamique discret (DyDW). Le DN et le DyDW
sont obtenus par une modification de la configuration des ``fl\`{e}ches''
droites {\it ou\/} gauches
qui composent le r\'{e}seau discret. 
Pour le DN, un m\'{e}canisme de ramification est introduit (en permettant
des fl\`{e}ches droites {\it et\/} gauches simultan\'{e}ment) alors que 
pour le DyDW, la direction des fl\`{e}ches est modifi\'{e}e 
(de droite \`{a} gauche et vice-versa).
Dans cet article, nous montrons qu'il existe une structure g\'{e}om\'{e}trique
analogue dans le cas continu. Plus pr\'{e}cis\'{e}ment, la direction
des fl\`{e}ches dans le cas discret est remplac\'{e}e par la direction
des points $(1,2)$ du r\'{e}seau Brownien  (en un point $(1,2)$ se trouvent
un chemin entrant et deux chemins sortants,
l'un d'eux \'{e}tant la continuation du chemin entrant).
Nous montrons que 
les ramifications et changements de direction 
peuvent \^{e}tre
introduits dans le cas continu par 
un marquage de type Poisson des points $(1,2)$. Par l'interm\'{e}diaire 
de ce marquage,
nous donnons une construction compl\`{e}te du DyBW et
une construction
alternative \`{a} celle de Sun et Swart du BN.

\end{abstract}


%\newpage
%\tableofcontents 
%\newpage
\section{Introduction}
\label{intro}
In ~\cite{FNRS07}, the present authors and L.~R.~Fontes obtained
some results about exceptional times for a dynamical model of coalescing
one-dimensional random walks (the ``dynamical discrete web" (DyDW)).
Underlying those results was the idea that there should
be a natural continuum limit of the DyDW, the ``dynamical Brownian web"
(DyBW) for which corresponding results
would be valid,
provided such a continuum system actually exists. The DyBW was also proposed in
a paper of Howitt and Warren~\cite{HW07}, where the DyDW was first discussed, and some
of its properties were analyzed, assuming its existence.

The main purpose of the present paper is to develop
a Poissonian marking of
certain nongeneric points (called (1,2) points, as
we will explain) in the (static) Brownian web (BW) which
we then use to give
the first complete construction of the DyBW.
In a revised version \cite{FNRS08} of \cite{FNRS07}, 
this construction will be used 
to argue that exceptional time results
derived earlier
for the DyDW  should extend to the DyBW.
As we shall see, this marking technology is
natural and has other applications besides the DyBW.
One of those, which we explore in detail in this paper, is
an alternative construction of the ``Brownian net"
(BN) of Sun and Swart~\cite{SS07}. A future application \cite{NRS08}, which we
discuss briefly in Subsection \ref{nucleation} below, is
to scaling limits of
one-dimensional voter models
in which there is ``nucleation along boundaries.'' That will extend,
in a nontrivial way, earlier work~\cite{FINR05} on scaling limits
in which nucleation ``in the bulk'' was treated by using marking
of nongeneric (0,2) points of the BW, which are simpler to
deal with than (1,2) points.
Another
model closely related to the marking of the Brownian web 
is a class of stochastic
flows of kernels introduced by Howitt and Warren \cite{HW07}.
This is the subject of
ongoing work \cite{ScheSS08.2}.

In addition to direct applications of Poissonian markings
of  BW (1,2) points, we believe that these constructions
are of interest as special examples of an approach that
is relevant beyond the  Brownian web setting.
Indeed, the idea of using Poissonian marking of nongeneric
double points in the context of the Schramm-Loewner Evolution $SLE(6)$, was
proposed in~\cite{CFN06a, CFN06b}
as an approach to the continuum scaling limits
of both ``near-critical'' and dynamical two-dimensional
percolation models. In that setting, the critical scaling limit
is analogous to the BW, dynamical percolation to  the DyDW
and near-critical percolation to a discrete web with
small nonzero drift. Progress in applying that approach
has been reported by Garban, Pete and Schramm \cite{GPS07,G09}; 
for other results on scaling limits of near-critical
percolation, see \cite{N07,NW07,CJM08}.

%\section{Introduction-Background}
%\subsection{The Problem}

\subsection{Arrows, Switching and Branching}
\label{arrows}

\begin{figure}
\centering
\includegraphics[scale=0.3]{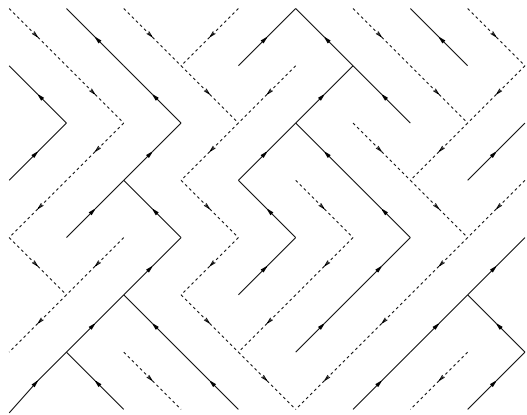}
\caption{Forward coalescing random walks 
(full lines) and their
dual backward walks (dashed lines).
}
\label{explain2}
\end{figure}

{\bf The Discrete Web.} The discrete web is a collection of coalescing 
one-dimensional simple random walks starting 
from every
point in the discrete 
space-time domain 
$\Z_{even}^2=\{(x,t)\in\Z^2: x+t \ \txt{is even}\}$. 
The Bernoulli percolation-like structure is highlighted by defining $\xi_{x,t}$ for
$(x,t)\in \mathbb{Z}_{even}^2$ to be the increment of the random walk
at location $x$ at time $t$. These Bernoulli variables
are symmetric and independent and the paths of all the coalescing 
random
walks can be reconstructed by assigning to each point $(x,t)$
an arrow from  $(x,t)$ to $\{x+\xi_{x,t},t+1\}$
and
considering
all the paths starting from arbitrary points in $\mathbb{Z}_{even}^2$
that follow the arrow configuration $\aleph$.
We note that there is also a set of dual (or backward) paths defined
by the same $\xi_{x,t}$'s but
with arrows from $(x,t+1)$ to $(x-\xi_{x,t},t)$.
The collection of all dual paths is a system of backward (in time) coalescing
random walks that do not cross any of the forward paths (see Figure \ref{explain2}).

There are two natural
variants of the discrete web; one is the dynamical discrete web (DyDW)
which involves {\it switching\/} of arrows and the other is the
discrete net (DN) which involves {\it branching\/}
(or equivalently, adding) of arrows. Each of
these is constructed by a straightforward modification of the arrow
structure in the standard discrete web. The
essence of this paper is a construction of analogous
modifications in the continuum space-time setting.

{\bf The Dynamical Discrete Web.} 
In the DyDW, there is, in addition to the random
walk discrete time parameter, an additional (continuous)
dynamical time
parameter $\tau$. The system starts at $\tau=0$ as an ordinary
DW and then evolves in $\tau$ by randomly switching the direction of each
arrow at a fixed rate (say $\lambda$),
independently of all other arrows.  This
naturally defines a dynamical arrow configuration
$\tau\rightsquigarrow\aleph(\tau)$
. If one follows the arrows
starting from the (space-time) origin at $(0,0)$, this begins
at $\tau=0$ as a simple symmetric random walk and then evolves dynamically
in $\tau$ in a different way than the ``dynamical random walks''
studied in~\cite{BHPS03}. As noted in~\cite{FNRS07}, the nature
of exceptional dynamical times is quite different in this
situation than in that of~\cite{BHPS03}.
For example,
the dynamical random walk constructed from 
$\aleph(\tau)$ violates the law of the iterated logarithm on
a set of $\tau$'s of Hausdorff dimension one.

{\bf The Discrete Net.}
In the DN,
space-time points have at any point arrows
of {\it both\/} directions with probability $p\in [0,1]$,
independently of
other points --- i.e., individual points have either both directions
(with probability $p$), corresponding to points where there is
branching of paths, or only a left arrow (with probability $(1-p)/2$)
or only a right arrow (with probability $(1-p)/2$).
The DyDW and DN models
may be coupled together by taking the DyDW,
declaring that there are both arrows at a point if at least
one switch occurred up to dynamical time $\tau$ and otherwise
declaring that
there is only one arrow whose direction is that of the DyDW
at dynamical time $0$. This yields the DN with $p = 1 - e^{-\lambda \tau}$.

Under diffusive scaling, individual
random walk paths converge to Brownian motions and the entire
collection of discrete paths in the DW converges in an appropriate
sense (see~\cite{FINR04}) to the continuum Brownian web (BW).
We review in Section~\ref{Brownian_Web} some of the
basic features of the BW, which developed from the work
of Arratia~\cite{A81} and of T\'{o}th and Werner~\cite{TW98},
but meanwhile we briefly comment
on its structure.
The BW is a random collection of paths (with specified starting points)
in continuum space-time with one or more paths starting from
every point. Furthermore, although generic (e.g., deterministic)
space-time points have only
$m_{\mbox{{\scriptsize out}}}=1$
outgoing (to later times)
paths from that point
and $m_{\mbox{{\scriptsize in}}}=0$
incoming paths passing through that point (from earlier times),
there are non-generic points with other values of
$(m_{\mbox{{\scriptsize in}}},m_{\mbox{{\scriptsize out}}})$.
In this paper, a dominant role is played by the $(1,2)$ points as
we shall explain.

It is natural that there
should also exist scaling limits of the DyDW (including of
the random walk from the origin
evolving in $\tau$) and of the DN
(with appropriate scaling of $\tau$ and $p$ along with
space-time). Indeed, this has
been studied by Sun and Swart~\cite{SS07} for the case of the net
and by Howitt and Warren~\cite{HW07} for the case of the dynamical web.
%(and its path from $(0,0)$).
The focus of this paper is on how
to construct these continuum objects directly
from the BW in a way that parallels the discrete construction.
%What makes the answer to this question less than obvious at
%first glance
A priori, this appears difficult since
the discrete construction is entirely
based on modifying the discrete arrow structure of the DW,
while in the BW it is unclear whether there even is
any arrow structure to modify.

The main themes of this paper are thus: ``Where is the
arrow structure of the BW?'' and ``How is it modified
to yield the BN and the DyBW (including a dynamically
evolving Brownian motion
from the origin)?''. As we will see, the answer to
the first question is that the arrow structure of the BW comes
from the $(1,2)$ points, each of which is equipped with
a left or right parity according to which of the two
outgoing paths is the
continuation of the single incoming
path --- see Figure \ref{fig12} below.
The answer to the second question is based
on a Poissonian marking of the
$(1,2)$ points, which can then be used either to
create  branching or to switch parity at marked points.

%\subsection{Connection with some nucleation model}
\subsection{Nucleation on Boundaries}
\label{nucleation}

The discrete-time one-dimensional voter model starts at time zero
with colors assigned to each odd integer site and then evolves
in time by assigning a color to the space-time point
$(i,j+1)$ with $i+j+1$ odd as that of $(i-1,j)$ or $(i+1,j)$ with probability
$1/2$ each, independently of other space-time points. The genealogy
of colors (looked at backwards in time) is described by coalescing
random walks (on these odd space-time points)
regardless of the initial state of the system. One often
considers the case where there are $q$ possible colors ($q = 2,3, \dots$);
%and the initial colors are chosen uniformly and indepently;
then the
boundaries between sites of different colors evolve forward in time (on the even space-time points) 
--- in the case $q=2$
as {\it annihilating\/} random walks,
as mixed annihilating-coalescing walks for $3 \leq q < \infty$
and in the limit $q \to \infty$ (with each site having its own unique
color at time zero) as coalescing random walks. Since the finite $q$
case can essentially be recovered from the $q=\infty$ model by projection,
one can restrict attention to the case of both forward and backward
coalescing random walks.

Naturally, the continuum scaling limit of
voter models is described by the BW. Indeed, in the voter model
as just described, it suffices to consider (as did Arratia~\cite{A81})
the collection of all outgoing BW paths from time zero. However,
if one modifies the voter model to allow for small noise, i.e., at
each space-time point there is a probability $p$ that rather than take
on the color of a neighboring spatial point one time step earlier,
a random color (out of $q$ possibilities, or a wholly new color
for $q=\infty$) is chosen (or nucleated),
then much more of the BW structure comes
into play in the scaling limit (in which also $p$ is properly scaled).
As analyzed in~\cite{FINR05}, this model in the scaling limit is one
in which new colors are nucleated on $(0,2)$ points of the BW
and it can be constructed by means of a Poissonian marking of those
points. The reason $(0,2)$ points are relevant is because a newly nucleated 
color in the voter model inside a cluster of some other color
creates two new boundaries which need
to persist for a macroscopic amount of time
before coalescing
in order to be seen in the scaling limit.

There are natural settings, namely the so-called
$q$-state stochastic Potts models of Statistical Physics,
such that for $q\geq 3$ (we recall that $q=2$ corresponds
to the Ising model) one needs to consider a more complex noise
structure in which
the probability of nucleation of new colors may depend on the color
of the site in question and its neighbors. For example, one may require
for nucleation that a site have a different color than its left
(respectively, right) neighbor. For that type of noise, it turns
out that the construction of the scaling limit naturally involves
the Poissonian marking of left (respectively, right) $(1,2)$ points.
The reason $(1,2)$ points are relevant here is that
the newly nucleated 
color in the voter model
is just to the right (or left)
of a previously
existing boundary
and creates a new boundary
that needs to persist
in the scaling limit.
This type of application of our marking of $(1,2)$
points
will be carried out 
in a future paper \cite{NRS08}.

\subsection{Outline of the Paper}
The remainder of the paper is organized as follows. In
Section~\ref{Brownian_Web}, we give a review of the basic structure of
the Brownian web and its dual (or backward) web,
with special emphasis on the $(1,2)$ points. In Section~\ref{marking_process},
we explain precisely how to mark $(1,2)$ points, which are points
where backward and forward BW paths touch, by first
defining for
finitely many backward and forward paths
a local time measure for touching to serve as a Poisson intensity
measure. The overall marking process
is then the limit as the number of forward and backward paths
tends to infinity. In Subsection~\ref{modify_the_Web}, we give
a preliminary
explanation of how the marking process will be used to construct
the BN and the DyBW.

In Section~\ref{modify_one_path}, we consider
the special marking process (and resulting
modified Brownian web path) constructed from a {\it single\/} forward
BW path and {\it all\/} backward paths that touch it from the right.
In particular,
we show that the resulting modified
forward path is related to the original
BW path by sticky reflection. Brownian motions
with a sticky interaction will also play an important role in later
sections as they do in \cite{SS07} and \cite{HW07}.
In Section~\ref{Brownian_Net}, we review
the construction from~\cite{SS07} of the BN and then prove that
our alternate construction using marked $(1,2)$ points is
equivalent.
In Section~\ref{dynamical_brownian_web},
we construct the DyBW and prove some elementary properties
of this object. 
Section \ref{proof} contains the proofs of many 
of the results
stated in previous sections along with
some propositions and lemmas that are needed for those proofs.
We note in particular
that Section \ref{sec-excursions}
contains a number of key results about the structure of excursions 
in the Brownian web from a single web path.

\section{The Brownian Web}
\label{Brownian_Web}

\subsection {The Forward Brownian Web}
\label{Forward_Web}
The (forward) Brownian web is the scaling limit of the discrete web
under diffusive space-time scaling; it is a random
collection of paths with specified starting points in space-time.
The (continuous) paths take values in a metric space $(\br^2,\rho)$
which is a compactification of $\r^2$.  $(\o,d)$ denotes the space
whose elements are paths with specific starting points. The metric
$d$ is defined as the
maximum of the sup norm of the distance between
two paths and the distance between their respective starting points.
The Brownian web takes
values in a metric space $(\h,d_\h)$,
whose elements are compact collection of paths in
$(\o,d)$ with $d_\h$ the induced Hausdorff metric. Thus the
Brownian web is an
 $(\h,\f_\h)$-valued random variable,
 where $\f_\h$ is the Borel $\s$-field  associated to the metric
 $d_\h$. The next theorem, taken from~\cite{FINR04}, gives some of the key
properties of the BW.

          \bteo
          \label{teo:char}
          There is an \( ({\cal H},{\cal F}_{{\cal H}}) \)-valued random variable
          \(
          {\W} \)
          whose distribution is uniquely determined by the following three
          properties.
          \begin{itemize}
                    \item[(o)]  from any deterministic point \( (x,t) \) in
          $\r^{2}$,
                    there is almost surely a unique path \( {B}_{(x,t)} \)
          starting
                    from \( (x,t) \).

                    \item[(i)]  for any deterministic, dense countable subset
        \(
          {\cal
                    D} \) of \( \r^{2} \), almost surely, \( {\W} \) is the
          closure in
                    \( ({\cal H}, d_{{\cal H}}) \) of \( \{ {B}_{(x,t)}: (x,t)\in
                    {\cal D} \}. \)
  \item[(ii)]  for any deterministic  $n$ and \((x_{1}, t_{1}),
        \ldots,
                    (x_{n}, t_{n}) \), the joint distribution of \(
                    {B}_{(x_{1},t_{1})}, \ldots, {B}_{(x_{n},t_{n})} \) is that
                    of coalescing Brownian motions (with unit diffusion
        constant).

          \end{itemize}
          \eteo

Note that (i) provides a practical construction of the Brownian web. For $\cal D$ as defined above,  construct coalescing Brownian motion paths starting from $\cal D$. This defines a {\it skeleton} for the Brownian web. $\W$ is simply defined as the closure of this precompact set of paths.

\subsection{The Backward (Dual) Brownian Web}
\label{backwardBW}
We have considered in Subsection \ref{arrows} the backward discrete web as the set of all
coalescing random walks starting from $\Z_{odd}^2$ running backward in
time without crossing the forward discrete web paths. The backward (dual)
BW $\hat{\mathcal{W}}$ may be defined analogously as a
functional of the (forward) BW $\mathcal{W}$.
More precisely for a countable dense deterministic set of space-time points,
the backward BW path from each of these is the (almost surely) unique
continuous curve (going backwards in time) from that point that
does not cross (but may touch) any of the (forward) BW paths;
$\hat{\mathcal{W}}$
is then the closure of that collection of paths. The
first part of the next proposition
states that the ``double BW'', i.e., the pair
$(\mathcal{W},\hat{\mathcal{W}})$, is the diffusive
scaling limit of the corresponding discrete pair
$(W^\delta,\hat{W}^\delta)$ (as the scale parameter $\delta \to 0$).
Convergence in the sense of weak convergence of probability
measures on $(\h,\f_\h) \times (\hat \h,\hat \f_\h)$ was proved in~\cite{FINR04};
convergence of finite dimensional distributions and
the second part of the proposition were already contained
in~\cite{TW98}.

\bprop
\label{forwardandbackward1}
\begin{enumerate}
\item Invariance principle : $(W^\delta,\hat{W}^\delta) \to
(\mathcal{W},\hat{\mathcal{W}})$ as $\delta \to 0$.
\item For any (deterministic)
pair of points $(x,t)$ and $(\hat x, \hat t)$ there is
almost surely a unique forward path $B$ starting from $(x,t)$ and a
unique backward path $\hat B$ starting from $(\hat x,\hat t)$.
\end{enumerate}

\eprop

The next proposition, from~\cite{STW00}, which gives the
joint distribution of a single forward and single backward BW path,
has an extension to the joint distribution of finitely many
forward and backward paths. We remark that that extension
can be used to give a characterization (or construction)
of the double Brownian web $(\W,\hat \W)$ analogous to the one for
the (forward) BW from
Theorem~\ref{teo:char} --- see~\cite{STW00, FINR05} for more details.
%using independent forward and backward Brownian motions starting
%from  a countable deterministic dense set of space-time
%points --- see~\cite{STW00, FINR05} for more details.

\bprop
\label{forwardandbackward2}
\begin{enumerate}
\item \underline{Distribution of $(B,\hat B)$}: Let $(B_{ind},\hat B_{ind})$
be a pair of independent forward and backward Brownian motions starting
at $(x,t)$ and $(\hat x,\hat t)$ and let $(R_{\hat B_{ind}}(B_{ind}),\hat B_{ind})$ be
the pair obtained after reflecting (in the Skorohod sense) $B_{ind}$ on
$\hat B_{ind}$, i.e., $R_{\hat B_{ind}}(B_{ind})$ 
is the following function of $u\in[t,\hat t]$:
 \begin{equation}
R_{\hat B_{ind}}(B_{ind})=\left\{ \begin{array}{ll}
B_{ind}(u)-0 \wedge \min_{t \leq v \leq u } ( B_{ind}(v)-\hat B_{ind}(v))\
\ \ \ \textrm{on $\{B_{ind}(t) \geq \hat B_{ind}(t) \}$}, \\
B_{ind}(u)-0 \vee \max_{t \leq v \leq u } ( B_{ind}(v)-\hat B_{ind}(v))\
\ \ \ \textrm{on $\{B_{ind}(t) < \hat B_{ind}(t) \}$}. \end{array}
\right.
\end{equation}
Then
\begin{equation}
(R_{\hat B_{ind}}(B_{ind}),\hat B_{ind})=(B,\hat B) \ \ \textrm{in law},
\end{equation}
where $B$ is the path in $\W$ starting at $(x,t)$ and $\hat B$ is 
the path in $\hat W$ starting at $(\hat x,\hat t)$.
\item Similarly,
\begin{equation}
(B_{ind},R_{B_{ind}}(\hat B_{ind}))=(B,\hat B) \ \ \textrm{in law}.
\end{equation}
\end{enumerate}

\eprop

\subsection{$(1,2)$ Points of the Brownian Web.}
\label{special_points}

While there is only a single path from any deterministic point in
$\R^2$ in both the forward and backward webs, there exist random
points $z \in \R^2$ with more than one path
passing through or starting from~$z$.

        We now describe the ``types'' of points $(x,t)\in\r^2$,
        whether deterministic or not. 
	We say that two paths $B,B'\in\W$ are equivalent paths entering 
	$z=(x,t)$, denoted by $B=_{in}^z B'$, iff $B=B'$ on $[t-\e,t]$
	for some $\e>0$. The relation $=^z_{in}$ is a.s. an equivalence relation on the set
        of paths in $\W$ entering the point $z$ 
	and we define $m_{in}(z)$ as the number of those equivalence 
	classes. ($m_{in}(z)=0$ if there are no paths entering $z$.)
        $m_{out}(z)$ is defined as the number of distinct paths starting from $z$.
	For $\hat \W$, $\hat m_{in}(z)$ and $\hat m_{out}(z)$ are defined similarly.

          \bdf
          The type of $z$ is the pair $(m_{in}(z),
          m_{out}(z))$.
          \edf
 %%%%%%%%%%%%%%%%%%%%%%%%%%%%%%%%%%%%%%%%%%%%%%%%
          %%%%%%%%%%%%%%%%%%  FIGURE 3 %%%%%%%%%%%%%%%%%%%%%%%%
          %%%%%%%%%%%%%%%%%%%%%%%%%%%%%%%%%%%%%%%%%%%%%%%%

	\begin{figure}[!ht]
         \begin{center}
          \includegraphics[width=6cm]{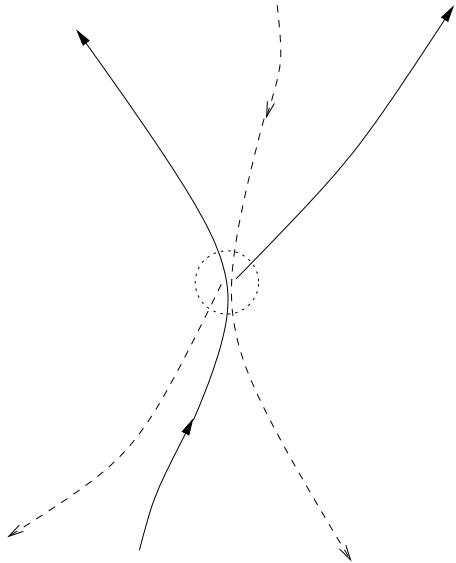}
          \caption{A schematic diagram of a left $(m_{\mbox{{\scriptsize\emph{in}}}},m_{\mbox{{\scriptsize
          \emph{out}}}})=(1,2)$ point 
          with necessarily also
          $(\hat m_{\mbox{{\scriptsize\emph{in}}}},\hat m_{\mbox{{\scriptsize
          \emph{out}}}})=(1,2)$.
          In this example the incoming forward path connects to the leftmost
          outgoing path
          (with a corresponding dual connectivity for the backward paths), the right 
          outgoing path is a newly born path.
      	  }\label{fig12}
          \label{12} %don't know what this does
          \end{center}
          \end{figure}

          The following results from \cite{TW98} (see also ~\cite{FINR05}) specify what types
          of points are possible in the Brownian web.

 \bteo
          \label{teo:types}
          For the Brownian web, almost surely, every $(x,t)$ has one of
          the
          following types, all of which occur: $(0,1)$, $(0,2)$, $(0,3)$,
          $(1,1)$,
          $(1,2)$, $(2,1)$.
          \eteo

\bprop
          \label{prop:dual}
          For the Brownian web, almost surely for {\em every}
          $z$ in $\r^2$, $\hat m_{in}
          (z)=m_{out}(z)-1$
          and
          $\hat m_{out}(z)=
          m_{in}(z)+1$. See Figure~\ref{fig12}.
          \eprop

%          \brm
It is important to realize that points of type $(1,2)$
can be characterized
in two ways, both of which will play a crucial role
in our construction of the DyBW and BN. 1)  $z\in \r^2$ is of
type $(1,2)$ precisely if both a forward and a backward path pass
through
$z$. 2) A single incident path continues along exactly
one
of the two outward paths --- with the choice determined intrinsically.
%See Figure~\ref{fig12} for a schematic
%diagram of a ``left-handed'' continuation.
It is either left-handed or
right-handed according to
whether the continuing
path is to the left or the right of the
incoming (from later time) backward
path. For a left $(1,2)$ point $z$, the right (resp, left) outgoing path
will be referred to as the {\it newly born path} starting from $z$.
See Figure~\ref{fig12} for a schematic
diagram of the ``left-handed'' case.
Both varieties occur and it is
known ~\cite{FINR05} that each of the two varieties, 
as a subset of $\R^2$, has Hausdorff dimension~$1$.
As noted in Section~\ref{intro}, the two varieties of $(1,2)$ points
play the same role in the continuum that left and right arrows play
in the discrete setting. 
In particular, one can change the direction of the ``continuum'' 
arrow at a given $(1,2)$
point $z$ by simply
connecting the incoming path to the newly born path starting from $z$. In the discrete 
picture, this amounts to changing the direction of an arrow
whose switching induces a ``macroscopic'' effect in the web.

%\section{The Marking Process on the Brownian Web}
\section{Marked $(1,2)$ Points on the Brownian Web.}
\label{marking_process}

\subsection{The Local Time Measure}
\label{localtimesec}
 Recall that the $\phi$-Hausdorff outer measure of an
 arbitrary subset $E$ of $\R$ for $\phi:(0,\infty)\raw (0,\infty)$ is defined as 
\begin{equation}
 m_\phi(E)=\lim_{\delta\downarrow 0} \ \inf\{\sum \phi(|b_i-a_i|) \ | \ E \subset \bigcup_i [a_i,b_i], |b_i-a_i|<\delta \}.
\end{equation}
In the following, we set $\phi(t)=\sqrt{2t\log(|\log(t)|)}$ and we denote the Lebesgue measure of $E$ by $|E|$.
Restricted to Borel subsets $E$ of $\R$, $m_\phi$ is a measure.

\bprop
\label{localtime}
\begin{enumerate}
 \item Let $(B,\hat B)$ be defined as in Proposition \ref{forwardandbackward2}.
 For almost every realization of $\W$, for every $t \leq u \leq \hat t$
\begin{equation}
\lim_{\epsilon \downarrow 0} \frac{1}{2\epsilon} |\{v: t\leq v\leq u, \frac{|B(v)-\hat B(v)|}{\sqrt{2}} \leq \epsilon\}|
\end{equation}
exists and will be denoted by $L_{B,\hat B}(u)$.
\item For a Borel set $A\subset \R$
\begin{equation}
\label{haus::exp}
\int_{u\in A} d L_{B,\hat B}(u) \ = \  m_{\phi}\left(\{ u \in A \  | \ B(u)=\hat B(u)\}\right) 
\end{equation}
\item \underline{Distribution of $L_{B,\hat B}$}:  $L_{B,\hat B}$ is a stochastic process on $[t,\hat t]$ which is identical in law to $\bar L_{B,\hat B}$ defined as follows:
\begin{equation}
 \bar L_{B,\hat B}(u)=\left\{ \begin{array}{ll}
               -\  \ 0 \wedge \min_{t \leq v \leq u } ( B_{ind}(v)-\hat B_{ind}(v))/\sqrt{2}\ \ \ \ \textrm{on $\{B_{ind}(t) \geq \hat B_{ind}(t) \}$},\\
             \  \  0 \vee \max_{t \leq v \leq u } ( B_{ind}(v)-\hat B_{ind}(v))/\sqrt{2}\ \ \ \ \textrm{on $\{B_{ind}(t) < \hat B_{ind}(t) \}$}, \end{array}
	      \right.
\end{equation}
where $(B_{ind},\hat B_{ind})$ are defined as in Proposition \ref{forwardandbackward2}.

\end{enumerate}

\eprop

Note that the third statement is analogous to the famous property discovered
by
L\'{e}vy that the local time (at the origin) of a one-dimensional
Brownian motion is identical in law with its record time process
(see, e.g., \cite{KS91}). Statement 2 is analogous to the fact that the 
measure induced by the local time at $0$ of a standard Brownian
motion coincides with the $\phi$-Hausdorff measure of its zero-set (see Theorem 1 in \cite{P81}).

Let  us consider a family of $n$ forward paths $\{B_{i}\}_{i=0}^{n-1}$ and a family of $m$ backward
paths $\{\hat B_{j}\}_{j=0}^{m-1}$. We will generally choose theses paths so that $B_i$ and $\hat B_i$ have the same starting point $z_i$ with $\mathcal{D}=\{z_i\}_{i=0}^\infty$ some dense deterministic set of points in $\R^2$ as defined in Subsection \ref{Forward_Web}; also for consistency with other notation, we will generally assume that $z_0$ is the origin in $\R^2$. In non-ambiguous contexts, $\{B_{i}\}_{i=0}^{n-1}$ and $\{\hat B_{j}\}_{j=0}^{m-1}$ will also refer to the union
of their respective traces in $\R^2$.

%\bdf
%$\W_n$ and $\W_m$ intersect at $(x,t)$ iff there exists $B_i\in\W_n$ and $\hat B_j\in\hat \W_m$ s.t. $(x,t)\in \W_n\bigcap \hat \W_m$.
%\edf
The expression for $L_{B,\hat B}$ given in (\ref{haus::exp}) can be easily generalized to the family $\{B_{i}\}_{i=0}^{n-1}$ and $\{\hat B_{j}\}_{j=0}^{m-1}$.
E.g., for a Borel $A\subset\R$, we simply define $L_{n,m}(A)$ by
\begin{eqnarray*}
\int_{u\in A} d L_{n,m}(u) \ & = & \  m_{\phi}\left(\{ t \in A \ | \ \exists x\in \R\ \  \textrm{s.t.} \ \ (x,t)\in \{B_i\}_{i=0}^{n-1} \cap \{\hat B_j\}_{j=0}^{m-1}\}\right) \\
& = & 
m_\phi\left( A\cap \mathcal{P}(\{B_i\}_{i=0}^{n-1} \cap \{\hat B_j\}_{j=0}^{m-1})\right),
\end{eqnarray*}
where $\mathcal{P}$ denotes the projection onto the t-axis.

Finally, we can extend $L_{n,m}$ to be a
measure acting on $\R^2$ in the following way, which implicitly
uses the a.s. property of $\W$ that if a forward and a backward family meet at some $t$, they do so only
at a single value of $x$.
\bdf
\label{localtimemeasure}{\bf [Local time measure]}

For the forward family $\{B_{i}\}_{i=0}^{n-1}$ and the backward family $\{\hat B_{j}\}_{j=0}^{m-1}$, we define the local time (outer) measure $\mathcal{L}_{n,m}$ on $\R^2$ as follows. For a 
general \textrm{space-time domain} $O$,
\begin{equation}
\mathcal{L}_{n,m}(O)=m_\phi\left(\mathcal{P}(\{B_i\}_{i=0}^{n-1} \cap \{\hat B_j\}_{j=0}^{m-1} \cap O)\right).
\end{equation}
In particular, $\mathcal{L}_{n,m}$ is supported on the space-time points where the
forward family touches the backward family. Finally, we define
an outer measure
\begin{equation}
\label{lo-def}
\mathcal{L}(O)=m_\phi\left(\mathcal{P}(\{B_i\}_{i=0}^{\infty} \cap \{\hat B_j\}_{j=0}^{\infty} \cap O) \right).
\end{equation}
$\Lo(O)$ will be referred to as the local time outer measure of $O$.
\edf
Both $\mathcal{L}_{n,m}$ and $\mathcal{L}$
are measures when restricted to Borel sets
but may take the value $\infty$.
We note that for any open set $O\subset\R^2$, $\Lo(O)=\infty$. However, we will later encounter (see e.g., Subsection \ref{rrrrrrrr}) some very natural subsets $O\subset\R^2$ with finite $\Lo$-measure.
See Section 4i of \cite{Ts04} for a similar discussion.

\subsection{The Marking Process}
\label{the_marking_process}

Let us consider the Poisson point process on $\mathbb{R}^2\times\mathbb{R}^+$ with intensity measure
$$I_{n,m}(O\times[0,\tau])=\sqrt{2}\mathcal{L}_{n,m}(O)\cdot\tau,$$ where $O$ is any open subset of $\mathbb{R}^2$. We define
the {\it partial marking} process $\tau\to\mathcal{M}_{n,m}(\tau)$ as
\begin{equation}
 \mathcal{M}_{n,m}(\tau)=\{z\in\R^2 : (z,u) \ \ \textrm{is a Poisson point for some} \ \ \ u\leq\tau\}.
\end{equation}

Heuristically, $\mathcal{M}_{n,m}(\tau)$
consists of the locations
of the switching (in the DyBW) between dynamical times $0$ and $\tau$
if one restricts the
dynamics to
the ``arrows'' at the intersection of the
forward family $\{B_i\}_{i=0}^{n-1}$
and the backward family $\{\hat B_j\}_{j=1}^{m-1}$,
while other arrows remain frozen.
In order to introduce a ``full dynamics"
we will couple the sequences $\{\mathcal{M}_{n,m}(\tau)\}_{n,m}$
in such way that for $n'\geq n$ and $m'\geq m$,
$\mathcal{M}_{n,m}(\tau)\subseteq \mathcal{M}_{n',m'}(\tau)$.
To achieve this, we define the point process $\mathbb{M}$ as follows:

\bdf
$\mathbb{M}$ is the four-dimensional Poisson point process  on
$\mathbb{R}^2\times \mathbb{R}^{+}\times \mathbb{N}\times\mathbb{N}$
with (locally finite and random) intensity measure $I$ defined  by
\begin{equation}
 I(O\times[0,\tau]\times\{0,...,n-1\}\times\{0,...,m-1\})= \sqrt{2} \mathcal{L}_{n,m}(O) \cdot \tau,
\end{equation}
where $O$ is any open subset of $\mathbb{R}^2$.
\edf

We can then define
$\mathcal{M}(\tau)$ as
\begin{equation}
\label{m_tau}
 \mathcal{M}(\tau)=\{z: \ \ \textrm{$(z,s,n',m')$  is in $\mathbb{M}$ for some $n',m'$ and some $s\leq \tau$}\}.
\end{equation}
and $\mathcal{M}_{n,m}(\tau)$ is simply obtained by adding the restriction to (\ref{m_tau}) that $n'\leq n-1$ and $m'\leq m-1$.
 
\bigskip

Informally, $\{\M(\tau)\}_{\tau\geq0}$ can be seen as a Poisson Point Process on $\R^2\times\R$ with intensity measure
$ \sqrt{2} \Lo(dz) \times d \tau$. In particular, for a Borel $O\subset\R^2$ with $\Lo(O)<\infty$, 
$\M(\tau)\cap O$ is a Poisson point Process on $\R^2\times\R$ with intensity measure $\sqrt{2} \ 1_{z\in O} \ \Lo(dz) \times d \tau$.

\subsection{Modifying the Web Using Marking}
\label{modify_the_Web}

\subsubsection{Constructing the Brownian Net}
\label{con::net}
Let $\tau>0$. We define a partial Brownian net
$\mathcal{N}_{n,m}(\tau)$ by
having branching  at the points of the partial
marking $\mathcal{M}_{n,m}(\tau)$ .(Later we will write 
$\Net_n(\tau)$ for $\Net_{n,n}(\tau)$.)
For example, if the $(1,2)$ point in Figure~\ref{fig12}
is marked, then the Brownian net will include not only paths that
connect to the left outgoing path (as in the original web)
but also ones that connect to the right outgoing path.
More formally, the set of paths in $\mathcal{N}_{n,m}(\tau)$ starting from $z\in\R^2$ is the set of
paths interpolating the set $\mathcal{S}$ of points $\M_{n,m}(\tau)\cup\{z\}\cup\{+\infty\}$
with paths in $\W$---i.e., between any {\it consecutive} pair of points in $\pi\cap\mathcal{S}$,
$\pi$ coincides with a path in $\W$. 
%concatenating
%finitely many portions of paths of the Brownian web at points in
%$\mathcal{M}_{n,m}(\tau)$, i.e., the set of paths $\pi$ of the form
%\begin{equation}
%\pi= \{B(t)\ \ \textrm{on}  \ \ [t_0,t_{mark,1})\}\bigcup \ \
%\left(\bigcup_{i=1}^{m-1} \{B_{mark,i}(t) \ \textrm{on}  \ \
%[t_{mark,i},t_{mark,i+1}) \} \right)
%\end{equation}
%where $t_0 \in {\mathbb{R}}$, $B \in \mathcal{W}$ starts at time
%$t_0$, $m \in {\mathbb{N}}$ and
%\begin{enumerate}
%\item $t_0<t_{mark,1}<t_{mark,2}<...<t_{mark,m}=\infty$ ,
%$t_0$ is the starting time of $B\in\mathcal{W}$.
%\item $B_{mark,i}$ is either of the two outgoing paths in
%$\mathcal{W}$ starting from the marked point
%$(x_{mark,i},t_{mark,i})$ in the partial marking $\mathcal{M}_{n,m}(\tau)$,
%\item $B_{mark,i}(t_{mark,i+1})=x_{mark,i+1}$ and
%$B(t_{mark,1})=x_{mark,1}$ (i.e., $\pi$ is continuous).
%\end{enumerate}

Finally, we define $\mathcal{N}_{mark}(\tau)$ as the closure of
$\bigcup_{n,m=1}^\infty \mathcal{N}_{n,m}(\tau)$. 
In other words,
$\mathcal{N}_{mark}(\tau)$
is defined by allowing branching
at every marked $(1,2)$ points in the Brownian
web $\mathcal{W}$. Analogously, we can define a backward partial Brownian Net
$\hat{\mathcal{N}}_{n,m}(\tau)$ by allowing branching at the points $\mathcal{M}_{n,m}(\tau)$
in the dual Web $\hat{\mathcal{W}}$
and define $\hat{\mathcal{N}}_{mark}(\tau)$ as the closure of
$\bigcup_{n,m=1}^\infty \hat{\mathcal{N}}_{n,m}(\tau)$.
In Section~\ref{Brownian_Net}, we prove the
equivalence of $\mathcal{N}_{mark}(\tau)$ to the Brownian net construction
of Sun and Swart~\cite{SS07}, which by their results (see Theorem 1.1 in \cite{SS07}) then implies
convergence of the properly rescaled discrete net to $\mathcal{N}_{mark}(\tau)$
in an appropriate topology.

\subsubsection{Constructing the Dynamical Brownian Web}
\label{modify_the_Web2}

We can construct a partial
dynamical Brownian web $\W_{n,m}(\tau)$, at dynamical time $\tau$,
to replace the original $\W$ by
switching the direction of all the
marked $(1,2)$ points in $\mathcal{M}_{n,m}(\tau)$. Formally,
$\pi$ is in $\W_{n,m}(\tau)$ iff $\pi$ is in the the partial net $\mathcal{N}_{n,m}(\tau)$ and
at
each time $t=\bar t_i$ that $\pi$ hits a point $(\bar x_i,\bar t_i)\in\mathcal{M}_{n,m}(\tau)$,
it then follows $B_{new}^i$, the newly born path of $\W$ starting from $(\bar x_i,\bar t_i)$, on
$[\bar t_i,\bar t_i+a]$ for some $a>0$. 
A nontrivial question is the existence of a limit
for $\W_{n,m}(\tau)$
as $n,m \to \infty$. It will be shown in
Section~\ref{dynamical_brownian_web}
that for almost all realizations of
the web and its marking, a limit $\W(\tau)$ exists for every $\tau$ (see Proposition \ref{existence_of_Br}).

\section{Sticky Brownian Motion by Marking a Single Path}
\label{modify_one_path}
From here through Section \ref{Dynamical::BW}, $\tau$ will denote a 
fixed deterministic number and the marking will refer to $\mathcal{M}(\tau)$.

We first recall the definition of a one-dimensional sticky (at the origin) Brownian motion.
\bdf\label{def-1d-sbm}
$B_{stick,x}$ is a $(1/\bar \tau)$-sticky Brownian motion starting at $x$ iff there exists a one-dimensional standard Brownian motion $B$ s.t.
\begin{equation}
\label{equation_of_sticky}
\forall t\geq 0, \ \  d B_{stick,x}(t)=1_{B_{stick,x}(t)\neq 0} d  B(t) + {\bar \tau} 1_{ B_{stick,x}(t)=0} dt.
\end{equation}
and $B$ is constrained to stay positive as soon it first hits zero.
\edf
It is known that (\ref{equation_of_sticky}) has a unique (weak) solution. 
Furthermore,
for $x=0$
this solution can be constructed
from  
a time-changed reflected Brownian motion. More precisely, consider 
$$
t\rightsquigarrow|\bar B|(C(t)), \ \ \ \txt{with $C^{-1}(t)=t+\fr{1}{\bar \tau} L_0(t)$},
$$  
where $|\bar B|$ is a reflected Brownian motion and $L_0$
is
its local time at the origin. Then there exists a Brownian motion $B$ such that $(|\bar B|(C(\cdot)),B)$ is a solution of (\ref{equation_of_sticky}) (see, e.g., \cite{W}). In words, the sticky Brownian motion is obtained from the reflected one by ``transforming'' local time into real time. In particular,
it spends a positive Lebesgue measure of time at the origin
and the larger the ``degree of stickiness'' $1/\bar\tau$ is, the more
the path sticks to the origin.

In this section
we consider the path $\rz$ starting at $z\in\mathcal{D}$ and constructed
by  switching only the direction of the {\it left} $(1,2)$ points in $\mathcal{M}(\tau)$ on 
$B_0$, the path of $\mathcal{W}$ starting from the origin. As we shall see, unlike in the complete DyBW,
it is not difficult to construct $\rz$ and the law of the pair $(\rz,B_0)$ can be characterized explicitly. In particular, if we set $\ro\equiv\rz$ for $z=(0,0)$ then it readily follows from Proposition \ref{cornerstone} below that $(\ro-B_0)/\sqrt{2}$ is a $(\sqrt{2}/\tau)$-sticky Brownian motion.
This will be very useful in the rest of the paper (see Sections ~\ref{Brownian_Net} and ~\ref{dynamical_brownian_web}) where 
the analysis
of paths that result from switching left and right 
$(1,2)$ points is a direct
extension of the analysis here. Our
construction of a sticky Brownian motion
using the marked excursions defined next is 
similar to Warren's construction in \cite{W02}
using the excursions of a single Brownian motion.

\bigskip

\bdf\label{df-excursion}{\bf [Excursions]}
Let $B_{new}$ be the newly born path emerging from a $(1,2)$ point $z=(x,t)$ on any path $B\in\W$.
The segment of $B_{new}$ before coalescence with $B$ is called an excursion from $B$.

$D(e)$ is the time duration of the excursion $e$, $|e|\equiv \sup\{|B-e|(s):t\leq s\leq t+D(e)\}$ is  its diameter, $T(e)\equiv t$ its starting time, $(T(e),T(e)+D(e))$ its lifespan.

If an excursion $e$ starts from a marked point, $e$ is called a marked excursion. 

A right marked excursion $e$ is called nested iff there exists another right marked excursion $e'$
s.t. $T(e)$ belongs to the lifespan of $e'$. An analogous definition holds for left marked excursions.

If a marked excursion $e$ is not nested, $e$ is said to be a maximal excursion. 
\edf

${_{[1]}r_0}$ may be defined as the path obtained after joining together all the right maximal excursions from $B_0$.
Stated differently, $\ro$ is the path whose excursions  (in the standard sense) from $B_0$ coincide with the right maximal excursions from $B_0$ in the marked Brownian web.  
We note that everytime $\ro$ hits a left $(1,2)$ point on $B_0$ it then follows the newly born path starting from it. 
(Among all the marked left $(1,2)$ points $\ro$ only hits the starting points of maximal excursions since nested excursions are ``straddled'' by some maximal excursions). Thus $\ro$ is consistent with the informal definition
in terms of switching
given earlier in this section.

Next, we recall that
for any deterministic point $z\in\R^2$, $B_z\in\W$ is the path starting from $z$. We define ${_{[1]}r_{z}}$ as the path starting from $z$ obtained by switching all the left marked $(1,2)$ points on $B_0\cap B_z$.
(This informal definition may be made precise as was done for $\ro$ by using the right 
maximal excursions from $B_{z'}$, where $z'$ is the coalescing
point between $B_0$ and $B_z$.)
Note that ${_{[1]}r_z}$ is a continuous path. To prove that,
it is clearly enough to show that
for fixed $T,\epsilon\in(0,\infty)$ the process ${_{[1]}r_z}$ only performs finitely many excursions of diameter 
$\geq\epsilon$ away from $B_0$ on the interval $[0,T]$. If that were not the case, there would exist a sequence 
of marked excursions $\{e_{k}\}$  from $B_0$ such that $e_{k}$ would make an excursion away from $B_0$ with 
diameter $\epsilon$ and duration $t_k$, with $t_k\rightarrow 0$. 
But that would violate the compactness of $\W$. 

\bigskip

We now set up some notation. For a path $\pi$ in $(\Pi,d)$ starting from $z$, we denote by $t_\pi$, the starting time of $\pi$. For two paths $\pi_1,\pi_2$, $T_{\pi_1,\pi_2}\equiv \inf\{t>t_{\pi_1}\vee t_{\pi_2}\ : \ \pi_1(t)=\pi_2(t)\}$ denotes the first meeting time of $\pi_1$ and $\pi_2$, which may be
$+\infty$.
In Subsection \ref{proof:sticky} we show the following proposition.

\bprop
\label{cornerstone}
For any deterministic $z\in\R^2$, almost surely, there exists $B_z^{(1)}$, a standard Brownian motion starting at $z$ so that $_{[1]} r_z$ satisfies the following SDE.
\begin{eqnarray}
\label{target::SDE}
d _{[1]} r_z(t)& = & d B_z^{(1)}(t) + 1_{_{[1]}  r_z(t)= B_0(t)} \ \tau \   dt, \nonumber \\[0.1in]
d B_{0}(t) \ d B_z^{(1)}(t) \ & = & \ 1_{\rz(t)=B_0(t)} dt, \nonumber \\[0.1in]
\forall t\geq T_{\rz,B_0}, & & \ \ \ \ \rz(t) \ \geq \ B_0(t)  \label{equation_stick_BM}
\end{eqnarray}
\eprop
\bigskip

Here  $d B_{0}(t) d B_z^{(1)}(t)$ denotes $d \langle B_0,B_z^{(1)}\rangle(t)$, where $\langle B_0,B_z^{(1)}\rangle(t)$ is the cross-variation process of $B_0$ and $B_z^{(1)}$ at time $t$.
The second part of Equation (\ref{equation_stick_BM})
amounts to saying that away from the diagonal $\{t:\rz(t)=B(t)\}$, $B_0$ and $\rz$ evolve independently  while on the diagonal they are perfectly correlated. In particular, without the drift on the diagonal to ``unstick'' $\rz$ from $B_0$, $\rz$ and $B_0$ would coalesce rather than stick when they meet.

Adopting the usual terminology, we will say that $_{[1]} r_z$ is distributed as a Brownian motion stickily reflected off $B_0$ with a degree of stickiness $1/\tau$. In particular, for $z=(x,0)$
the process $\{(\rz -B_0)/\sqrt{2}\}$ is a $(\sqrt{2}/\tau)$-sticky Brownian motion 
(see Definition \ref{def-1d-sbm}).

In \cite{SS07}, Sun and Swart studied a similar 
equation but with the difference that $\ro$ (resp., $B_0$) is replaced by a right (resp., left) drifting Brownian motion (see Equation (\ref{lrsde})). For that equation, they established existence and uniqueness of a weak solution (see Proposition 2.1 in \cite{SS07}). Since (\ref{target::SDE}) and (\ref{lrsde}) only differ by their drift terms, existence and uniqueness for (\ref{target::SDE}) follows from their result and the  Girsanov Theorem.
In particular, a (weak) solution $(\ro,B_0)$ of (\ref{target::SDE}) is a strong Markov process.

\section{The Brownian Net by Marking}
\label{Brownian_Net}

In Subsection \ref{remainder}, we outline the construction of the Brownian
net given by Sun and Swart \cite{SS07} and  state some related results. 
The presentation we give of that construction is taken from~\cite{SSS08}.
As  will be seen, this construction of the Brownian net is different in
spirit to the one using marking given in Subsection \ref{con::net}.
However, we will show in Theorem \ref{equivalence_Net} that the two 
constructions lead to the same object.

\subsection{The Brownian Net as Introduced by Sun and Swart}
\label{remainder}
We now recall the {\em left-right Brownian web} $(\W_l, \W_r)$, which is
the key intermediate object in the construction of the Brownian net in
\cite{SS07}. Following \cite{SS07}, we call $(l_1, \ldots, l_m; r_1,
\ldots, r_n)$ a collection of {\em left-right coalescing Brownian
motions}, if $(l_1, \ldots, l_m)$ is distributed as coalescing
Brownian motions each with drift $-\tau$, $(r_1, \ldots, r_n)$ is
distributed as coalescing Brownian motions each with drift $+\tau$, paths
in $(l_1, \ldots, l_m; r_1, \ldots, r_n)$ evolve independently when
they are apart, and the interaction between $l_i$ and $r_j$ when they
meet is a form of sticky reflection. More precisely, for any $L \in
\{l_1, \ldots, l_m\}$ and $R\in\{r_1, \ldots, r_n\}$, the joint law of
$(L,R)$ at times $t>t_L\vee t_R$ is characterized as the unique
weak solution of

\begin{eqnarray}
d L(t)\ & = &  \ d B_l - \tau dt, \nonumber \\[0.1in]
d R(t)\ & = &\ d B_r +\tau d t, \nonumber \\[0.1in]
d \langle B_l,B_r\rangle(t)\ & = & \ 1_{L(t)=R(t)} \ dt, \nonumber \\[0.1in]
\forall t\ \geq \ T_{R,L}, & & R(t) \ \geq \ L(t), \label{lrsde}
\end{eqnarray}
where $B^{l}, B^{r}$ are two standard
Brownian motions. We then have the following characterization of the left-right
Brownian web from \cite{SS07}.

\label{T:lrwebchar}{\bf [Characterization of the left-right Brownian web]}\\
There exists an $(\mathcal{H}^2, \mathcal{F}_{\mathcal{H}^2})$-valued random variable $(\W_l,
\W_r)$, called the standard left-right Brownian web (with parameter $\tau>0$), whose distribution
is uniquely determined by the following two properties:
\begin{itemize}
\item[{\rm (a)}] $\W_l$, resp.\ $\W_r$, is distributed as the standard
Brownian web, except tilted with drift $-\tau$, resp.~$+\tau$.

\item[{\rm (b)}] For any finite deterministic set $z_1, \ldots, z_m,
z'_1, \ldots, z'_n\in\R^2$, the subset of paths in $\W_l$ starting from
$z_1, \ldots, z_m$, and the subset of paths in $\W_r$ starting from
$z'_1, \ldots, z'_n$, are jointly distributed as a collection of
left-right coalescing Brownian motions.
\end{itemize}

Similar to the Brownian web, the left-right Brownian web $(\W_l, \W_r)$
admits a natural dual $(\hat\W_l, \hat\W_l)$ which is equidistributed
with $(\W_l, \W_r)$ modulo a rotation by
$180^{\rm o}$ of $\R^2$. In particular, $(\W_l,\hat\W_l)$ and $(\W_r, \hat\W_r)$ are
pairs of tilted double Brownian webs.

\bigskip

Based on the left-right Brownian web, \cite{SS07} gave three
equivalent characterizations of the Brownian net, which are called
respectively the {\em hopping, wedge, and mesh characterizations}. We
first recall what is meant by hopping, wedges and
meshes.\medskip

{\bf Hopping:} Given two paths $\pi_1, \pi_2\in\Pi$, let $t_1$ and $t_2$ be the starting times of those paths.
For $t >
t_1 \vee t_2$ (note the strict inequality),
 $t$ is called
an intersection time of $\pi_1$ and $\pi_2$ if $\pi_1(t)=\pi_2(t)$. By
hopping from $\pi_1$ to $\pi_2$, we mean the construction of a new
path by concatenating together the piece of $\pi_1$ before and the
piece of $\pi_2$ after an intersection time. Given the left-right
Brownian web $(\W_l, \W_r)$, let $H(\W_l\cup\W_r)$ denote the set of paths
constructed by hopping a finite number of times between paths in
$\Wl\bigcup\Wr$.\medskip

\noindent
{\bf Wedges:} Let $(\hat\W_l, \hat\W_r)$ be the dual left-right Brownian
web almost surely determined by $(\W_l, \W_r)$. For a path $\hat\pi \in
\hat\Pi$, let $t_{\hat \pi}$ denote its (backward) starting
time. Any pair $\hat l\in\hat\W_l$, $\hat r\in \hat\W_r$ with 
$\hat r(t_{\hat l}\wedge t_{\hat r}) < \hat
l(t_{\hat l}\wedge t_{\hat r})$ defines an open set
\be\label{wedge}
W(\hat r, \hat l) = \{(x,u)\in\R^2 : T<u< t_{\hat l}\wedge
t_{\hat r},\ \hat r(u) < x < \hat l(u) \},
\ee
where $T = \sup\{t<t_{\hat l}\wedge t_{\hat r}: \hat
r(t) =\hat l(t)\}$ is the first (backward) hitting time of $\hat r$
and $\hat l$, which might be $-\infty$. Such an open set is called a
{\em wedge} of $(\hat\Wl, \hat\Wr)$.\medskip

\noindent
{\bf Meshes:}  By definition, a {\em mesh} of $(\Wl,\Wr)$ is an open set of
the form
\be\label{mesh}
M=M(r,l)=\{(x,t)\in\R^2:t_l<t<T_{l,r},\ r(t)<x<l(t)\},
\ee
where $l\in\Wl$, $r\in\Wr$ are paths such that $t_l=t_r$,
$l(t_l)=r(t_r)$ and $r(s)<l(s)$ on
$(t_l,t_l+\epsilon)$ for some $\e>0$.  We call
$(l(t_l),t_l)$ the bottom point, $t_l$ the bottom time,
$(l(T_{l,r}),T_{l,r})$ the top point, $T_{l,r}$ the top
time, $r$ the left boundary, and $l$ the right boundary of $M$.\medskip

Given an open set $A\subset \R^2$ and a path $\pi\in\Pi$, we say $\pi$
{\em enters} $A$ if there exist $t_\pi <s<t$ such that $\pi(s)
\notin A$ and $\pi(t)\in A$. We say $\pi$ {\em enters $A$ from
outside} if there exists $t_\pi <s<t$ such that $\pi(s) \notin
\bar A$, the closure of $ A$, and $\pi(t)\in A$. We now recall the following
characterization of the Brownian net from \cite{SS07}.

\bteo\label{T:netchar}{\bf [Characterization of the Brownian net]}\\
There exists an $(\mathcal{H}, \mathcal{F}_{\mathcal{H}})$-valued random variable $\Net$, the
standard Brownian net (with parameter $\tau$), whose distribution is uniquely determined by
property {\rm (a)} and any of the three equivalent properties {\rm
(b1)--(b3)} below:
\begin{itemize}
\item[{\rm (a)}] There exist $\Wl, \Wr\subset \Net$ such that $(\Wl, \Wr)$ is
distributed as the left-right Brownian web.
\item[{\rm (b1)}] Almost surely, $\Net$ is the closure of $H(\Wl \cup\Wr)$ in
$(\Pi, d)$.
\item[{\rm (b2)}] Almost surely, $\Net$ is the set of paths in $\Pi$ which do
not enter any wedge of $(\hat\Wl, \hat\Wr)$
  from outside.
\item[{\rm (b3)}] Almost surely, $\Net$ is the set of paths in $\Pi$ which do
not enter any mesh of $(\Wl, \Wr)$.
\end{itemize}
\eteo 

As pointed out in \cite{SSS08}, the construction of the Brownian net from the left-right Brownian web can be regarded as an outside-in approach because
$\mathcal{W}_l$ and $\mathcal{W}r$ are the ``outermost'' paths among all paths in $\mathcal{N}$. On the other hand, the marking construction of the Brownian net  can be regarded as an inside-out approach. We start from a standard Brownian web, which consist of the
``innermost'' paths in the Brownian net, and construct the rest of the Brownian net paths by allowing branching at a
Poisson  set of  marked points  in the Brownian web.

\subsection{Equivalence of the Constructions}
\label{equivalence_between_constructions}
The main ingredient in the construction we just described is the pair $(\mathcal{W}_r,\mathcal{W}_l)$.
In order to prove the equivalence between the two constructions we first prove that the sets
of leftmost and rightmost paths of  $\mathcal{N}_{mark}$
(as defined in Subsection \ref{con::net})
are distributed as such a pair (see Proposition \ref{sticky_pairss}).

In Section \ref{modify_one_path}, $_{[1]} r_z$ was constructed from $B_z$ 
by switching
all the marked left $(1,2)$ points on $B_0$, the path of $\W$
starting from the origin. Analogously, we can define 
$_{[n]} r_z$ after switching all the marked left $(1,2)$ points on $B_0,B_1, ..., B_{n-1}$, where $B_k$ is the path starting from $z_k$.
As can easily be seen, the interaction between 
$_{[n]} r_z$ and the family $\{B_{i}\}_{i=0}^{n-1}$ is  local. Hence, Proposition \ref{cornerstone} implies that $\rz$ evolves like 
an independent Brownian motion away from $\{B_i\}_{i=0}^{n-1}$ and the interaction
between ${_{[n]}r_z}$ and $B_i$ when they meet is a sticky reflection.
More precisely,
we have the following immediate generalization of Proposition \ref{cornerstone}.
\bprop
\label{lll_stick}
For any deterministic $z$, there exists $B_z^{(n)}$, a standard Brownian motion
starting at $z$, so that ${_{[n]}r_z},\{B_k\}_{k=0}^{n-1}$ 
satisfy the following SDE.
\begin{eqnarray}
d _{[n]} r_z & = & d B_z^{(n)}(t) + 1_{\bigcup_{k=0}^{n-1} \{{_{[n]}  r_z}(t)= B_k(t)\}} \ \tau \   dt, 
\nonumber \\[0.1in]
d  B_k(t) \ d B_z^{(n)} (t) \ & = & \ 1_{_{[n]} r_z(t)=B_k(t)} dt, \nonumber \\[0.1in]
\forall t\geq T_{_{[n]} r_z,B_k}, & & _{[n]} r_z(t)\ \geq \ B_k(t).
\label{equation_stick_BM'}
\end{eqnarray}
\eprop
We now motivate the next proposition.
As $n\raw\infty$, $\{B_k\}_{k=0}^{n-1}$ ``fills'' more and more space of $\R^2$ 
and because ${_{[n]} r_z}$ sticks to the family $\{B_k\}_{k=0}^{n-1}$ it is intuitively clear that $1_{\bigcup_{k=0}^{n-1}\{_{[n]}  r_z= B_k\}}\approx1$ (see Lemma \ref{stick_on_B} for a precise 
version
of this statement). Hence for large $n$, the first part of
(\ref{equation_stick_BM'}) becomes
\begin{eqnarray}
 d {_{[n]} r_z}(t)&=& d B_z^{(n)}(t)+  1_{\bigcup_{k=0}^{n-1} \   \{_{[n]}  r_z= B_k\} }  \ \tau \   dt \label{lotofstickyness'} \label{e:ww}\\
 & \approx & d B_z^{(n)}(t)+  \tau \   dt.
\end{eqnarray}
Hence, for any $k\in\N$, we expect $({_{[n]} r_z},B_k)$ to converge as $n\raw\infty$ in distribution to a
pair $(r_z,B_k)$ satisfying
the following SDE.
\begin{eqnarray}
d r_z& = & d B_z^r + \tau dt, \ \  \nonumber \\
 \ \ d B_k  \ d  B_z^r  (t) \ & = &  \  1_{r_z(t)=B_k(t)}dt, \ \nonumber\\
\forall t\geq T_{r_z,B_{k}},& & r_z \ \geq \ B_k, \label{equ::needed}
\end{eqnarray}
where $B_z^r$ is a Brownian motion starting from $z$.

We recall that $\{z_i\}_{i=0}^{\infty}$ is a dense deterministic subset of $\R^2$. Let $i\in\N$. In the following, we write ${_{[n]}r_i}$ for ${_{[n]}r_{z_i}}$. Since $\{_{[n]} r_i\}_n$ is clearly increasing in $n$, the sequence $\{_{[n]} r_i\}_n$ actually converges {\it pathwise} and the 
limit is a drifting Brownian motion.
(Although it is not even clear a priori that the sequence of paths
is bounded, this will follow from the fact, as motivated by 
(\ref{e:ww})-(\ref{equ::needed}), that there is convergence in distribution.)
This pathwise
limit
will be referred to as $r_i$; it corresponds to the rightmost path of the net $\mathcal{N}_{mark}$ starting from $z_i$. In particular, any path of any partial net $\mathcal{N}_n(=\Net_{n,n})$ starting at $z_i$
is always to the left of $r_i$ (i.e., $\leq r_i$). This motivates the following proposition, whose proof is given in Subsection
\ref{pr3}.
\bprop
\label{Drifting_Brownian}
$_{[n]} r_i$ converges pointwise to a continuous path $r_i$ starting from $z_i$ with $(r_i,B_k)$ satisfying the
three-part SDE (\ref{equ::needed}).
\eprop

Analogously, using the set of marked {\it right}  $(1,2)$ points of $\W$, we can define $\{l_j\}_j$ a family of left-drifting Brownian motions reflected in a sticky way on the paths of $\W$. In Subsection \ref{pr3}  
we prove the following extension of Proposition \ref{Drifting_Brownian}.

\bprop
\label{sticky_pairss}
 $\{r_j\}_j$   (resp. $\{l_j\}_j$) is a family of coalescing right- (resp., left-) drifting Brownian motions with drift $\tau$ (resp., $-\tau$). The pair $(\mathcal{W}_l,\mathcal{W}_r)$, defined as the closures of $\{l_j\}_j, \{r_j\}_j$ respectively, is distributed as a left-right Brownian web.
\eprop

Now,
let $\mathcal{N}_{hop}$ denote the net obtained from  $(\mathcal{W}_r,\mathcal{W}_l)$ by the hopping
construction given in Section \ref{remainder}. In Subsection \ref{pr3}, we prove

\bteo
\label{equivalence_Net}
\begin{equation}
 \mathcal{N}_{hop}=\mathcal{N}_{mark}.
\end{equation}
\eteo

\section{The Dynamical Brownian Web}
\label{Dynamical::BW}
In order to describe the dynamical web, we will need the following notion of stickiness.
\bdf
\label{def::stick}{\bf [Stickiness]}
Let $\pi_1,\pi_2$ be in the net $\Net$ with $x=\pi_1(t)=\pi_2(t)$. We say that $\pi_1$ sticks to $\pi_2$ at $z=(x,t)$, or equivalently $\pi_2 \eq \pi_1$,  iff
for any $\epsilon>0$,
$$\int_{t}^{t+\epsilon} 1_{\pi_1(u)=\pi_2(u)} \ du >0 \ \txt{and} \  \int_{t-\e}^{t} 1_{\pi_1(u)=\pi_2(u)} du>0.$$
\edf

We now set up some notation. We say that a path enters a point $z=(x,t)$ if $t_\pi<t$ and 
$\pi(t)=x$. Let $z$ be a $(1,2)$ point
in $\Net_{mark}$. 
For any $B\in\W$ entering $z$, we denote by $B_{switch}$ the path obtained from $ B$ after switching the direction of $z$.  Since for any paths $\pi\in\Net_{mark}$ and $\bar B,B\in\W$ entering $z$, 
$\pi\sim^z B$ iff $\pi\sim^z \bar B$ ,
we will sometimes
write $\pi\sim^z B$ without specifying $B$ to mean that there exists a $B\in\W$ such that $\pi\sim^z B$.
Analogously, we will write $\pi\sim^z B_{switch}$, without specifying the path $B$ from which $B_{switch}$
was constructed.

Recall the partial dynamical web $\{\W_{n,m}(\tau)\}_{\tau\geq0}$ given in Subsection \ref{modify_the_Web2}. In the following, $\Net_{mark}(\tau)$ is the net constructed from $\M(\tau)$. The proof of the next proposition
is given in Subsection \ref{edep}. That proof makes clear
that the three parts of Proposition \ref{existence_of_Br} correspond
to 
three alternative constructions of the dynamical Brownian web. 
\label{dynamical_brownian_web}
\bprop
\label{existence_of_Br}
\begin{enumerate}
\item[(1)] There exists $\{\W(\tau)\}_{\tau\geq0}$ in $(\mathcal{H},d_{\mathcal{H}})$ s.t. almost surely 
$$\forall \tau \geq 0, \lim_{n,m\uaw\infty}d_{\mathcal{H}}\left(\W_{n,m}(\tau),\W(\tau)\right)=0.$$
\item[(2)]  
$ \W(\tau)   =  \{\pi\in\Net_{mark}(\tau) \ : \ \textrm{every time $\pi$ enters a point $z$ in $\M({\tau})$, $\pi \sim^z B_{switch}$}  \}.$
\item[(3)] Almost surely, $\W(\tau)$ satisfies the two following conditions (of Theorem \ref{teo:char}) for every $\tau\geq0$.
\begin{enumerate}
\item[(o)] From any deterministic point \( z \) in
          $\r^{2}$,
                    there is a unique path \( {B}_{z}^\tau \in \W(\tau)\)
          starting
                    from \( z \).
\item[(i)] \( {\W}(\tau) \) is the
          closure in
                    \( ({\cal H}, d_{{\cal H}}) \) of \( \{ {B}_{i}^\tau\} \) where ${B}_{i}^\tau$ is the unique path in $\W(\tau)$ starting from $z_i\in\mathcal{D}$. 
\end{enumerate}
\end{enumerate}
\eprop

To motivate item (2), note that in the partial dynamical web $\W_{n,m}(\tau)$,
any path $\pi$ entering a point $z\in\M_{n,m}(\tau)$ 
locally coincides with any path $B\in\W$ entering $z$ and then connects
to the newly born path starting from $z$. Hence, $\pi$ locally coincides with $B_{switch}$ and therefore obviously sticks to it. 
However, if $z$ belongs to $\M(\tau)\setminus\M_{n,m}(\tau)$, then $\pi\sim^z B$.
In the limit $n,m\raw\infty$, $\pi\sim^z B_{switch}$
for every $z$ in $\M(\tau)$.

We now turn to the description of some properties of the dynamical Brownian web.
We start with a definition.

\bdf
$(B,B')$ is a
$(1/\tau)$-sticky pair of Brownian motions iff 
\begin{enumerate}
\item $B$ and $B'$ are both Brownian motions starting at $(x_B,t_B)$ and $(x_{B'},t_{B'})$ that move 
independently when they do not coincide.
\item For $t\geq0$, define $B_{stick}(t)\equiv|B-B'|(t+t_B\vee t_{B'})/\sqrt{2}$. Conditioned on $x=B_{stick}(0)$, $\{B_{stick}(t)\}_{t\geq0}$ is a $(\sqrt{2}/\tau)$-sticky Brownian motion (see Definition \ref{def-1d-sbm}).
\end{enumerate}
\edf

We call $(B_1, \ldots, B_m; B_1',
\ldots, B_n')$ a collection of {\em ($1/\tau$)-sticking-coalescing Brownian
motions}, if $(B_1, \ldots, B_m)$ and $(B_1', \ldots, B_n')$ are
each distributed
as a set of 
coalescing
Brownian motions and for any $B \in
\{B_1, \ldots, B_m\}$ and $B'\in\{B_1', \ldots, B_n'\}$, $(B,B')$ is a
$(1/\tau)$-sticky pair of Brownian motions.

We will say that $(\W,\W')$ is a {\it $1/\tau$-sticky pair of Brownian webs} iff $(\W,\W')$ satisfies the following properties
\begin{itemize}
\item[{\rm (a)}] $\W$, resp.\ $\W'$, is distributed as the standard
Brownian web.

\item[{\rm (b)}] For any finite deterministic set $z_1, \ldots, z_m,
z'_1, \ldots, z'_n\in\R^2$, the subset of paths in $\W$ starting from
$z_1, \ldots, z_m$, and the subset of paths in $\W'$ starting from
$z'_1, \ldots, z'_n$, are jointly distributed as a collection of
($1/\tau$)-sticking-coalescing Brownian
motions.
\end{itemize}
Note that $(\W,\W')$ is defined in a similar way as $(\W_l,\W_r)$ except that  in (a) there is no drift and in
(b)
the collection of left-right coalescing
Brownian motions is replaced by the
collection of ($1/\tau$)-sticking-coalescing Brownian
motions.  
We are now ready to state the main result of this section whose proof is postponed to Subsection~\ref{rrrrrrrr}. 
\bteo
\label{real::stick}
 \begin{itemize}
\item[(a)]  $(\W,\W({\tau}))$ is a $1/(2\tau)$-sticky pair of Brownian webs.
\item[(b)]  (a Markov property). For $\tau_1\leq\tau_2$ and given $(\W,\{\M(\tau)\}_{\tau\leq\tau_1})$, the 
distribution of the 
pair $(\W({\tau_1}),\W({\tau_2}))$ only depends on $\W({\tau_1})$.
\item[(c)]  (Stationarity). For $\tau_1\leq\tau_2$, $(\W({\tau_1}),\W({\tau_2}))$ and $(\W,\W({\tau_2-\tau_1}))$ are equidistributed.
\item[(d)]  For any fixed deterministic time $t_0>0$,
the process $\tau\rightarrow B^{\tau}_0(t_0)$ is piecewise constant.
 \end{itemize}

\eteo

We remark that  existence of a
consistent family of finite dimensional
distributions for the process
$\W(\tau)$ follows from the results of \cite{HW07}---
see in particular Theorem 9 there.

\section{Proofs}\label{proof}
This section is organized as follows. In Subsection \ref{about_the _Brownian_web},  
we recall some useful properties of the Brownian web.
In Subsection \ref{existence-lo-time}, we complete the construction of the local time measure 
outlined in Subsection \ref{localtimesec}. In Subsection \ref{sec-excursions}, we 
carefully study some quantities related to the marked excursions of the web. Those results, 
whose proofs can be skipped at first reading, will be the key 
ingredients
in the proofs of Proposition \ref{cornerstone} (in Subsection \ref{proof:sticky})
and Theorem \ref{real::stick} (in Subsection \ref{rrrrrrrr}).  In Subsection 
\ref{poof:nett}, we 
provide a proof of the results from Section \ref{Brownian_Net} on the equivalence between the marking and the 
hopping 
constructions of
the Brownian net. 
In Subsection \ref{separation_points}, 
we give a 
proof of a basic fact relating the BN to $(1,2)$ points of the BW --- that
every ``point of separation'' in the BN is (in our coupling of the
BW and BN) also a $(1,2)$ point of the BW. 
We study some elementary properties
of the separation points in the Brownian net, and apply those results 
in Subsection \ref{rrrrrrrr} to prove Proposition
\ref{existence_of_Br} about the
existence of the
dynamical Brownian web. 
We note that the results about separation points of the Brownian
net
had already been derived by one of us (E.~S.) jointly with Sun and
Swart and will also appear in a paper
\cite{SSS08} by those three authors.

\subsection{Some Results about the Brownian Web}
\label{about_the _Brownian_web}

 We start by defining the age of a point $(x,t)$ as
\begin{equation}
\label{age29}
\sup\{t-t_B: \textrm{$B\in\W$ and $B(t)=x$}\}. 
\end{equation}
The $\g$-age truncation of the Brownian web is the set of paths obtained after
shortening every path of $\W$ by removing (if necessary) the initial segment 
consisting of those points of age less than $\gamma$.
In [FINR06] it was proved that:
\bprop
\label{aging}
The $\gamma$-age truncation of $\W$ is ``locally sparse'' in the sense that for every bounded set $U$,
the intersection between $U$ and the $\gamma$-age truncation of $\W$ only consists of finitely many path segments. 
\eprop
Two corollaries of that proposition can be formulated as follows:
\bcor
\label{one_coalescence_time}
Given $B$ and $\{B_n\}$ in $\W$ so that $B_n\to B$ (in $(\Pi,d)$) then the coalescence 
time of $B_n$ and $B$ converges to the starting time of $B$. 
\ecor
\begin{proof}
Let $t$ be the starting time of $B$ and take any $\bar t>t$. Let us consider the points $z_n$ (resp., $z$) where $B_n$ (resp., $B$)
intersect the line $\mathbb{R}\times\{\bar t\}$. The toplogy of $(\Pi,d)$ (see [FIN05]) implies that the starting time of $B_n$ converges to $t$.
Hence, for $n\geq n_0$ with $n_0$ large enough, $z_n$ has an age larger than $(\bar t-t)/2>0$. Moreover, since $z_n\rightarrow z$, the sequence $\{z_n\}$ belongs to a bounded segment of the line. By Proposition 1, we get that 
$\{z_n\}_{n\geq n_0}$ consist of only finitely many points. Therefore, ${z_n}$ is fixed after a certain n
and $B_n$ coincides with $B$ at $\bar t$.
Since this is valid for any $\bar t>t$, the claim of Corollary \ref{one_coalescence_time} follows.
\end{proof}

\bcor
\label{one_coalescence_time_2}
Let $B$ be a path in $\W$ starting at $t_0$. For any $\mathcal{D}$ as in Theorem \ref{teo:char} and $t>t_0$, on $[t,\infty)$ the path $B$ coincides with a path of the skeleton (determined by $\mathcal{D}$).
\ecor

\begin{proof}
By definition, there exists $B_n$ in the skeleton converging to $B$. The conclusion immediately  follows from the previous 
corollary.
\end{proof}

\subsection{Existence of the Local Time Measure}
\label{existence-lo-time}
In this section, we prove Proposition \ref{localtime} on which is based the construction of the local time measure. For simplicity of notation, we assume $(x,t)=(0,0)$ .

Let $(\bar B_1 ,\bar B_2)$ be two independent standard Brownian motion paths starting at $(0,0)$.
We define  $(B_{ind},\hat B_{ind})$ as
\begin{eqnarray}
B_{ind}(u) & = & \bar B_1(u), \nonumber \\
\hat B_{ind}(u) & = & \hat x+\bar B_2(u)-\bar B_2(\hat t) \ \ \ \mbox{for $u\in[0,\hat t]$}.
\end{eqnarray}
Clearly, $(B_{ind},\hat B_{ind})$ is a pair of independent forward and backward Brownian motions and we construct the system of refelected paths $(B,\hat B)$ as in Proposition \ref{forwardandbackward2}, i.e  $(B,\hat B)=(R_{\hat B_{ind}}(B_{ind}),\hat B_{ind})$.

\bigskip

In the following, we will assume that $\hat B(0)(=\hat B_{ind}(0)=\hat x-\bar B_2(\hat t))<0$. The case $\hat B(0)>0$ can be treated analogously, and $\hat B(0)=0$ has zero probability. Let $R_0(\bar B_1-\bar B_2)$ (resp., $R_0(B_{ind}-\hat B_{ind})$) be the Skorohod reflection of $\bar B_1-\bar B_2$ (resp., $B_{ind}-\hat B_{ind}$) at zero, i.e.,
\begin{eqnarray}
R_0(\bar B_1-\bar B_2)(u)    & = & (\bar B_1-\bar B_2)(u)- \min_{[0,u]} (\bar B_1-\bar B_2) \label{reflection_zero} \\
R_0(B_{ind}-\hat B_{ind})(u) & = &  (B_{ind}-\hat B_{ind})(u)- 0\wedge\min_{[0,u]} (B_{ind}-\hat B_{ind})\\
                            & = & (B-\hat B)(u).
\end{eqnarray}
Let $T_0$ be the first time $(B_{ind}-\hat B_{ind})$ hits $0$. 
Since $(B_{ind}-\hat B_{ind})$ is a translation of $\bar B_1-\bar B_2$ by $-\hat B_{ind}(0)>0$, Equation ~(\ref{reflection_zero}) immediately implies that
\begin{eqnarray}
R_0(B_{ind}(u)-\hat B_{ind})(u)=R_0(\bar B_1-\bar B_2)(u) \ \ \forall u \geq T_0, \\
R_0(B_{ind}(u)-\hat B_{ind})(u)\neq 0 \ \ \forall u < T_0.
\end{eqnarray}

$(\bar B_1-\bar B_2)/\sqrt{2}$ is a standard Brownian motion and it is a well known result (see, e.g., \cite{KS91}) that
$R_0(\bar B_1-\bar B_2)/\sqrt{2}$ is distributed as the absolute value of a Brownian motion and its local time at $0$, defined as
\begin{equation}
L(u)=\lim_{\epsilon\downarrow 0} \frac{1}{2 \epsilon} |\{v\leq u: \frac{R_0(\bar B_1-\bar B_2)(v)}{\sqrt{2}}<\epsilon\}|
\end{equation}
is equal to $- \min_{[0,u]} (\bar B_1-\bar B_2)/\sqrt{2}$. This implies that the quantity 
\begin{equation}
 L_{B,\hat B}(u) = \lim_{\epsilon\downarrow0} \frac{1}{2\epsilon} |\{v\leq u: \frac{1}{\sqrt{2}}R_0(B_{ind}-\hat B_{ind})(v)=\frac{1}{\sqrt{2}}(B-\hat B)(v) < \epsilon\}|
\end{equation}
is well defined and moreover
\begin{eqnarray}
\label{T0}
 L_{B,\hat B}(u) & = &  L(u\vee T_0)-L(T_0) \\
& = &-\min_{[0,u\vee T_0]} \frac{\bar B_1-\bar B_2}{\sqrt{2}}+ \min_{[0,T_0]} \frac{\bar B_1-\bar B_2}{\sqrt{2}} \\
& = & - \ \ 0\wedge\frac{\min_{[0,u]}(B_{ind}-\hat B_{ind})}{\sqrt{2}}.
\end{eqnarray}
This completes the proof of items 1 and 3 of Proposition \ref{localtime}.
\bigskip

Finally, item 2 follows from the fact (see Theorem 1 in \cite{P81}) that almost surely, 
the local time measure at zero of a Brownian motion is the $\phi$-Hausdorff measure of its zero-set.

\subsection{Excursions}
\label{sec-excursions}

To motivate this section, let us consider the pair $(\ro,B_0)$ (see Section \ref{modify_one_path}). On any interval of $\{t: B_0(t)\neq \ro(t)\}$, $\ro$ coincides with some path of the Brownian web other than $B_0$. Therefore, away from $B_0$, $\ro$ evolves as a Brownian motion independent of $B_0$ (this is part of the proof in Subsection \ref{proof:sticky} below of Proposition \ref{cornerstone}, which describes the distribution of $(\ro,B_0)$). Hence, to determine the distribution of $(\ro,B_0)$, we will 
need to analyze how $\ro$ escapes from the diagonal $\{t: \ro(t)=B_0(t)\}$.

Let us define $t_\e^r=\inf\{s: (\ro-B_0)(s)=\sqrt{2} \e\}$, the first time the pair $(\ro,B_0)$ escapes from 
the
$\sqrt{2}\e$-neighborhood of the diagonal.  
By construction, 
$t_\e^r$ is also the first time any right marked excursion is at a spatial distance $\sqrt{2}\e$ from $B_0$.
In Subsection \ref{distr-tere}, we give an explicit expression for the distribution of $t_\e^r$.
In Subsection \ref{distr-tere2}, we obtain  asymptotics for $\E(t_\e^r)$ for small $\e$. This will be used to prove Proposition \ref{cornerstone}. Finally,
we present
Proposition \ref{rihght-left-e} in Subsection \ref{sec-excursions3}---a result relating left and right excursions from $B_0$. It
will be used to prove Theorem \ref{real::stick}(a) which describes the joint
distribution of the dynamical Brownian web at two different dynamical times.

\subsubsection{Distribution of $t_\e^r$}
\label{distr-tere}
In this subsection, we will prove the following proposition.
\bprop\label{dist-ter}
Let $|B|_0^\e(t)$ 
be a Brownian motion on $[0,\e]$, 
starting at $0$ and reflected at $0$
and $\e$ and let $l_\e(t)$ be its local time at level $\e$. Then
\beq
\P(t_\e^r\leq t)=\P(l_\e(t)\geq \txt{Exp}[\fr{1}{\sqrt{2}\tau}]),
\eeq
where $\textrm{Exp}[ 1/(\sqrt{2}\tau)]$ is an exponential random variable with mean $1/(\sqrt{2}\tau)$, independent of $|B|_0^\e$.
\eprop

By definition, $t_\e^r\leq t$ iff a marked excursion enters the region
 \begin{equation}
\label{Rregion}
U_{\epsilon,t}=\{(x,u): 0\leq u \leq t,  B_0(u)+\sqrt{2} \epsilon \leq x\}.
\end{equation}
Equivalently, this condition can be re-expressed using the dual Brownian web.
\blem
$t_\e^r\leq t$ iff  there exists a backward path $\hat B$ starting from $U_{\e,t}$ and hitting $B_0$
at a marked point.
\elem
\begin{proof}
To show the only if part of the lemma, assume
there exists a right marked excursion $e_r$ from $B_0$
and $0\leq s \leq t$ such that $(e_r(s),s)\in U_{\e,t}$. 
One can then construct a sequence $\{\hat B_n\}$ in $\hat \W$ such that $\hat B_n$ starts at $(\hat x_n,\hat t_n)$ with $B_0(\hat t_n)<  \hat x_n < e_r(\hat t_n)$ and $(\hat x_n,\hat t_n)\raw(e_r(s),s)$.
Since paths of the web and its dual do not cross, $\hat B_n$ is squeezed between $e_r$ and $B_0$ and thus  enters the marked starting point $z$ of $e_r$. By compactness of $\hat \W$,  $\hat B_n$ converges (along a subsequence) to some path $\hat B\in\hat \W$ starting at $(e_r(s),s)\in U_{\e,t}$ and entering the point $z$.  The converse argument to prove the if part of the lemma is similar.
\end{proof}
We denote by $L_{\e,t}([t_1,t_2])$ the local time measure of all the points in $\R\times[t_1,t_2]$ where $B_0$ meets a backward path starting from $U_{\e,t}$. This naturally defines a measure $L_{\e,t}$ 
on $\R$  and we set $L_{\e,t}([0,t])\equiv \tilde l_\e(t)$.
By definition, the set of marked points at the intersection
between $B_0$ and the set of backward paths starting from $U_{\e,t}$ is a Poisson point process with intensity
$\sqrt{2} \tau \tilde l_\e(t)$. Hence,
\begin{equation}
\label{vraiment_trop}
 \P(t_\e^r\leq t) \ \ = \ \ \P(\tilde l_{\e}(t)\geq\textrm{Exp}[\frac{1}{\sqrt{2}\tau}]), 
\end{equation}
where $\textrm{Exp}[1/(\sqrt{2}\tau)]$ is independent of $\W$.

\begin{figure}[!ht]
         \begin{center}
          \includegraphics[width=8cm]{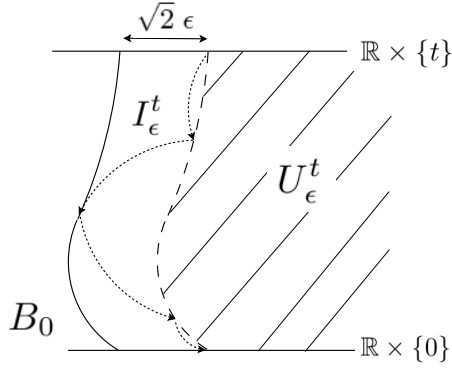}
         \caption{The process $I_\e^t$ is the left envelope of all the backward paths starting from the region $U_\e^t$.}\label{Iprocess} %don't know what this does
          \end{center}
          \end{figure}

To study the measure $L_{\e,t}$, we introduce the (backward) process $I_\e^t$ (see Figure \ref{Iprocess}) defined as
\begin{equation}
\forall s\in[0,t], \ \  I_\e^t(s)=\inf\{\hat B(s): \hat B \in\hat \W, z(\hat B) \in U_{\e,t}\},
\end{equation}
where $z(\hat B)$ denotes the starting point of $\hat B$.

Not surprisingly, the set of times when $I_\e^t$ and $B_0$ coincide is the support of $L_{\e,t}$. This claim can be verified as follows. Because of the compactness of
$\hat{\mathcal{W}}$, the time it takes for a path in $\hat \W$ starting from $U_{\epsilon,t}$ to reach
the curve $B_0$ is uniformly bounded away
from $0$. This means that the (backward) age of those paths (see (\ref{age29})) is strictly positive and the claim follows directly from Proposition ~\ref{aging}. Proposition \ref{dist-ter} directly follows from 
(\ref{vraiment_trop}) and the following lemma.

\blem
\label{identical_set_of_points}
\begin{enumerate}
The process $|B|_0^\e$ defined on $[0,t]$ by
$$|B|_0^\e(s)\equiv-\frac{1}{\sqrt{2}}\left( I_\e^t(t-s)-B_0(t-s)-\sqrt{2}\e\right)$$
is a Brownian motion on $[0,\e]$ starting at $0$ and reflected at $0$ and $\e$.
\end{enumerate}
\elem

\begin{proof}
Let  $\{\hat B_{k,n}\}_{n\in\N,k\in\{0,...,2^n\}}$ be the family of backward paths starting from  points of the form 
$z_{k,n}=(B_0(k t/2^n)+\sqrt{2}\epsilon, k t/2^n)$. We define
\begin{equation}
 \{\pi_{k,n}\}=\{\frac{-1}{\sqrt{2}}\left(\hat B_{k,n}(t-s)-B_0(t-s)-\sqrt{2}\epsilon \right): t-k t/2^n \leq s \leq t\}.
\end{equation}
Clearly, $\{\pi_{k,n}\}$ starts from $\{(0,t-k t/2^n)\}$ and is identical in law with a family of forward coalescing Brownian motions Skorohod reflected at $\epsilon$.

As can easily be seen, the process
\begin{equation}
_n |B|_0^\e(u)\equiv\sup\{\pi_{k,n}(u): k\in\{0,\dots,2^n\}\}
\end{equation}
converges pointwise to $|B|_0^\e$ as $n$ goes to $\infty$.

Now, let us decompose the process $|B|_0^\e$ into its up and downcrossings (the first upcrossing is the section of the path on $[0,t_\e^1]$, where $t_\e^1$ is the first time $|B|_0^\e$ hits $\e$; the first downcrossing is the section of the path between $t_1^\e$ and its return time to $0$). We aim to prove that an upcrossing (resp., downcrossing) is a copy of an independent Brownian motion starting at $0$ (resp., $\epsilon$), reflected at $0$ (resp., $\epsilon$) and stopped when it hits $\epsilon$ (resp., $0$). It is straightforward to show
the equidistribution and independence of the up and downcrossings.
The downcrossings have the required distribution because $|B|_0^\e$ coincides with $\pi_{k,n}$ for some $n$ and $k$ during a downcrossing. It remains to determine the law of the upcrossings.
\begin{figure}[!ht]
         \begin{center}
          \includegraphics[width=10cm]{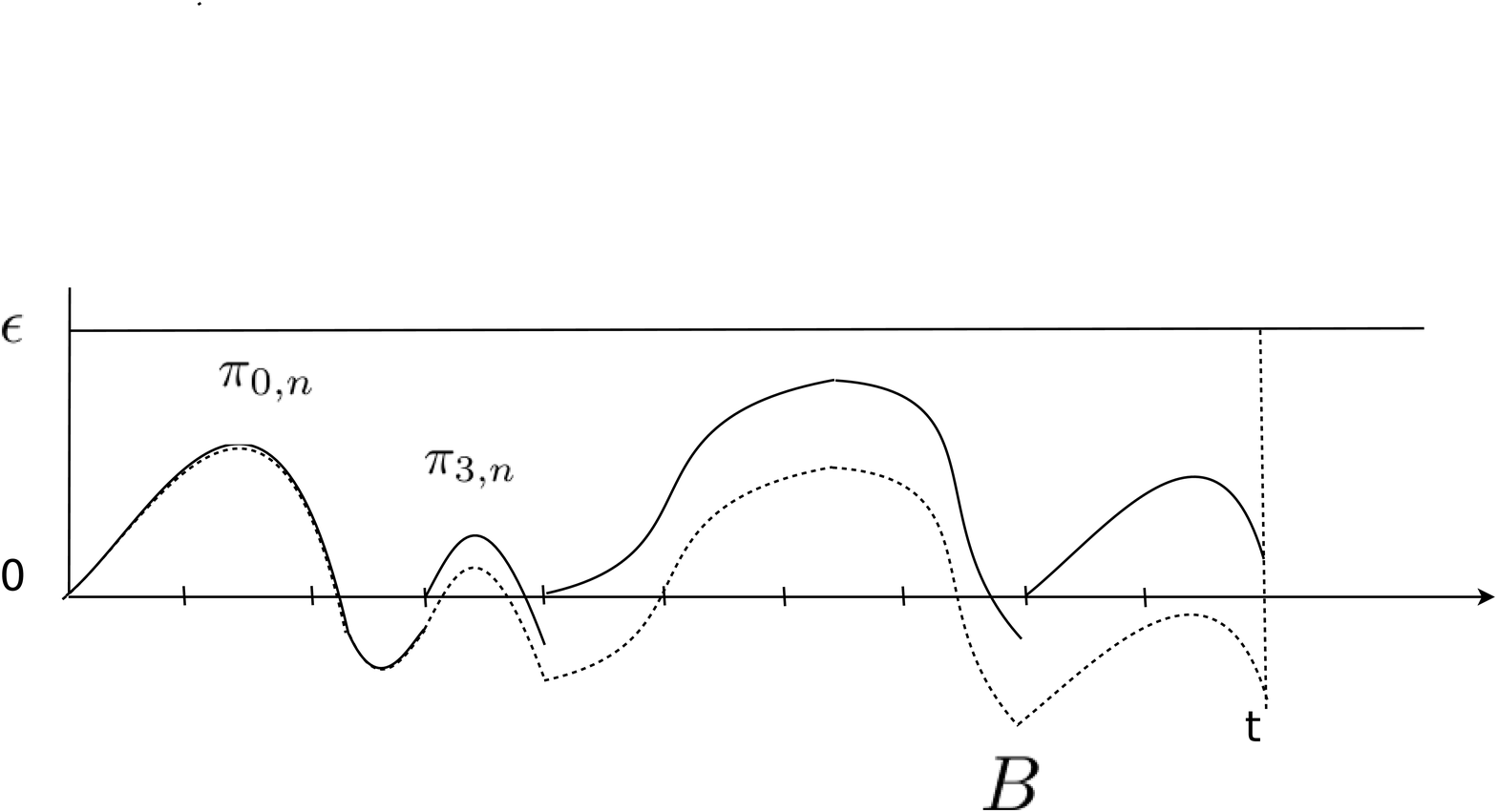}
          \caption{The continuous dashed path $B$ is constructed from the plain path $_n |B|_0^\e$.}\label{diagram2}%don't know what this does
          \end{center}
\end{figure}
Let $u_1$ (resp., $u_{1,n}$) be the first upcrossing of the process $|B|_0^\e$ (resp., $_n |B|_0^\e$). $u_{1,n}$
is simply made of pieces of Brownian motions stopped if they hit $\epsilon$. Let $B$ (which depends on $n$) be the continuous process starting at $0$ and obtained by gluing those pieces at their endpoints (see Figure \ref{diagram2}).  
By a simple induction, it is easy to see that
\begin{equation}
\label{reflex_contr}
  \forall \ t \in[k/2^n,(k+1)/2^n), \ \ \ \ \ u_{1,n}(t)=B(t)- \inf\{B(\frac{j}{2^n}): j=0,1/2^n,...,k/2^n\}.  
\end{equation}
and by the Markov property, $B$ is a Brownian motion stopped when $u_{1,n}$ hits $\epsilon$.
As $n\raw\infty$, the right hand side of ~(\ref{reflex_contr}) converges in law to
 \begin{equation}
B(t)- \inf_{[0,t]} B
\end{equation}
where $B$ is a Brownian motion stopped when $B(t)-\inf_{[0,t]} B$ hits $\epsilon$. On the other hand, the left hand side of ~(\ref{reflex_contr}) converges almost surely to $u_1$. Hence, the first upcrossing of $|B|_0^\e(s)$ is identical in distribution with that of a Brownian motion starting at $0$, Skorohod reflected at $0$ and stopped when it hits $1$.
\end{proof}

\subsubsection{Rate of Excursions From $B_0$}
\label{distr-tere2}
In this subsection, we prove 
\bprop \label{rate-ex-rpath}
$\lim_{\e\daw0} \ \E(t_\e^r)/\e=\fr{\sqrt{2}}{\tau} \ \ \mbox{and}
\ \ \ \ \E([t_\e^r]^2)=o(\e)$ as $\e\daw0$.
\eprop

We only prove the first claim. The second one 
can be proved along the same lines.

Let $\P_\W$ denote the probability distribution of the marked Brownian web 
conditioned on the web $\W$. By Proposition \ref{dist-ter}, we have the following. 
\begin{eqnarray*}
{{\mathbb{E}(t_\e^r)}}/{\e} & =\int_0^\infty \P(t_\e^r>\e t) dt  & = \int_0^\infty \E(\P_\W(t_\e^r>\e t)) dt
 \\
& =  \int_0^\infty  \E(\P_\W \left(l_\e(\e t)< \textrm{Exp}[\frac{1}{\sqrt{2}\tau}]\right))  \ dt
& =  \int_0^\infty \E\left(\exp[- \sqrt{2} \tau \dd l_\e(\e t)]\right) \ dt. \label{appper}
\end{eqnarray*}
To take the limit as $\e\raw0$, we will use the following lemma.

\blem
\label{conv::LT}
Let $t,\g>0$. There exist $c,C\in(0,\infty)$ such that
\begin{equation}
\label{proff::reflec}
 \P(|l_\e(\e t)-\frac{t}{2}|> \g t) \leq C \exp(- c \frac{t \g}{\e}). 
\end{equation}
Hence, $l_\e(\e t)$ converges in probability to $t/2$ as $\e\raw 0$.
\elem

\begin{proof}

We need to show that
\begin{eqnarray}
 \P(l_\e(\e t)-\frac{t}{2}>  t \g) \leq C \exp(-c\frac{t \g}{\e}),\\
\P(\frac{t}{2}-l_\e(\e t) >  t \g) \leq C \exp(-c\frac{t \g}{\e}).
\end{eqnarray}
We only prove the first inequality. 
The second one can be obtained using analogous arguments.
Using the scaling invariance of Brownian motion, the first inequality reduces to
\begin{equation}
\label{large::dev::local}
 \P(\e l_1(t/\e)-\frac{t}{2} >  t \g) \leq C \exp(-c\frac{t \g}{\e}),
\end{equation}
where $l_1(u)$ is identical in distribution to the local time accumulated on the set $\{x=2j+1\}_{j\in\Z}$ at time $u$ by a standard Brownian motion $B$. Define $t_0=\inf\{s:B(s)=\pm 1\}$ and for $k\geq 1$, $t_k=\inf\{t\geq t_{k-1}: |B(t)-B(t_{k-1})|=2\}$. $\D t_k=t_{k+1}-t_{k}$ has mean $4$. Furthermore, by excursion theory, the local times $\D l_k$ accumulated on $\{x=2j+1\}_{j\in\Z}$ during the time intervals $[t_{k},t_{k+1}]$, for $k\geq0$ are independent exponential random variables with mean 2.  

Define $N_\e(t)=\inf\{k: t_k \geq t/\e\}$. Then, if we set 
$ \gamma' = (1+\g)$ and $n=\frac{t}{4\e}(1+\g) =\frac{t \gamma'}{4\e}$,
\begin{eqnarray*}
\P\left(\e l_1(t/\e)-\frac{t}{2}>  t \g\right)& \leq & \P\left( N_\e(t)>n \right)+
\P\left([\e \sum_{k\leq n} \D l_k] - \frac{t}{2} >  t \g\right) \\
& \leq & 
\P\left(\sum_{k\leq n} \D t_k < \frac{t}{\e}\right) + 
\P\left(\e \sum_{k\leq n}[\D l_k-2]>  \frac{t \g}{2} \right) \\
&\leq &
\P\left(\frac{1}{n}\sum_{k\leq n} [4-\D t_k] 
% \geq \frac{4\g}{1+\g}-\frac{16 \e}{t(1+\g)} \right) + 
%\P\left(\frac{1}{n} \sum_{k\leq n} [\D l_k-2]>  \frac{2 \g}{1+\g}\right).
\geq \frac{4\g}{\gamma'}-\frac{16 \e}{t \gamma'} \right) + 
\P\left(\frac{1}{n} \sum_{k\leq n} [\D l_k-2]>  \frac{2 \g}{\gamma'} \right).
\end{eqnarray*}
Equation (\ref{large::dev::local}) follows by classical large deviation estimates.
\end{proof}

To complete the analysis of $\lim_{\e\raw0} \E(t^r_\e)/\e$, we use Lemma \ref{conv::LT} with $\g=1/4$ to see that 
\begin{eqnarray}
\E(\exp(-\sqrt{2} \tau \dd l_\e(\e t)))& \leq & \exp(-\tau \frac{\sqrt{2} t}{4})+\P(l_{\e}(\e t)\leq \frac{t}{4})  \\
& \leq &  \exp(-\tau \frac{\sqrt{2}t}{4})+  C \exp(- c \frac{t}{4 \e}).    
\end{eqnarray}
It follows that the family $\{\P(t_\e^r\geq \e \dd)\}_{\e\leq1}$ is uniformly integrable. Therefore, by Lemma \ref{conv::LT} 
\begin{eqnarray*}
\lim_{\e\daw0}  \int_0^\infty \E(\exp(-\sqrt{2} \tau \dd l_\e(\e u)) du  =  
\int_0^\infty  \lim_{\e\daw0}  \E(\exp(-\sqrt{2} \tau \dd l_\e(\e u)) du 
& = %\int_{0}^\infty \exp(-\sqrt{2} \tau \dd \frac{u}{2}) du
\int_{0}^\infty e^{-\sqrt{2} \tau \dd u/2} du
=\frac{\sqrt{2}}{\tau}.
\end{eqnarray*}
This completes the proof of Proposition \ref{conv::LT}.

\subsubsection{Marked Right and Left Excursions}
\label{sec-excursions3}
Let $e_l$ be a  left marked excursion from $B_0$. We say that $T(e_l)$ (the starting time of $e_l$)
is straddled by the right excursion $e_r$ iff $T(e_r)<T(e_l)<T(e_r)+D(e_r)$. In this subsection, we prove the following.
\bprop\label{rihght-left-e}
Let $e_{l,\e}$ be the first left marked  excursion from $B_0$ with diameter (see Definition \ref{df-excursion}) greater or equal to $\sqrt{2}\e$. Then
\beq
\P(\ T(e_{l,\e}) \ \txt{is straddled by some right marked excursion $e_r$}\ )\raw 0, \ \ \txt{as $\e\daw0$}.
\eeq
\eprop
Define $A_\e\equiv\{\ T(e_{l,\e}) \ \ \txt{is straddled by some right marked excursion $e_r$}\ \}$ and $t_\e^l$
as the left analog of $t_\e^r$ (so that $t_\e^l$ is the first time $t$ that $B_0(t)-e_{l,\e}(t)=\sqrt{2}\e$ and therefore $t_\e^l
\geq T(e_{l,\e})$).

For any $H>0$, 
\begin{eqnarray*}
\limsup_{\e\downarrow0}  \P\left( A_\e  \right) 
& \leq &
\limsup_{\e\downarrow0}  \P(A_\e \ , \ t_\e^l \leq \e H ) +\limsup_{\e\downarrow0} \P(t_\e^l\ >\  \e H) 	\\	 
& \leq & 
\limsup_{\e\downarrow0}  \P(A_\e \ , \ t_\e^l \leq \e H ) \ + \limsup_{\e\downarrow0} \fr{\E(t_\e^l)}{\e H} \\
& = &
\limsup_{\e\downarrow0}  \P(A_\e \ , \ t_\e^l \leq \e H)
+ \frac{\sqrt{2}}{\tau H} 
\end{eqnarray*}
where the equality follows from Proposition \ref{rate-ex-rpath} and the the identity $\E(t_\e^r)=\E(t_\e^l)$. Since $H$ can be made arbitrarily large, in order to prove Proposition \ref{rihght-left-e} it suffices
to show that for any $H>0$,
\begin{equation}
\label{unedernier}
 \limsup_{\e\downarrow0} \P\left(A_\e \ , \ t_\e^l \leq \e H \right)  =0.
\end{equation}
Let $\e_n=(\e H)/2^n$ for $n\geq0$ and let $\e_{-1}=+\infty$.
Breaking up $A_\e$ accordingly to the duration of the excursion $e_r$ straddling $T(e_{l,\e})$,
we have
\begin{eqnarray}
\label{finir}
\P(A_\e \ , \ t_\e^l \leq \e H) 
 & = &\P(\exists \textrm{a right marked excursion $e_r$ with} \\
&&  T(e_r) \leq T(e_{l,\e}) \leq T(e_r)+D(e_r)\ , \ t_\e^l\leq\e H) \nonumber \\
& = &
\sum_{n\geq -1}\P(C_{n}') \ \leq \   \sum_{n\geq -1}\P(C_{n}),              
\end{eqnarray}
where
\begin{eqnarray}
 C_n'= & \{\exists \ \txt{a  right marked excursion 
$e_r$ with $D(e_r)\in [{\e_{n+1}},{\e_{n}})$ s.t.} \nonumber \\
 & T(e_r) \leq T(e_{l,\e}) \leq T(e_r)+\e_n \ , \ t_\e^l \leq\e H\}, \\
C_n= & \{\exists \ \txt{a  right marked excursion 
$e_r$ with $D(e_r)\geq \e_{n+1}$ s.t.} \nonumber \\
 & T(e_r) \leq T(e_{l,\e}) \leq T(e_r)+\e_n, \ T(e_{l,\e}) \leq\e H \}.
\end{eqnarray}
Let $\P_{L,\W}$ be the probability distribution of the marked Brownian web conditioned on
$\W$ and the marking of the left $(1,2)$ points.
Since given the Brownian web, the markings of the left and the right $(1,2)$ points are independent, we get
\begin{eqnarray}
 \P_{L,\W}(C_n) & = & \sqrt{2} \tau \Lo_{D(e_r)\geq \e_{n+1}} \left(\R\times[T(e_{l,\e})-\e_n \ , \ T(e_{l,\e})]\right) \\
& = & \sqrt{2} \tau L_{D(e_r)\geq \e_{n+1}} \left([T(e_{l,\e})-\e_n \ , \ T(e_{l,\e})]\right) 
\end{eqnarray}
where $\Lo_{D(e_r)\geq \e_{n+1}}$ is the local time measure on the possible
starting points in $\R^2$ of a right excursion from $B_0$ with $D(e_r)\geq \e_{n+1}$, and $L_{D(e_r)\geq\e_{n+1}}$ is the projection of that measure along the $t$-axis. Let $n\geq0$. Since $T(e_{l,\e})\in[0,\e H]$, there exists $k\in\{-1,0,...,2^n-2\}$ such that
\beq
[T(e_{l,\e})-\e_n \ , \ T(e_{l,\e})]\subset T_{k,n}, \ \txt{with} \ \ T_{k,n}=[k\e_n, \ (k+2) \e_n].
\eeq
Hence,
\begin{eqnarray}
 \P_{L,\W}(C_n) & \leq & \sqrt{2} \tau \max_{-1\leq k\leq 2^n-2} L_{D(e_r)\geq \e_{n+1}} \left(T_{k,n}\right), \\[0.1in]
 & \leq & \sqrt{2} \tau \max_{0\leq k\leq 2^n} L_{D(e_r)\geq \e_{n+1}} \left(T_{k,n}\right) \label{last!!}
\end{eqnarray}
where we used the equality $\Lo(\R\times T_{k,n})=\Lo(\R\times[k\e_n\vee 0, (k+2)\e_n\vee 0])$ to deduce the second inequality. 
Note that
with the convention that $T_{0,-1}=[0,\e H]$, the formula above also remains valid for 
$n=-1$. Averaging over the realizations of $\W$ and the  marking of left $(1,2)$ points,
we obtain that for any $p\geq 1$,
\begin{eqnarray}
 \P(C_n)&\leq \sqrt{2}\tau \ \E(\max_{0\leq k\leq 2^n} L_{D(e_r)\geq \e_{n+1}} \left(T_{k,n}\right)) \\
& \leq C_p \  |\max_{0\leq k\leq 2^n} L_{D(e_r)\geq \e_{n+1}} \left(T_{k,n}\right)|_p \label{enfin!!!}
\end{eqnarray}
where $C_p$ is a finite positive constant and $|X|_p$ denotes the $L_p$ norm of $X$ w.r.t. $\P$.

\blem
\label{lemma:er}
For any $p\geq 1$  there exists $K<\infty$ s.t. for $n\geq-1$,
 \begin{eqnarray}
 |\max_{0\leq k\leq 2^n} L_{\{D(e_r)\geq \e_{n+1}\}}(T_{k,n})|_p & \leq & K \ 2^{n(\frac{1}{p}-\frac{1}{2})}\sqrt{\e H}.
\end{eqnarray}
\elem

\begin{proof}
We prove the lemma for $n\geq0$. The case $n=-1$ (where $T_{0,-1}=[0,\e H]$) can be treated analogously.
By translation invariance of the marked Brownian web,
\begin{eqnarray}
\P\left[L_{\{D(e_r)\geq \e_{n+1}\}}([k\e_{n},(k+2)\e_{n}])>x\right] \
& = &
\P\left[L_{\{D(e_r)\geq \e_{n+1}\}} ([0,2\e_{n}])>x\right]. 
\end{eqnarray}
Therefore,
\begin{eqnarray*}
\P(\max_{0 \leq k\leq 2^n} L_{\{D(e_r)\geq\e_{n+1}\}}(T_{k,n})>x) \
& \leq &
2^{n+1} \P(L_{\{D(e_r) \geq \e_{n+1}\}} ([0,2\e_{n}])>x). 
\end{eqnarray*}
The scaling invariance of the Brownian web under the mapping on paths, $B \rightsquigarrow \lambda^{-1/2} B(\lambda t)$, yields (for $a_0,b_0\geq0$) the equidistribution of $L_{\{D(e_r)\geq a_0 \lambda\}}([0,b_0\lambda])$ and $\sqrt{\lambda} L_{\{D(e_r)\geq a_0 \}}([0,b_0])$. Hence
$$
L_{\{D(e_r)\geq \e_{n+1}\}}([0,2\e_{n}])=_{d} \ \ 
\sqrt{\frac{\e H}{2^{n}}} \ \ L_{\{D(e_r)\geq \frac{1}{2}\}}([0,2]),
$$
which, using the standard identity that $(|X|_p)^{p}$ equals $\int_0^\infty p x^{p-1} \P(|X|>x) dx$, implies that
\begin{eqnarray*}
|\max_{0\leq k\leq 2^n} L_{\{D(e_r)\in T_n\}}(T_{k,n})|_p  & \leq & 
  2^{1/p} 2^{n/p}  \
\left(p \ \int_{0}^\infty x^{p-1} \P\left[L_{\{D(e_r)>\frac{1}{2}\}} ([0,2])>x\sqrt{\frac{2^n}{\e H}}\right] \  dx\right)^{\frac{1}{p}} \\
& = & 2^{1/p} 2^{n(\frac{1}{p}-\frac{1}{2})}\sqrt{\e H} \ \ \ |L_{\{D(e_r)>\frac{1}{2}\}} ([0,2])|_p .
\label{autrefin}
\end{eqnarray*}
To complete the proof, we need to show that for any $p\geq1$,
$|L_{\{D(e_r)>\frac{1}{2}\}} ([0,2])|_p<\infty$.

We use the fact (see, e.g., \cite{FINR04}) that for any $s>0$, there are two distinct dual Brownian paths starting from $(B_0(s),s)$, those two paths being separated by the path $B_0$. 
In order for $s\in[0,2]$ to be in the support of $L_{\{D(e_r)>\frac{1}{2}\}}$, $B_0$ must 
be hit  by a (dual) path of $\hat \W$ starting in the region $\{(x,t): x\geq B_0(t),\ \ t\geq s+1/2\}$.
At any such time $s$, there must be an integer $k$ in $\{1,...,10\}$ such that $B_0$ is hit by
$\hat B_{k/4}$, the dual path starting at $(B_0(k/4),k/4)$ and located to the right of $B_0$. This implies 
that $L_{\{D(e_r)>\frac{1}{2}\}}$ is bounded above by the local time measure induced
by the finite family of backward paths $\{\hat B_{k/4}\}_{k\leq 10}$. From \cite{STW00} (see Proposition \ref{forwardandbackward2} above), the process
\begin{equation}
  \ s\raw \hat B_{k/4}(k/4-s)-B_0(k/4-s)
\end{equation}
defined on $[0,k/4]$ is a Brownian motion reflected at $0$ and the local time 
measure $L_{B_0,\hat B_{k/4}}$ is just the usual local time measure at the origin of that reflecting Brownian motion. It is a standard fact that local time at the origin has all moments
and Lemma \ref{lemma:er} follows.
\end{proof}

Combining (\ref{finir}), (\ref{enfin!!!})
and Lemma \ref{lemma:er} for any $p>2$, there exists $C_p'< \infty$ s.t.
\begin{equation}
\P(A_\e, \ t_\e^l \leq \e H)\leq  C_p' \sqrt{\e H},
\end{equation}
so that (\ref{unedernier}) and hence Proposition \ref{rihght-left-e} follow.

\subsection{Distribution of $(B_0,\rz)$ (Proof of Proposition \ref{cornerstone})}
\label{proof:sticky}

First, we prove the following lemma.  
\blem
\label{Markov_Pair}
The family  $\{(B_0,{_{[1]}r_z})\}_{z\in\R\times\{0\}}$ of random pairs of
continuous paths is a
family of strong Markov processes with stationary transition probabilities. 

More precisely, for any stopping time $T$, conditioned on the past of the paths up to $T$, 
i.e.,  conditioned on $\{\mathcal{F}_T\}$ (where $\mathcal{F}_t$ is the $\sigma$-field generated by $\{\left(B_0(s),\ro(s)\right)\}_{s\leq t}$ and $\{\mathcal{F}_T\}$ is defined accordingly), 
$$\left(B_0(t+T)-B_0(T), \rz(t+T)-B_0(T)\right)_{t\geq 0}$$ is distributed like $(B_0, {_{[1]}r_{z(T)}})$ 
with $z(T)=\rz(T)-B_0(T)$.
\elem
\begin{proof}
We take $z=(0,0)$,
first prove the weak Markov property and then the strong Markov property. 
The proof can trivially be extended to any deterministic~$z$.

\bigskip

\textbf{Weak Markov Property}: 
Recall that $\mathcal{L}_{1,n}$ is the natural local time measure on the set $E_{n}$ defined as
$$
E_{n}=B_0 \bigcap \left(\bigcup_{i=0}^{n-1}\hat B_i\right).
$$ 
For the time being, $T>0$ is deterministic. In the following, $E_{n}^{-}$ will denote the subset of $E_{n}\bigcap\{t\leq T\}$ consisting of all the points on $B_0$
hit by a path $\hat B_i$ starting from $z_i=(x_i,t_i)$ with $i\leq n-1$ and $t_i\leq T$. $E_{n}^{+}$ will refer to $E_{n}\bigcap\{t\geq T\}$. 
Finally, we define ${_{(1,n)}\bar r_{0}}$ as the path constructed from $B_0$ by switching the direction of the marked left $(1,2)$ points in $E_{n}^{+}\bigcup E_{n}^{-}$.

Let $\mathcal{L}_{(1,n)}^{+}$ (resp., $\mathcal{L}_{(1,n)}^{-}$) be the measure $\mathcal{L}_{(1,n)}$ restricted to $E_n^+$ (resp., $E_n^-$). First, conditioned on the Brownian web, the markings of $E_n^{+}$ and $E_n^{-}$ are two independent Poisson point processes with respective intensity measure $\mathcal{L}_{(1,n)}^+$ and $\mathcal{L}_{(1,n)}^-$. 
Second, $(\mathcal{L}_{(1,n)}^+,\{{_{(1,n)}\bar r_{0}}(t)\}_{t\geq T})$ (resp., $\mathcal{L}_{(1,n)}^-$) is measurable w.r.t. 
$$\left(\ \W_{[T,\infty]} \ , \ B_0(T) \ , \ {_{(1,n)}\bar r_0}(T) \ \right) \ \ \ \textrm{(resp., $\W_{[-\infty,T]}$)}, $$
 where $\W_{[t_1,t_2]}$ is the set of paths in $\W$ starting in the window $[t_1,t_2]$ and stopped at $t_2$.  By independence of $\W_{[T,\infty)}$ and $\W_{[-\infty,T]}$,  the future evolution of $(B_0,_{(1,n)}\bar r_{0})$ is independent of its past given $(B_0(T),{_{(1,n)}\bar r_0}(T))$.
Assuming momentarily that $_{(1,n)}\bar r_{0}$ converges pointwise to $\ro$, it is straightforward
to show that $\ro$ also continues afresh at $T$ provided
that the distribution of $(B_0,_{[1]} r_{\bar z})$, with $\bar z=(\bar x,0)$, is continuous with respect to $\bar x$.
This we will do next.
The stationarity of transition probabilities in Lemma \ref{Markov_Pair} then simply follows from the translation invariance of the marked Brownian web.  

We now prove that $(B_0,_{[1]} r_{(\bar x, 0)})$ is continuous with respect to $\bar x$.
 Let $\bar z_n=(\bar x_n,0)\rightarrow \bar z$. We distinguish between two cases:
\begin{enumerate}
\item $\bar z=(\bar x,0)$ with $\bar x\neq0$. Before meeting $B_0$, ${_{[1]} r_{\bar z}}$ (resp., ${_{[1]} r_{\bar z_n}}$) follows 
$B_{\bar z}$ (resp., $B_{\bar z_n}$), the path in $\W$ starting from $\bar z$ (resp., $\bar z_n$).
For $n$ large enough, $B_{\bar z}$ and $B_{\bar z_n}$ coalesce at some time $\mu_n$ before either of those paths  meets $B_0$. 
Hence, ${_{[1]} r_{\bar z}}$ and ${_{[1]} r_{\bar z_n}}$ coalesces at time $\mu_n$ with $\mu_n\rightarrow0$ as $n\uparrow\infty$. 
\item $\bar z=(0,0)$. For any $\g>0$, we can always find a marked left $(1,2)$ point
at $(B_0(t),t)$ for some
$t\in[0,\gamma]$. Let $\hat B\in\hat \W$ pass through that mark and let $(x_M,t_M)$ be the earliest of the marks along $\hat B$. Since almost surely $(0,0)$ is not a $(1,2)$ point, $t_M$ is strictly positive
and for $n$ large enough $0<x_n<\hat B(0)$. 
For $n$ large enough, $B_{\bar z_n}$ coalesces with $B_0$ before $t_M$. 
By construction, ${_{[1]}r_{z_n}}$ and $\ro$ can only cross $\hat B$ at a marked point on $B_0\cap\hat B$.
Since $t_M$ is the earliest marked point on $\hat B$, ${_{[1]}r_{\bar z_n}}$ and $\ro$ are squeezed between $B_0$ and $\hat B$ on $[0,t_M]$ and thus they meet (and coalesce) by $t_M\leq \g$. 
\end{enumerate}

For the weak Markov property, it remains to prove that $_{(1,n)}\bar r_0$ converges to $\ro$. Recall that the excursions 
of $\ro$ from $B_0$ coincide with the maximal excursions from $B_0$ (see Definition \ref{df-excursion}). 
First, let $e$ be a maximal excursion starting at some
$z$. For $n$ large enough, it is clear that $z$ belongs to $E_{n}^{+}\bigcup E_{n}^{-}$. By definition of a maximal excursion, $_{(1,n)}\bar r_0$ hits $z$ and then follows $e$.
Second, let $z'$ be the starting point of a marked excursion $e'$ which is not maximal and hence is
nested in some maximal excursion. For $n$ large enough, $_{(1,n)}\bar r_0$ follows that maximal excursion and therefore misses the excursion $e'$. Hence, in the limit, the excursions of $_{(1,n)}\bar r_0$ coincide with the maximal excursions from $B_0$, and thus 
$_{(1,n)}\bar r_{0}$ converges pointwise to $\ro$.

\bigskip 

\textbf{Strong Markov Property}: Now let $T$ be a stopping time 
with respect to the right-continuous filtration $\{\mathcal{F}_t\}$ 
and let  $T_{n}$ be the following discrete approximation of $T$:
\begin{equation}
 \textrm{if $T\in[\frac{k}{2^n},\frac{k+1}{2^n})$}, \ \ \ T_n=\frac{k+1}{2^n}. 
\end{equation}
$T_n$ is a discrete stopping time and the weak Markov property implies that Lemma \ref{Markov_Pair} is also valid for $T_n$. $\{(B_{0}(t+T_n)-B_0(T_n),{_{[1]} r_0}(t+T_n))-B_0(T_n)\}_{t\geq 0}$ converges pathwise to $\{(B_{0}(t+T)-B_0(T),{_{[1]} r_0}(t+T)-B_0(T))\}_{t\geq 0}$ as $n\raw\infty$. 
The result now follows 
from the distributional continuity of 
$(B_0,_{[1]} r_{(\bar x,0)})$ with respect to 
$\bar x$ that we have already established.

\end{proof}

Next, we claim that the pair $(B_0,\rz)$ satisfies the three following properties.
\begin{enumerate}
\item[(1)] $B_0$ is a standard Brownian path starting at $(0,0)$. $\rz$ starts at $z$.
\item[(2)] Away from the diagonal $\{t:_{[1]} r_z(t)=B_0(t)\}$, the two processes evolve as
 two independent Brownian motions.
\item[(3)] Defining $t_\e^r\equiv\inf\{t>0:|\ro-B_0|(t)=\sqrt{2} \e\}$ satisfies
\begin{enumerate}
\item[(i)] $\P\left(\ (\ro-B_0)(t_\e^r)\ = \ +\sqrt{2}\e \right)=1$,
\item[(ii)] $\lim_{\e\daw0} \E(t_\e^r)/\e = \sqrt{2}/\tau$ and $\E([t_\e^r]^2)=o(\e)$ as $\e\daw0$. 
\end{enumerate}
%$B_{stick}=\frac{1}{\sqrt{2}}({_{[1]} r_0}-B_0)$ is a sticky Brownian motion (with an amount of stick $1/\tau$). More precisely, there exists a standard Brownian motion $B$ such that
%\begin{equation}
%\label{equation_of_sticky}
% d  B_{stick}(t)=1_{ B_{stick}(t)\neq 0} d  B(t) + \tau 1_{ B_{stick}(t)=0} dt.
%\end{equation}
\end{enumerate}
In words, (1) and (2) describe the pair $(B_0,\rz)$ away from the diagonal. (3) describes the splitting mechanism when $(B_0,\rz)$ is on the diagonal. 3(i) says that $\rz$ always escapes  the diagonal to the right.
(Note that the definition of $t_\e^r$ given in (3) is consistent
with the one given in Section \ref{proof:sticky} as the first time 
$\ro-B_0$ hits $+\sqrt{2}\e$.)
3(ii) specifies the rate at which $(B_0,\rz)$ escapes the diagonal. We note that this approach is very similar
to the one in \cite{HW06}.

We now turn to the verification of (1)-(3) for $(B_0,\rz)$. Property (1) is obviously satisfied. Property (2) 
follows directly from Lemma \ref{Markov_Pair} and the definition of $\rz$.
Property (3)(i) is obvious. Property 3(ii) is given by Proposition \ref{rate-ex-rpath} above. 

Next, we verify that if $(\bar B_0,\brz)$ is a solution of the SDE (\ref{target::SDE}), it
also satisfies conditions (1)-(3).

\blem
\label{zzzzza}
Let $(\bar B_0,\brz)$ be a solution of the SDE (\ref{target::SDE}). Then
$(\bar B_0,\brz)$ is
a strong Markov process
with stationary transition probabilities and it satisfies conditions (1)-(3). 
\elem
\begin{proof}
As discussed in Section \ref{modify_one_path}, 
the SDE (\ref{target::SDE}) has a unique weak solution
which implies that 
$(\bar B_0,\brz)$ is a 
strong Markov process. The stationarity property is obvious and 
$(\bar B_0,\brz)$ obviously satisfies
properties (1)-(2) and (3)(i). It remains to verify 3(ii).

Since $B_{stick}\equiv(\bro-\bar B_0)/\sqrt{2}$ is 
a $(\sqrt{2}/\tau)$-sticky Brownian motion, it is identical in law with 
$$
t\rightsquigarrow|B|(C(t)), \ \ \ \txt{where $C^{-1}(t)=t+\fr{\sqrt{2}}{\tau} L_0(t)$},
$$  
where $|B|$ is a reflected Brownian motion and $L_0$
is
its local time at the origin. Therefore, $t_\e^r$ (for $(\bro-\bar B_0)$) is distributed like
$$
\frac{\sqrt{2}}{\tau} L_0(T_\e) + T_\e,
$$
where $T_\e$ is the first time $|B|$ hits $\e$. By excursion theory, 
$L_0(T_\e)$ is an exponential random variable with mean $\e$. Since 
the distribution of $T_\e$ is that of $\e^2 T_1$, we indeed get
\begin{equation}
 \E(t_\e^r)/\e \raw \frac{\sqrt{2}}{\tau} \ \ \mbox{and} \ \ \E([t_\e^r]^2)=o(\e). 
\end{equation}

\end{proof}

Finally, we prove the following uniqueness result which is the last ingredient
needed to prove Proposition \ref{cornerstone}. This result is analog to 
Proposition 16 in \cite{HW06}.

\blem \label{laplace_transform} Let
$\{(B_0,{_{[1]}r_z})\}_{z\in\R\times\{0\}}$ and  $\{(\bar
B_0,\brz)\}_{z\in\R\times\{0\}}$ be two families of strong Markov
processes, with stationary transition probabilities, satisfying
properties (1)-(3). For $z=(x,0)$,
$B_{stick,x}\equiv(B_0-\rz)/\sqrt{2}$ and $\bar B_{stick,x}\equiv(\bar
B_0-\brz)/\sqrt{2}$ are equidistributed. 
\elem

\begin{proof}
 By stationarity of the transition probabilities and the 
Markov property, 
$B_{stick,x}$ or $\bar B_{stick,x}$ can be decomposed into two independent parts. The first part is a Brownian motion stopped when it hits zero while the second one is distributed like $B_{stick,0}$ or $\bar B_{stick,0}$.
Hence, it is enough to show that $B_{stick}\equiv B_{stick,0}$ and $\bar B_{stick}\equiv \bar B_{stick,0}$ are equidistributed. 
Also by the Markov property and the stationarity of the transition probabilities,
it is enough to show that 
for any $s\geq0$,  $B_{stick}(s)$ and $\bar B_{stick}(s)$ are equidistributed.
We also note that by property (3)(i), $B_{stick}$ and $\bar B_{stick}$ are $\geq0$.

\bigskip

In the following, $X$ denotes either $B_{stick}$ or $\bar B_{stick}$ and $f$ is a positive bounded continuous function vanishing on the interval $[0,\epsilon_0]$, with $\e_0>0$. For any $\e<\e_0$, define $t^0_\e=0$ and, for any $k\geq0$,
\begin{eqnarray}
t^{2k+1}_\e\equiv\inf\{t>t^{2k}_\e: |X|(t)=\e \}  \ \ , \ \ 
t^{2k+2}_\e\equiv\inf\{t>t^{2k+1}_\e: X(t)=0 \}. 
\end{eqnarray}
We have
\begin{eqnarray}
 \E(\int_0^\infty e^{-\lambda s} f(X(s))ds) & = & \sum_{k=1}^\infty\mathbb{E}\left(\int_{t^{2k-1}_\epsilon}^{t^{2k}_\epsilon} f(X(s)) \ e^{-\lambda s} ds\right), \\
                                 & = & \mathbb{E}\left(\int_{t_{\e}^1}^{t_{\e}^2} f(X(s)) \ e^{-\lambda s} ds \right)  \ \ \sum_{k=0}^\infty \mathbb{E}(e^{-\lambda t^{2k}_\epsilon}).  
\end{eqnarray}
Next, $$t^{2k}_\e=\sum_{i=0}^{k-1} \left([t^{2i+2}_\e-t^{2i+1}_\e]+[t^{2i+1}_\e-t^{2i}_\e]\right).$$
By stationarity and the Markov property, we get that
\begin{equation}
 \mathbb{E}(e^{-\lambda t^{2k}_\epsilon}) = \left(\E(e^{-\lambda t^1_\e})\right)^k \ \left(\E(e^{-\lambda [t^2_\e-t_\e^1]})\right)^{k}.
\end{equation}
This implies that
\begin{eqnarray}
\label{first1}
 \E(\int_{0}^\infty e^{-\lambda s} f(X(s))ds)  =  {\mathbb{E}\left(\int_{t_{\e}^1}^{t_{\e}^2} f(X(s)) \ e^{-\lambda s} ds\right)} 
\frac{1}{1-\E(e^{-\lambda t^1_\e})\dd\E(e^{-\lambda [t^2_\e-t^1_\e]})} .
\end{eqnarray}
Moreover, since
\begin{eqnarray}
\lim_{\e\daw0} \frac{\E(t^r_\e)}{\e}=\frac{\sqrt{2}}{\tau}  \ \ \mbox{and} 
 \ \ \E([t_\e^r]^2)=o(\e)
\end{eqnarray}
and $t_\e^1=t_\e^r$, it follows that
\begin{equation}
\label{second1}
 \E(e^{-\lambda t_\e^1})=1-\frac{\sqrt{2}\lambda}{\tau}\e+o(\e).
\end{equation}

During $[t^1_\e,t^2_\e]$, the process coincides with a Brownian motion starting at $\e$ and stopped when it hits $0$. By standard computations, we get that
\begin{equation}
\label{third1}
 \E(e^{-\lambda [t^2_\e-t^1_\e]})=e^{-\sqrt{2\lambda}\e}.
\end{equation}
Combining Equations (\ref{first1}),(\ref{second1}) and (\ref{third1}), we obtain
\begin{eqnarray}
 \E(\int_{0}^\infty e^{-\lambda s} f(X(s))ds)  & = & \frac{\mathbb{E}\left(\int_{t_{\e}^1}^{t_{\e}^2} f(X(s)) \ e^{-\lambda s} ds\right)}{\e}
\sqrt{2}\left( \sqrt{\lambda}+\frac{\lambda}{\tau}+o(1)\right)^{-1}. \label{benzema}
\end{eqnarray}
Since the left-hand side of the equality does not depend on $\e$,
\begin{eqnarray}
 \frac{\mathbb{E}\left(\int_{t_{\e}^1}^{t_{\e}^2} f(X(s)) \ e^{-\lambda s} ds\right)}{\e}
\end{eqnarray}
has a limit $l(X)$, depending on $f$, as $\epsilon\raw0$ and 
\begin{eqnarray}
 \E(\int_{0}^\infty e^{-\lambda s} f(X(s))ds)  = \int_{0}^\infty e^{-\lambda s} \E(f(X(s)))ds = \sqrt{2} l(X) \left( \sqrt{\lambda}+\frac{\lambda}{\tau}\right)^{-1}.
\end{eqnarray}
Futhermore,
using the various defining properties 
of $B_{stick}$ and $\bar B_{stick}$,
\begin{eqnarray*}
l(X) & = & 
\lim_{\e\daw0} \ \ \e^{-1}\   \E\left(e^{-\lambda t_\e^1} \int_{0}^{t^2_\e-t^1_\e} \ f(X(u+t_\e^1))  \ e^{-\lambda u} \ du \ \right)  \\
& = &
\lim_{\e\daw0} \ \ \e^{-1} \ \E\left(e^{-\lambda t_\e^1}\right) \  \E\left(\int_{0}^{t^2_\e-t^1_\e} \ f(X(u+t_\e^1)) \ e^{-\lambda u} \
du \right) \\
& = &
\lim_{\e\daw0} \ \ \e^{-1}\ \ \E\left(\int_{0}^{T} \ f(\e+B(v))  \ e^{-\lambda v} \ dv\right)\ 
\end{eqnarray*}
where $B$ is a standard Brownian motion and $T$
is the first time it hits $-\e$. (The second equality follows from
the strong Markov property while the third one follows from
$\lim_{\e\daw0} \E\left(e^{-\lambda t_\e^1}\right)=1$
and also from
the fact that on $[t_\e^1,t_\e^2]$, $X$ 
evolves like a Brownian motion.)

Thus $l(B_{stick})=l(\bar B_{stick})$ and therefore
$$\int_{0}^\infty e^{-\lambda s} \E(f(B_{stick}(s)))ds =\int_{0}^\infty e^{-\lambda s} \E(f(\bar B_{stick}(s)))ds.$$
Inverting the Laplace transform yields that for every $s$ and  every positive bounded continuous function $f$ vanishing on the interval $[0,\epsilon_0]$, $\E(f(B_{stick}(s)))=\E(f(\bar B_{stick}(s)))$. By the monotone convergence theorem, we can remove the constraint $f(x)=0$ for  $x\in[0,\epsilon_0]$ which implies that $B_{stick}(s)$ and $\bar B_{stick}(s')$ are equidistributed.

\end{proof}

Lemma \ref{laplace_transform} shows that the
distribution of  $(B_0-\rz)/\sqrt 2$ is
determined by the three properties stated above.  By Lemma \ref{zzzzza}, it follows that $(\ro-B_0)/\sqrt{2}$ is a ($\sqrt{2}/\tau$)-sticky
Brownian motion. The proof of
Proposition \ref{cornerstone} is a consequence of
the following
observation. Let $z=(x,0)$. A pair $(B_0,\rz)$ satisfying properties
(1)-(2) and such that $(\rz-B_0)/\sqrt{2}$ is a $(\sqrt{2}/\tau)$-sticky Brownian motion 
satisfies the SDE (\ref{target::SDE}).

\bigskip

Proposition \ref{cornerstone} being now established, we end this section with a possibly 
surprising theorem about the exit time of a sticky Brownian motion.
Combining Proposition \ref{dist-ter} and Proposition \ref{cornerstone}, we have
\bteo
Let $B_{stick}$ be a Brownian motion starting at $0$ and stickily reflected at $0$ with an amount of stick $\bar \tau$.
If $t_\e$ is the first $\e$ hitting time of $B_{stick}$: 
$$\P(t_\e\leq t) \ = \ \P(l_\e(t)\geq\textrm{Exp}[\frac{1}{2\bar \tau}])$$
where $l_\e$ is the local time at level $\e$ at time $t$ of a Brownian motion on $[0,\e]$, starting at $0$, reflected at $0$ and $\e$, and $\textrm{Exp}\left(1/(2\bar \tau)\right)$ is an independent exponential random variable with mean $1/(2\bar \tau)$.
\eteo

\subsection{The Brownian Net by Marking}\label{pr3}
\label{poof:nett}
The heuristics described in Subsection \ref{equivalence_between_constructions}
are made rigorous in this subsection.

\bigskip

{\bf [Proof of Proposition \ref{Drifting_Brownian}]}

We set $i=0$ (with $z_i=0$) as the proof for general $i$ is essentially the same. 
Recall that $\nro$ and $\{B_k\}_{k\leq n-1}$ are coupled via the SDE (\ref{equation_stick_BM'}).
We start by proving that for such a coupling we have

\blem
\label{stick_on_B}
$
 \forall t\geq 0, \ \  \mathbb{P}(t\in \bigcap_{j=0}^{n-1} \{s: r_0(s)\neq B_{j}(s)\}) \rightarrow 0 \ \ \textrm{as $n\rightarrow\infty$}.
$
\elem

\begin{proof}
Let $\e$ be a fixed positive number. We define
 $$x_n^\e=\sup\{\ B_{j}(t-\e): \ B_{j}(t-\e) \leq r_0(t-\e) \ \txt{for} \ j\leq n-1 \ \}.$$
Let $B^\e$ be  the path in $\{B_{i}\}_{i=0}^{n-1}$ such that $B^\e(t-\e)=x_n^\e$. For any $s\geq t-\e$,
 we define
$$
\D^\e(s)=\frac{1}{\sqrt2}(r_0-B^\e)(s).
$$
By (\ref{equ::needed}), conditioned on the past of $(r_0,B_0,...,B_{n-1})$ up to time $t-\e$, $\D^\e$ solves the following SDE, where $B$ is a standard Brownian motion.
\begin{eqnarray}
d \D^\e(s)=1_{\D^\e\neq 0} d B(s) + \fr{\tau}{\sqrt{2}}  ds  \ \ \ , \ \ \D^\e(t-\e)=x^\e_n.
\end{eqnarray} 
$\D^\e$ is a drifting Brownian motion stickily reflected at $0$ and 
\begin{eqnarray*}
 \P(t\in \bigcap_{j=0}^{n-1} \{s: r_0(s)\neq B_{j}(s)\})
& \leq & \P(r_0(t)\neq B^\e(t))=
\P(\D^\e(t)\neq 0\ | \ {\D^\e(t-\e)=x^\e_n}).
\end{eqnarray*}
Since $x_n^\e\raw0$ as $n\raw\infty$,
%\begin{eqnarray*}
$$
\limsup_{n\raw\infty} \P(t\in \bigcap_{j=0}^{n-1} \{s: r_0(s)\neq B_{j}(s)\})
%& \leq & 
\, \leq \, \P(\D^\e(t)\neq 0\ | \ \D^\e(t-\e)=0)=\P(\D^\e(\e)\neq 0\ | \ \D^\e(0)=0).
$$
%\end{eqnarray*}
Note that the process $\tilde \D^\e$ defined by $d\tilde \D^\e=d \D_\e-\fr{\tau}{\sqrt{2}}1_{\D_\e\neq0}dt$ is a $(\sqrt{2}/\tau)$-sticky
Brownian motion. For such a process, it is known (see e.g., \cite{BS}) that $\P(\tilde \D^\e(\e)\neq 0\ | \ \tilde \D^\e(0)=0)\rightarrow 0$ as $\e\raw0$. By a straightforward application of the Girsanov theorem, we see that
$\P(\D^\e(\e)\neq 0\ | \ \D^\e(0)=0)\rightarrow 0$ as $\e\raw 0$ and Lemma \ref{stick_on_B} follows.
\end{proof}

Let $t>0$ and $\mathbb{P}_{_{[n]}  r_0, B_k}^t$ be the probability measure induced by the pair $(_{[n]}  r_0, B_k)$  on the space of continuous functions on $[0,t]$ endowed with its usual Borel $\s-$ algebra. $\mathbb{P}_{r_0, B_k}^t$ is defined analogously as the distribution of the pair
satisfying (\ref{equ::needed}). We first prove that
\begin{equation}
  \mathbb{P}_{_{[n]}  r_0, B_k}^t\Longrightarrow \mathbb{P}_{ r_0,  B_k}^t \ \  \textrm{as $n\raw\infty$}.
\end{equation}
We define  $_n \chi(t)=1_{t\in \bigcap_{j=0}^{n-1}\{s: r_0(s)\neq  B_j(s)\}}$. 
Lemma \ref{stick_on_B} above and Fubini's Theorem imply that
\begin{equation}
\label{int_to_0}
 \mathbb{E}(\int_{0}^{t} \ _n \chi (t') \ dt') \rightarrow 0.
\end{equation}
For $n\geq k$,  the SDEs (\ref{equation_stick_BM'}) and (\ref{equ::needed}) only differ by their drift term. 
By the Girsanov Theorem, $\P_{_{[n]}  r_0, B_k}^t$ is absolutely continuous with respect to $\mathbb{P}_{ r_0, B_k}^t$ and
\begin{equation}
 d \mathbb{P}_{_{[n]}  r_0, B_k}^t= d \mathbb{P}_{ r_0, B_k }^t\ \exp\left(-\tau \int_{0}^t \  _n \chi(t') \ d r_0(t')+\frac{\tau^2}{2} \int_{0}^t  \ _n \chi(t') \ dt'\right).
\end{equation}
Since $r_0$ is a (drifting) Brownian motion,  (\ref{int_to_0}) and standard arguments imply that the term in the exponential tends to zero in probability. It follows that
\begin{equation}
  \mathbb{P}_{_{[n]}  r_0, B_k}^t\Longrightarrow \mathbb{P}_{ r_0,  B_k}^t \ \  \textrm{as $n\raw\infty$.}
\end{equation}

The pointwise convergence of $\nro$ to $r_0$ was already explained in Section \ref{equivalence_between_constructions} 
by the fact that $\nro$
is monotonic in $n$. This completes the proof of Proposition \ref{Drifting_Brownian}.

\bigskip

{\bf [Proof of Proposition \ref{sticky_pairss}]}

It is easy to see from Proposition \ref{Drifting_Brownian} that $\W_r$ (resp., $\W_l$) is a right-drifting (resp., left-drifting) Brownian web (it is enough to check that two paths in $\W_r$ evolve independently when they are apart; this can been done by simple locality arguments). 
It remains to prove that $\W_r$ and $\W_l$ interact in the
sticky way of a left-right Brownian web (see \cite{SS07} and Subsection \ref{remainder} above). This boils down to proving that $(r_i,l_j)$ satisfies the four-part SDE (\ref{lrsde}). For simplicity, let us take $i=j=0$. Other cases can be treated similarly. We already know that $r_0$ and $l_0$ satisfy
\begin{eqnarray}
d r_0= d B_0^r  \ +  \tau \ dt \\
d l_0= d B_0^l  \ -  \tau \ dt,
\end{eqnarray}
and that $r_0\geq l_0$. It remains to show that 
$d\langle B_0^r,B_0^l\rangle(t)\ = \ 1_{r_0 = l_0}(t) d\langle B_0,B_0\rangle (t)\ = \ 1_{r_0 = l_0}(t) dt.$
As can easily be seen, $r_0$ and $l_0$ evolve independently away from each other. Therefore,
\begin{eqnarray}
 d\langle B_0^r,B_0^l \rangle(t)& = & 1_{r_0 =l_0}(t) d\langle B_0^r,B_0^l\rangle(t) + 1_{r_0 \neq l_0}(t) d\langle B_0^r,B_0^l \rangle(t) \\
 & = &    1_{r_0 = l_0}(t) d \langle B_0^r,B_0^l\rangle(t). \label{ffff} 
\end{eqnarray}
$B_0$ is squeezed between $r_0$ and $l_0$. Hence, $r_0(t)=l_0(t)$ implies that $B_0(t)=r_0(t)=l_0(t)$.
Since, by Proposition \ref{Drifting_Brownian}, 
\begin{eqnarray}
d\langle B_0,B_0^r\rangle(t) &  = &  1_{r_0(t)=B_0(t)} \ dt \\
d\langle B_0,B_0^l\rangle(t) & = &   1_{l_0(t)=B_0(t)} \ dt,
\end{eqnarray}
(\ref{ffff}) implies, as desired, that
\begin{equation}
 d\langle B_0^r,B_0^l\rangle(t)\ = \ 1_{r_0 = l_0}(t) d\langle B_0,B_0\rangle (t)\ = \ 1_{r_0 = l_0}(t) dt.
\end{equation}

\bigskip

{\bf [Proof of Theorem \ref{equivalence_Net}]}

In the proof
we will also consider
$\Net_{wedge}$, the net obtained from
 $(\W_r,\W_l,\hat \W_r,\hat \W_l)$ by the wedge construction
of Subsection \ref{remainder}.
Here $\hat \W_r$ and $\hat \W_l$ are respectively the dual (backward) webs
of $\W_r,\W_l$ (constructed by marking) which can be constructed
using the dual versions of Propositions \ref{Drifting_Brownian} and \ref{sticky_pairss}.

Since $\Net_{hop}=\Net_{wedge}$ (see Theorem \ref{T:netchar}), it suffices to show that (i) $\mathcal{N}_{mark} \supset \mathcal{N}_{hop}$
and (ii) $\mathcal{N}_{mark}\subset\mathcal{N}_{wedge}$.

In order to prove (i), we need to show that a path obtained
by hopping from $\mathcal{W}_r$ to $\mathcal{W}_l$ (or $\mathcal{W}_l$ to $\mathcal{W}_r$)
is still in $\mathcal{N}_{mark}$.
Take two paths $r_i$ and $l_j$ intersecting at time $t$; we need to show that the concatenation of $r_i$ (before $t$) with $l_j$ (after $t$) is in $\mathcal{N}_{mark}$ and similarly for the other concatenation. First,
if we consider the analogous question in a partial net $\mathcal{N}_n$
the result is obviously true. Indeed, if $_{n} r_{i}$ and $_{n} l_{j}$
are respectively the right- and left-most paths of 
$\mathcal{N}_n$ starting from $z_i$ and $z_j$, the path constructed
by hopping from one path to the other at some meeting point is in
$\mathcal{N}_n$. Let $\epsilon>0$ be fixed. Almost surely, there
is some $u\in[t,t+\epsilon]$ such that $r_i(u)>l_j(u)$.  Taking $n$
large enough, we get $_{n}r_{i}(u) > {_{n} l_j}(u)$. On
the other hand, $_n r_i \leq r_i$ and $_n l_j \geq l_j$ so 
that $_n r_i(t)\leq {_{n}l_i}(t)$. Consequently, there exists $v\in[t,u]$
where $_n l_j$ and $_n r_i$ intersect. Now consider the path obtained
by hopping from $_n r_i$ to $_n l_j$ at time $v$. This path is in $\mathcal{N}_n$
and approximates the one obtained by hopping from $r_i$ to $l_j$ at time $t$
except on $[t,t+\epsilon]$. Since $\epsilon$ is arbitrary, the latter path
is approximated by paths in $\bigcup_{n} \mathcal{N}_n$ and therefore it
also belongs to $\mathcal{N}_{mark}$.

We now prove (ii). Consider a wedge constructed from a  pair
$(\hat r_i,\hat l_j)$ starting at $\left((x_i,t),(x_j,t)\right)$
with $x_i<x_j$ and let us assume there exists $\pi\in\mathcal{N}_{mark}$
entering this wedge from outside and show that this
leads to a contradiction. Again, we can approximate $(\hat r_i,\hat l_j)$
by $(_{n}\hat r_i,{_{n}\hat l_j})\in\hat{\mathcal{N}}_n\times\hat{\mathcal{N}}_n$
and $\pi$ by $\pi_n\in\mathcal{N}_n$. Since $_n {\hat r_i} \geq {\hat r_i}$
 and $_n {\hat l_j} \leq {\hat l_j}$, the pair $(_n{\hat r_i},\ _n{\hat l_j})$ forms
a ``partial wedge'' approximating the original wedge from inside.
 Hence, for $n$
large enough, $\pi_n$ would enter this partial wedge from outside. By considering separately
the cases where the putative entering is at a marked $(1,2)$ point of $\mathcal{M}_n$ or not, 
such an entry is seen to be impossible.

\subsection{Separation Points in the Brownian Net}
\label{separation_points}
In Section \ref{modify_the_Web} we defined the
dynamical Brownian web as the limit of partial dynamical webs.
In this subsection, we give a series of results which will guarantee the
existence of such a limit. These are essentially identical to results in \cite{SSS08}. However, in 
\cite{SSS08}
the proofs rely on the hopping construction of the Brownian net, while in this paper we show the results
by using the marking approach. As we shall see, the two 
points of view are rather different.
We start with a definition.

\bdf{\bf [Separation points]}
\label{separation_points_1}
 Two paths $\pi_1$ and $\pi_2$ in $\Net$ starting respectively at $(x_1,t_1)$ and $(x_2,t_2)$ separate at $z=(x,t)$ iff  $t>t_1\vee t_2$ with $\pi_1(t)=\pi_2(t)$ and there exists $a>0$ such that $\pi_1$, $\pi_2$ do not touch on $(t,t+a]$. A point $z$ is called a separation point of $\Net$
iff there is some $\pi_1,\pi_2\in\Net$ that separate at $z$.
\edf

Note that in the partial Brownian net $\mathcal{N}_n$ paths separate at marked $(1,2)$ points. That remains valid in the Brownian net (i.e. when $n\raw\infty$). Indeed, in Subsection \ref{Separation_points_are_marked}, we prove the following result. 

\bprop
\label{separation_point}
The set of separation points in $\mathcal{N}_{mark}$ and the set of marked $(1,2)$ points of the Brownian web coincide.
\eprop
 
Furthermore, in Subsection  \ref{pr:pnew} we prove the following proposition,
which uses the notation $\pi\sim^z B$ and $\pi\sim^z B_{switch}$ introduced in Section \ref{Dynamical::BW}..

\bprop
\label{SP::new}
Let $z=(x,t)$ be a separation point in $\Net_{mark}$, $B$ be any path of $\W$ passing through $z$, and  $\Net_{\leq t-\e}$ be the set of paths in $\Net_{mark}$ starting before or at time $t-\e$. For any $\e\geq0$, define the following (which will
not depend on the choice of $B\in\W$).
\begin{eqnarray*}
\left[\sim^z_\e B_{switch}\right] & = & \{\pi\in\Net_{\leq t-\e} \ : \ \textrm{$\pi$ enters $z$ and $\pi \eq B_{switch}$}\}, \\
\left[\sim^z_\e B \right] & = & \{\pi\in\Net_{\leq t-\e} \ : \ \textrm{$\pi$ enters $z$ and $\pi \eq B$}\}, \\
\nparallel_\e^z & = & \{\pi\in\Net_{\leq t-\e} \ : \ \textrm{$\pi$ does not enter $z$}\}.
\end{eqnarray*}
\begin{enumerate}
\item Let $E_z$ be the set of paths in $\Net_{mark}$ entering $z$. $\sim^z$ is an equivalence relation on $E_z$,
and  $E_z$ can be decomposed into the two equivalence classes $\left[\sim^z_0 B_{switch}\right]$ and $\left[\sim^z_0 B\right]$.
\item For $\e>0$ (note the strict inequality), $\Bsw \ , \  \Bol$ and $\nB$ are disjoint elements of $\mathcal{H}$.  
\item $\exists \bar z\in\R^2$ and $\e>0$ s.t. every path of $\W$
starting in the ball $B(\bar z,\e)$ enters $z$.
\end{enumerate}
\eprop

We note that in the partial net, each path entering a marked point $z$ coincides either with $B$ or $B_{switch}$ for a positive interval of time. In the full net limit, a path coincides either with $B$ or $B_{switch}$ for a positive Lebesgue measure of time.

\subsubsection{Proof of Proposition \ref{separation_point}}
\label{Separation_points_are_marked}

By construction, marked points are separation points so we only need to prove the
converse.

\bdf{\bf [$(T_1,T_2)$ Separation Points]}
\label{Tseparation_points}
 $(x,t)$ with $T_1<t<T_2$ is said to be a $(T_1,T_2)$ separation point iff there are two paths $\pi_1$ and $\pi_2$ in the net starting from $\mathbb{R}\times\{T_1\}$ and separating at $(x,t)$ which do not touch on $(t,T_2]$.
\edf

Let $T_1,T_2$ be two rational numbers. It suffices to prove that if $(x,t)$ is a  $(T_1,T_2)$ separation point of $\Net_{mark}$, then it is a marked $(1,2)$ point. Let
$\pi_1$ and $\pi_2$ be two paths as described in Definition \ref{Tseparation_points}. Since the net is closed under hopping, we can assume without loss of generality that $\pi_1$ and $\pi_2$ have been chosen to coincide up to $t$.

By construction, there exist $\{\pi_i^n\}_{i=1,2}$ with $\pi_i^n$ in the partial net $\mathcal{N}_{n}$($=\Net_{n,n}$; see Subsection \ref{con::net}) so that
$\{\pi_i^n\}$ converges to $\pi_i$. Let us take two numbers $T_1<q_1<q_2\leq t$ where $q_2$ is arbitrarily close to $t$. Proposition \ref{SunSwart} below (for $S=q_1$ and $T=q_2$),
implies that $\pi_1^n(q_2)=\pi_2^n(q_2)$ for large enough $n$.

Hence, for large enough $n$,  $\pi_1^n$ and $\pi_2^n$ start below $\R\times\{q_1\}$ and separate at a point arbitrarily close to $(x,t)$. Since the set of 
$(q_1,T_2)$ separation points is locally finite (see Proposition \ref{(0,1)S.P_locally_finite} below), 
$\pi_1^n$ and $\pi_2^n$ separate at $(x,t)$ for large enough $n$. By construction, $\pi_1^n$ and $\pi_2^n$ only separate at marked points and Proposition \ref{separation_point} follows.

\bprop(\cite{SS07})
\label{SunSwart}
For any $S,T$ with $S<T$, the set of intersection points
 between the line $\R\times T$ and the set paths of $\Net$
 starting on or below $\R\times \{S\}$ is (almost surely) locally finite.
\eprop

\bprop(\cite{SSS08})
\label{(0,1)S.P_locally_finite}
For any $S,T$ with $S<T$, the set of $(S,T)$-separation points is
(almost surely)  locally finite.
\eprop

\subsubsection{Proof of Proposition \ref{SP::new}}
\label{pr:pnew}

\label{entering_1}
\begin{figure}  [!ht]
         \begin{center}
          \includegraphics[width=3cm]{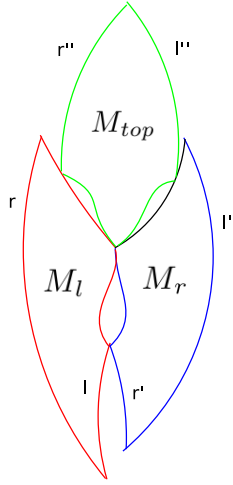}
          \caption{Structure of meshes around a separation point.}\label{mesh.eps}
         %don't know what this does
          \end{center}
          \end{figure}

In the following, for any paths $\pi_1,\pi_2$ in $(\Pi,d)$ entering a point $z$, we will write $\pi_1\sim^z_{out} \pi_2$ (resp., $\pi_1\sim^z_{in} \pi_2$) iff for any $\e>0$,
$$\int_{t}^{t+\epsilon} 1_{\pi_1(u)=\pi_2(u)} \ du >0 \ (\txt{resp.,}  \int_{t-\e}^{t} 1_{\pi_1(u)=\pi_2(u)} du>0).$$
Note that $\pi_1\sim^z\pi_2$ iff $\pi_1\sim^z_{out} \pi_2$ and $\pi_1\sim^z_{in} \pi_2$.
In order to prove Proposition \ref{SP::new}, we will use the following result from \cite{SSS08}. 
Since this result is part of a much larger theorem there, we provide a direct proof.
For a ``pictorial'' representation of the result, see Figure \ref{mesh.eps}.

\bteo(\cite{SSS08})
\label{str::meshes}
Let $z=(x,t)$ be a separation point in $\Net$ and let $\e>0$. There exist three distinct meshes $M_l(r,l)$, $M_r(r',l')$ and $M_{top}(r'',l'')$ such that
\begin{enumerate}
\item  The bottom times of $M_l(r,l)$, $M_r(r',l')$ are in $(t-\e,t)$ and their top times are in $(t,\infty)$.
Moreover,  $l\leq r'$ (at coexistence times of $M_r$ and $M_l$), $l(t)=r'(t)=x$ and $l\sim^z_{in} r'$.
\item $z$ is the bottom point of $M_{top}(r'',l'')$. $M_{top}(r'',l'')$ is squeezed between $M_l(r,l)$ and $M_r(r',l')$ (i.e., $l\leq r''$ and $l''\leq r'$ at respective 
coexistence times). Moreover, $r''\sim^z_{out} l$, $r'\sim^z_{out} l''$.
\end{enumerate}
\eteo

\begin{proof}

In the following, we say that two paths $\pi_1$ and $\pi_2$ meet at time $\bar t$ iff $\pi_1(\bar t)=\pi_2(\bar t)$ but $\pi_1<\pi_2$ or $\pi_1>\pi_2$ on $(\bar t-a,\bar t)$ for some $a>0$.

{\bf [Construction of $M_r$ and $M_l$]}
 
Recall that we constructed the Brownian net by marking a non-drifting Brownian web. 
There is an alternative marking construction of the net which can be described as follows. 
Start with a {\it left-drifting} Brownian web $\mathcal{W}_l$,
 with drift $-\tau$. 
Mark the $(1,2)$ points of $\mathcal{W}_l$ and construct $\mathcal{N}_{l}$
 by branching at all the left $(1,2)$ points (of $\W_l$) in $\mathcal{M}_l(2\tau)$, the set of marks whose dynamical 
time coordinate is $\leq 2\tau$ (the factor $2$ compensates for 
the $-\tau$ drift in $\W_l$). On the one hand, 
repeating step by step what was
 done in  Section \ref{Brownian_Net} (see Theorem \ref{equivalence_Net}), 
one can show that $\mathcal{N}_{l}$ is identical in law to the usual $\mathcal{N}_{hop}$, as in Subsection \ref{remainder}.
On the other hand, following the proof of Proposition \ref{separation_point},  separation points  of the net $\mathcal{N}_{l}$ must be marked left  $(1,2)$ points of $\mathcal{W}_l$. Hence, separation points of the net $\mathcal{N}_{hop}$ are left $(1,2)$ points of $\mathcal{W}_l$ and symmetrically they are also right $(1,2)$ points of $\mathcal{W}_r$.

One consequence is that $z=(x,t)$ must be a separation point for two paths
$\bar l\in\W_l$ and $\bar r\in\W_r$ starting from deterministic points.
Lemma 6.5 in \cite{SS07} analyzes 
meshes to the left of a path $\bar l\in\W_l$.
Using that lemma and the fact that points on $\bar l$ where other paths 
from $\W_l$ coalesce with $\bar l$ from the left are dense in $\bar l$ (along with the analogous results
for $\bar r$), it follows that
there exists 
a mesh
$M_l(r,l)$ 
(resp., $M_r(r',l')$) with bottom time in $(t-\e,t)$ and top time $>t$ such that $l(t)=\bar l(t)=x$ (resp., $r'(t)=\bar r(t)=x$). 

By Corollary \ref{one_coalescence_time_2},
 $l$ and $r'$ coalesce with some paths 
$l_i$ and $r_j$ (in the skeleton of $\W_l$ and $\W_r$ respectively) before entering the point $z$.
The pair $(l_i,r_j)$ satisfies the SDE (\ref{lrsde}) and
in particular,  $l_i \leq r_j$ from the first time they meet.
It is clear that $l_i$ and $r_j$ do not meet and separate at the same point. Hence,
there exists $a'>0$ so that $ l_i \leq  r_j$ on $[t-a',\infty)$ and a sequence $t_n\uaw t$ s.t. $l_i(t_n)=r_j(t_n)$. It immediately follows that there exists $a''$ with $a'\geq a''>0$ so that
 $ l \leq  r'$ on $[t-a'',\infty)$ and a sequence $t_n'\uaw t$ s.t. $l(t_n')=r'(t_n')$.
Lemma \ref{touch::stickyness} below then immediately implies that $l\sim^z_{in} r'$.

{\bf [Construction of $M_{top}$]}

Up to reversal of the time coordinate, the backward Brownian net is distributed as the Brownian net (see Subsection \ref{remainder}). Hence, by what has been just proved, $z$ is a separation point for two paths $(\hat l, \hat r)\in(\hat{\mathcal{W}}_l,\hat{\mathcal{W}}_r)$ and there exists $a>0$ such that $\hat r \leq \hat l$ on $(-\infty,t+a]$. Let $r''$ (resp., $l''$) be the newly born path of $\W_r$ (resp. $\W_l$) starting from $z$. Since $(x,t)$ is a right $(1,2)$ point for $\mathcal{W}_r$ and  a left $(1,2)$ point for $\mathcal{W}_l$, we get that on $(t,t+a]$
\begin{equation}
\label{not_dollar}
 r''\leq \hat r \leq \hat l \leq l'' \ \ \ \txt{and} \ \ r''\leq r', \ l \leq l''.
\end{equation}
$M_{top}$ is defined as the mesh $M_{top}(r'',l'')$ formed by $r''$ and $l''$.

The second part of (\ref{not_dollar}) implies that  $M_{top}$ is either squeezed between $M_r$ and $M_l$ or it contains either $l$ or $r'$. Since
paths of $\Net$ do not enter meshes from outside, we get that on $(t,t+a]$
\begin{equation}
\label{not_dollar_2}
  l \leq r''\leq \hat r \leq \hat l \leq l'' \leq r'. 
\end{equation}
Recall the construction of the net $\mathcal{N}_{l}$ (described at the beginning of this proof) based on the marking of a left drifting Brownian web    
 and let $(l,\hat l)$ be a pair of paths in $(\W_l,\hat \W_l)$. As can be easily seen, the set
$$
\{(x,t) : l(t)=\hat l(t)=x \ \txt{and} \  \exists a>0 \ \txt{s.t.} \ \forall s\in(t,t+a), \ l(s)<\hat l(s) \}
$$
has zero local time measure. (By Proposition \ref{forwardandbackward2} and taking the difference between $l$ and $\hat l$, this follows from the fact that, for a standard Brownian motion $B$, the set
$$
\{t : B(t)=0 \ \txt{and} \  \exists a>0 \ \txt{s.t.} \ \forall s\in(t,t+a), \ |B(s)|>0 \}
$$
has zero local time measure). Since in $\Net_l$, separation points are left marked $(1,2)$ points, the  argument just given implies
that for every {\it marked} point, there exists $t_n\daw t$ such that $l(t_n)=\hat l(t_n)$. By (\ref{not_dollar_2}), $l\leq r''\leq \hat l$, implying that $l(t_n)=r''(t_n)$. By Lemma \ref{touch::stickyness} below we have
that $l\sim^z_{out} r''$ and by a similar argument, we get $r'\sim^z_{out} l''$. 
\end{proof}
\blem
\label{touch::stickyness}
Let $(l,r)\in(\W_l,\W_r)$ be 
such that for some $t>t_{r}\vee t_{l}$, $l(t)=r(t)$. 
For any $\e>0$, $\int_{t-\e}^{t+\e} 1_{l(s)=r(s)} ds>0$.
\elem
\begin{proof}
Choose any $t'$ with $t_{r}\vee t_{l}<t'<t-\e$. By Corollary \ref{one_coalescence_time_2}, 
on $[t',\infty)$, the pair
$(l,r)$
coincides with a pair $(L,R)$ of $(\W_l,\W_r)$, starting from deterministic points and satisfying the SDE (\ref{lrsde}). 
Lemma \ref{touch::stickyness} then follows from the fact (see Proposition 3.1 in \cite{SS07})
that the support of the measure $\mu$, defined as $\mu([t_1,t_2])=|\{t\in[t_1,t_2]: L(t)=R(t)\}|$, coincides with $\{t: L(t)=R(t)\}$.
\end{proof}

We now prove the first two claims of Proposition \ref{SP::new} for a  separation point $z=(x,t)$. Note that if claim 2 holds for a given $\e$,
it immediately holds for any $\e'>\e$. Hence, w.l.o.g.,  we can take $\e>0$ small enough such that there is a
path $B\in\W$
entering $z$ and starting
at $t'\leq t-\e$. In the following, $\tilde E_\e$ will denote 
the subset of $\Net_{\leq t-\e}$ consisting of all the paths entering $z$.

Recall that paths of $\mathcal{N}$ do not enter meshes 
(see Theorem \ref{T:netchar} (b3) in Subsection \ref{remainder}). Hence, for any given mesh $M$ with bottom time in $(t-\e,\infty)$, we can partition $\Net_{\leq t-\e}$
into $\{R(M),L(M)\}$, where $R(M)$ (resp., $L(M)$) is the compact subset of $\Net_{\leq t-\e}$ consisting of all
the paths passing to the right (resp., left) of $M$. Let $M_r,M_l$ and $M_{top}$ be 
as in Theorem \ref{str::meshes} and let us define 
\begin{eqnarray}
 \tilde E^r_\e=[L(M_r)\cap R(M_l)]\cap R(M_{top}) \ \ \ , \ \  \tilde E^l_\e=[L(M_r)\cap R(M_l)]\cap L(M_{top}),
\\
\tilde E^{c}_\e=R(M_r)\cup L(M_l).
\end{eqnarray}
In particular $\{\tilde E^r_\e,\tilde E^l_\e\}$ (resp., $\{\tilde E^r_\e,\tilde E^l_\e,\tilde E^c_\e\}$) defines a natural partition of $\tilde E_\e$ (resp., $\Net_{\leq t-\e}$) into elements of $\mathcal{H}$. 

By definition, paths in $\tilde E^{l}_\e$ are squeezed between $l$ and $r'$ below $z$ while they are squeezed between $l$ and $r''$ above $z$.
Hence, Theorem \ref{str::meshes}
immediately implies that for any two paths $\pi_1,\pi_2\in \tilde E^{l}_\e$, $\pi_1\sim^z\pi_2$. The same property holds for $\tilde E^{r}_\e$. 
Conversely, if $\pi_l\in \tilde E^{l}_\e$ and 
$\pi_r\in \tilde E^{r}_\e$, the two paths 
separate at $z$. 
This implies that $\sim^z$ is an equivalence
relation on $\tilde E_\e$ and the corresponding 
equivalence classes are given by
$\tilde E^{r}_\e$ and $\tilde E^{l}_\e$. 
Since $B$ and
$B_{switch}$ separate at $z$   
they do not belong to the
same equivalence class and claims 1 and 2 of
Proposition \ref{SP::new} follow.

Next, we say that the ball $B(\bar z,\bar \e)$ with $\bar z=(\bar x,\bar t)$ is squeezed between
between $l$ and $r'$ iff $\bar t-\bar\e \geq t_{l}\vee t_{r'}$ and for
every $(x',t')\in B(\bar z,\bar \e)$, $l(t')\leq x' \leq r'(t')$. It is clear that one can find such a ball below the point $z$ and that
any path starting from that ball
is squeezed between $l$ and $r'$ and so
is forced to enter the 
point $z$. Claim 3 of Proposition of \ref{SP::new} follows.

\subsection{The Dynamical Brownian web}
\label{rrrrrrrr}

\subsubsection{Proof of Proposition \ref{existence_of_Br}}
\label{edep}
In the following, we use the notation of Proposition \ref{SP::new}.

By compactness of $\Net(\tau)$, $\{\W_{(n,m)}(\tau)\}_{(n,m)}$
is a precompact subset of $\mathcal{H}$. 
Let  $\W_{1}$
be any subsequential limit of $\{\W_{(n,m)}(\tau)\}_{(n,m)}$
as 
$n,m \rightarrow \infty$  and let 
$$\W_{2}(\tau)=\{\pi\in\Net_{mark}(\tau) \ : \ \textrm{every time $\pi$ enters a point $z$ in $\M({\tau})$, $\pi \sim^z B_{switch}$}  \}.$$
be as in item (2) of Proposition \ref{existence_of_Br}. We first prove

\bigskip

{\bf (i) $\W_{1}(\tau) \subset \W_{2}(\tau)$.}

Let $z=(x,t)\in\M(\tau)$, and let $\pi\in\W_{1}(\tau)$ start at $t-2\e$ with $\e>0$, and pass through $z$. 
By definition, there exists a sequence $\{\pi_N\}_{N\geq0}$ so that $\pi_{N}$ belongs to $\cup_{n,m>N}\W_{n,m}(\tau)$ and $\{\pi_N\}$ converges to $\pi$.
Taking $N$ large enough, we can assume w.l.o.g. that $\pi_N$ belongs to $\Net_{\leq t-\e}$ and $(x,t)\in\mathcal{M}_{n,m}(\tau)$  
for $n,m>N$.
By Proposition \ref{SP::new}(2),
$\{\pi_N\}$ enters $z$ for $N$ large enough. By construction, $\pi_N \sim^z B_{switch}$ and since 
$\pi_N\raw\pi$, Proposition \ref{SP::new}(2) implies that
$\pi\sim^z B_{switch}$. Hence, $\W_{1}(\tau) \subset \W_{2}(\tau)$.

\bigskip

Next, we prove that $\W_{2}(\tau)$ satisfies (3)(o). We first claim that when two paths of
$\W_{2}(\tau)$ meet, they coalesce. Let $\pi_1,\pi_2\in\W_{2}(\tau)$ start at $t_1,t_2$ respectively and meet at $t'>t_1\vee t_2$ and let us assume that $\pi_1$ and $\pi_2$ separate at $z=(x,t)$ with $t\geq t'$. By Proposition \ref{SP::new}(1), either $\pi_1\sim^z B$ or $\pi_2\sim^z B$.
This contradicts the definition of $\W_2(\tau)$ and we conclude 
that $\W_{2}(\tau)$ is a coalescing collection of paths. Let $z_i\in\mathcal{D}$. Any path in $\W_{2}(\tau)$ starting at $z_i$ is squeezed between $r_{i}$ and $l_{i}$, the paths in $\W_r$ and $\W_l$
respectively starting from $z_i=(x_i,t_i)$. Since there exists a sequence $t_n'\daw t_i$ s.t. $l_{i}(t_n')=r_{i}(t_n')$ and since paths in $\W_{2}(\tau)$ coalesce, there must 
be a unique path in $\W_{2}(\tau)$ starting from $z_i$.
We call this path $B_i^\tau$ and define $\W_{3}(\tau)$ as $\overline{\{ {B}_{i}^\tau\}}$.
We continue to prove:

\bigskip

{\bf (ii) $\W_{2}(\tau) \subset \W_{3}(\tau)$.}

Let $\pi\in\W_{2}(\tau)$ start at $(x',t')$ and let $\e>0$. We claim that $\pi$ hits a path in $\W_r\cup\W_l$ in $(t',t'+\e]$. To see this, let $a\in(t',t'+\e)$ and let $\{r_n\}_n\subset\W_r$ (resp., $\{l_n\}_n\subset\W_r$) start at $z_n^r$ (resp., $z_n^l$) with $z_n^r$ (resp., $z_n^l$) converging to $(\pi(a),a)$ from the left (resp., from the right) of $\pi$.
If there is not any path in $\{r_n,l_n\}$  meeting $\pi$ on $(a,t'+\e)$, $\{r_n\}$ and $\{l_n\}$ converge (along a subsequence) to $r\in\W_r$ and $l\in\W_l$ respectively, both starting at $(\pi(a),a)$ and s.t. $r<\pi<l$ on $(a,t'+\e)$. In other words, $\pi$ enters a mesh from outside, yielding a contradiction to Theorem \ref{T:netchar}.

$\M(\tau)$, or equivalently the set of separation points in $\Net(\tau)$, is dense along any path $\pi'$ in $\W_r\cup\W_l$. Since once $\pi$ touches some $\pi'$, they can only separate at a point
in $\M(\tau)$, it follows that 
$\pi$ enters some point $z\in\M(\tau)$ before $t+\e$. By virtue of Proposition \ref{SP::new} (3), there exists a ball $B(\bar z,\e')$ such that any path in $\Net(\tau)$ starting in
$B(\bar z,\e')$ enters the point $z$. Hence, any path
$B^{\tau}_{i}$ such that $z_i$ belongs to $\mathcal{D}\cap B(\bar z,\e')$ hits $z$. It follows that $\pi$ coalesces with some $B^{\tau}_{i}$ before time $t'+\e$. As a consequence, $\W_{2}(\tau) \subset \W_{3}(\tau)$.
Finally, we prove:

\bigskip

{\bf (iii) $\W_{3}(\tau) \subset \W_{1}(\tau)$.}

It is clear that there is at least one path $\pi_i\in\W_1(\tau)$ starting from $z_i$. Since $\W_{1}(\tau)\subset \W_2(\tau)$, property 3(o) for $\W_2(\tau)$ (which we have already proved) implies that $\pi_i=B_{i}^\tau$. Since $\W_{1}(\tau)$ is compact, it follows that $\W_{3}(\tau)\subset\W_{1}(\tau)$ and from (i), (ii) above, we get that $\W_{1}(\tau)=\W_2(\tau)=\W_3(\tau)$. This shows
that all subsequence limits of
$\{\W_{n,m}\}$ agree and Proposition \ref{existence_of_Br} follows.

\subsubsection{$(\W,\W({\tau}))$ is a $1/(2\tau)$-Sticky Pair of Brownian Webs}
\label{proof::stickiness::webs}
In the remaining subsections of the paper we prove the four parts
of Theorem \ref{real::stick}.
In this subsection and the next, the term marking will refer to the set $\mathcal{M}(\tau)$.
We already showed in the proof of Proposition \ref{existence_of_Br} that $\W(\tau)$ is a coalescing
set of paths. By
a simple locality argument, it is not hard to see that for $i\neq j$, $B^\tau_i$ and $B^\tau_j$ move independently when they are apart. In the following, we 
prove that $(B_i,B_j^\tau)$ is a $1/(2\tau)$-sticky pair of Brownian motions. This ensures that 
each $B_j^\tau$ is a Brownian motion and since the paths 
of $\W(\tau)$ are coalescing, it follows that
$\W(\tau)$ is a Brownian web and furthermore that the interaction
between $\W$ and $\W(\tau)$ is $(1/2\tau)$-sticky as claimed.

\bigskip

We now prove that $(B_i,B_j^\tau)$ is a $1/(2\tau)$-sticky pair of Brownian motions.
Since the distribution of the Brownian net is invariant
under translation in the space time domain, 
Proposition \ref{existence_of_Br}(2) implies that
$\W(\tau)$ is also translation invariant. Hence, it
suffices to prove that
$(B_0,B_j^\tau)$ is a $1/(2\tau)$-sticky pair of Brownian motions.

Define ${_{(n,m)}B_j^\tau}$
as the path obtained from $B_j$ after switching the directions of the points in $\mathcal{M}_{(n,m)}(\tau)$.
By parts (1) and (3)-(o) of Proposition \ref{existence_of_Br}, we have
\begin{equation}
 \lim_{n\uaw\infty}\lim_{m\uaw\infty} d({_{(n,m)}B_j^\tau},B_j^\tau)=0.
\end{equation}
In the following, we will denote by ${_{[n]}B_j}\equiv {_{[n]}B_j^\tau}$ the limit of ${_{(n,m)}B_j^\tau}$ as $m\raw\infty$. Informally, ${_{[n]}B_j}$ is the path constructed from $B_j$ after switching the direction of all the (left and right) $(1,2)$ points in $\mathcal{M}(\tau)$ that lie on $\{B_i\}_{i=0}^{n-1}$.

In order to prove that $(B_0,B_j^\tau)$ is a $1/(2\tau)$-sticky
pair of Brownian motions, we 
claim that it is enough to prove the following lemma 
(which is done in Subsection \ref{Lemma::esdd} below).

\blem
\label{intermed}
 $(B_0,{_{[1]}B_j})$ is  $1/(2\tau)$-sticky pair of Brownian motions.
\elem
The sufficiency of Lemma \ref{intermed} follows from the observation that the law of $(B_0,{_{[n]}B_j})$ is identical to the one of $(B_0,{_{[1]}B_j})$. 
For example, for $n=2$, one may consider a revised marked Brownian web $\W^{*}$ 
in which all the marked $(1,2)$ points along the finite segment of $B_1$ before it
coalesces with $B_0$ have been switched. In $\W^*$ the marks along $B_0^*(\equiv B_0)$ are the same as in the original web. The following lemma (for $k=1$ and $l=0$) 
states that this $\W^*$ is equidistributed with the original
Brownian web. On the other hand, the pair $(B_0^*,{_{[1]} B_j^{*}})$ for $\W^*$ is { {\it identical}
to the pair $(B_0,{_{[2]}B_j})$ for the original marked web. Since ${_{[n]}B_j}$  almost surely converges to $B^\tau_j$, 
$(B_0,B_j^\tau)$ is  $1/(2\tau)$-sticky pair of Brownian motions.

\blem
Let $_{[1]}\W$ denote the web resulting from switching all the marked $(1,2)$
points in the original web $\W$ along $B_0$;
then $_{[1]}\W$ is equidistributed as the original web.
Similarly, if for some fixed $k,l$ with $k\neq l$,
$\W^*$ denotes the marked web
resulting from switching the original web along 
the finite
segment
of $B_k$ before it coalesces
with $B_l$. Then $\W^*$ is equidistributed with $\W$.
\elem
\begin{proof}
To prove the first part of the lemma, it suffices
to show that $\{_{[1]}B_j\}$ are coalescing Brownian motions.
Lemma \ref{intermed} implies that each individual $_{[1]}B_j$
is a Brownian motion and their construction shows that they are independent
before meeting.
The proof that they coalesce upon meeting is basically the same as that 
given for the paths of $\W_2(\tau)$ in Subsection \ref{edep}.
For the second part of the lemma
w.l.o.g., set $k=0$. Then the paths $B_j^*\in\W^*$
starting from $z_j$ coincide with $_{[1]}B_j$ for times
before the coalescence time of $B_0$ and $B_l$ and 
afterward coincide with paths in $\W$. 
It follows
that $\{B_j^*\}$ are coalescing Brownian motions and thus that
$\W^*$ is equidistributed with $\W$.
\end{proof}

\subsubsection{Proof of Lemma \ref{intermed}}
\label{Lemma::esdd}
We prove the result for $j=0$. The result can then
be trivially extended to any $j$.
Our proof follows along
the lines of the proof of
Proposition \ref{cornerstone} given in Subsection \ref{proof:sticky}, except of course that
here both right and left marked $(1,2)$ points along $B_0$ are switched leading to
${_{[1]}B_0}$ rather than ${_{[1]}r_0}$. 
Here, it is enough to prove that $(B_0,\Bzz)_{z\in\R\times\{0\}}$ is a family of strong Markov processes
with stationary transition probabilities and that
the pair $(B_0,\Bzz)$ satisfies the following three  properties.
\begin{enumerate}
\item[(1)] $B_0$ is a standard Brownian path starting at $(0,0)$. $\Bzz$ starts at $z$. 
\item[(2)] Away from the diagonal $\{t:\Bzz(t)=B_0(t)\}$, the two processes evolve as
 two independent Brownian motions.
\item[(3)] Defining $t_\e=\inf\{t>0:|\Bo-B_0|(t)=\sqrt{2} \e\}$, one has
\begin{enumerate}
\item[(i)] $\P\left(\ (\Bo-B_0)(t_\e)\ = \ \sqrt{2}\e \right)=\frac{1}{2}$
\item[(ii)] $\lim_{\e\daw0} \E(t_\e)/\e = \sqrt{2}/(2\tau)$ and $\E([t_\e]^2)=o(\e)$. 
\end{enumerate}
\end{enumerate}
The strong Markov property and the stationarity of the transition probabilities
can be shown as
in Lemma \ref{Markov_Pair}. 
Those two properties and the definition of $\Bzz$ 
easily imply Properties (1)-(2). Property (3)-(i) is clearly true by right-left symmetry.
It remains to prove (3)-(ii). Recall the definition of 
$\ro$ given in Section \ref{modify_one_path}. We define ${_{[1]}l_0}$ analogously, i.e., ${_{[1]}l_0}$ is obtained from $B_0$ by switching all the marked right $(1,2)$ points in $\M(\tau)\cap B_0$. We also define  
\begin{eqnarray}
t_\e^l=\inf\{t: {_{[1]}l_0}(t)=B_0(t)-\sqrt{2}\e\}, \\
t_\e^r=\inf\{t: {_{[1]}r_0}(t)=B_0(t)+\sqrt{2}\e\}.
\end{eqnarray}

$t_\e^r$, which was carefully studied in Subsection \ref{sec-excursions}, (resp., $t_\e^l$) is the first time a right (resp., left) marked excursion away from $B_0$ hits $B_0 + \sqrt{2}\e$ (resp., $B_0 -\sqrt{2}\e$). 
In order to verify the first part of (3)-(ii), we will
prove that $\lim_{\e\daw0}\E(t_\e)/\e$ coincides with
$\lim_{\e\daw0}\E(t_\e^r\wedge t_\e^l)/\e$ and that $\E(t_\e^r\wedge t_\e^l)/\e$ has the desired limit. 
The second part can be proved similarly. 
We first use the following lemma.

\blem
\label{excursion::bo}
$\Bo$ is obtained by joining together marked excursions from $B_0$.
\elem

\begin{proof}
Let $z$ be a point at which $\Bo$ separates from $B_0$. By Proposition \ref{separation_point}, $z$ is a marked point of the original Brownian web $\W$ 
and there is a marked excursion $e$ from $B_0$ starting at $z$.
By the
structure of the separation points given in Proposition \ref{SP::new}
and since  
$_{(1,n)}B^\tau_0\raw \Bo$ as $n\raw\infty$, we see that
$_{(1,n)}B^\tau_0$ follows the excursion $e$ for sufficiently large $n$. As a consequence, $\Bo$ also follows $e$. Since this is true for every such $z$, the lemma follows.
\end{proof}

Lemma \ref{excursion::bo} immediately implies that 
\begin{eqnarray}
\label{ref::3333}
 t_\e \geq t_\e^r \wedge t_\e^l. 
\end{eqnarray}
Continuing with our proof of Property (3)-(ii), we define
$$ 
T_\e=\inf\{t \geq t_\e^r\wedge t_\e^l:  |\Bo(t)-B_0(t)|= 0 \}.
$$
Using $T_\e^{(0)}\equiv T_\e$ as a (first) stopping time increment, 
%considering
%the segments $\{B_0^{(0)},\ro^{(0)},\lo^{(0)},\Bo^{(0)}\}$ 
denoting the segments 
of $\{B_0,\ro,\lo,\Bo\}$ up to time
$T_\e^{(0)}$ 
by $\{B_0^{(0)},\ro^{(0)},\lo^{(0)},\Bo^{(0)}\}$
and then translating $\left(B_0(T_\e),T_\e\right)$ onto $(0,0)$, 
we may inductively define 
$$
\{B_0^{(n)},\ro^{(n)},\lo^{(n)},\Bo^{(n)},T_\e^{(n)}\},
$$
which, as in the proof of Lemma \ref{Markov_Pair}, are i.i.d. Next, define
\begin{eqnarray}
K_\e=\inf\{ k: \  \exists t\in[0, T_\e^{(k)}], \ \ \ |\Bo^{(k)}-B_0^{(k)}|(t)=\sqrt{2} \e\}
\end{eqnarray}
and also, $$\tilde T_\e^{(n)}=T_\e^{(n)}\wedge\inf\{t\in[0,T_{\e}^{(n)}]: |\Bo^{(n)}(t)-B_0^{(n)}(t)|= \sqrt{2}\e \}.$$
Then, (letting $\tilde T_\e\equiv \tilde T_\e^{(0)}$) we have
\begin{eqnarray}
\E(t_\e)   &  = &\ \sum_{n\geq 0} \ \E(\tilde T^{(n)}_\e \ 1_{K_{\e}\geq n}) 
  =  \ \sum_{n\geq 0} \ \E(\tilde T_\e) \ \P(K_{\e}\geq n) \\
 & = &   \E( \tilde T_\e)  \ \sum_{n\geq0} \P\left(\forall t\in[0,T_\e], \ \ \ |\Bz-B_0|(t)<\sqrt{2} \e\right)^n \\
& = & \ \frac{\E(\tilde T_\e)}{ \P\left(\exists t\in[0,T_\e], \ \ \ |\Bz-B_0|(t)=\sqrt{2} \e\right)} \\
& \leq & \ \frac{\E(\tilde T_\e)}{ \P\left(t_\e^l \wedge t_\e^r = t_\e\right)}, \label{wq1}
\end{eqnarray}

Next we prove the following  three lemmas.

\blem
\label{1:0}
$\E(\tilde T_\e-t_\e^r \wedge t_\e^l)/\e\raw 0$ as $\e\daw0$.
\elem
\begin{proof}
The path $\Bz$ evolves like a Brownian motion when it is away from $B_0$.
It follows that  $\E(\tilde T_\e-t_\e^r \wedge t_\e^l)\leq \sup_{x\in[0,\e]}(\E(S_x))$ where $S_x$ is the time a standard Brownian motion starting at $x$ exits the interval $[0,\e]$. This yields 
the claimed result.
\end{proof}

\blem
\label{1:2}
$
 \E(t_\e^r\wedge t_\e^l)/\e \rightarrow \sqrt{2}/(2\tau).
$
\elem

\begin{proof}
Conditioned on $\W$ (but not the marking $\M(\tau)$), $t_\e^r$ and $t_\e^l$ are independent. If we denote by  $\P_\W$ the probability distribution of the marked Brownian web conditioned on a realization of the web $\W$, and by $\E$ expectation with respect to the distribution $\P$ of $\W$, we have 
\begin{eqnarray}
 \E(t_\e^r\wedge t_\e^l)/\e & = &  \int_{0}^\infty \E(\P_{\W}(t_\e^r\wedge t_\e^l \geq \e t)) \ dt=
\int_{0}^\infty \E(\P_{\W}(t_\e^r \geq \e t) \cdot \P_{\W}(t_\e^l \geq \e t))  \ dt. 
\end{eqnarray}
By Proposition \ref{dist-ter},
\begin{equation}
 \P_\W(t_\e^r \geq \e t)=\P_\W\left( L_{\e,\e t}([0,\e t])\leq \textrm{Exp}(1/(\sqrt{2}\tau))\right)=\exp(-\sqrt{2}\tau l_{\e}(\e t)),
\end{equation}
By Lemma \ref{conv::LT},
we know that in probability
$l_{\e}(\e t) \raw t/2$,
implying that $\P_{\W}(t_\e^r \geq \e t)\raw e^{-\tau t/\sqrt{2}}$. By symmetry, $\P_{\W}(t_\e^l \geq \e t)\raw e^{-\tau t/\sqrt{2}}$. 
In Subsection \ref{distr-tere2}, we showed that $\{\P(t_\e^r\geq \e \cdot)=\E(\P_{\W}(t_\e^r\geq \e \cdot))\}_{\e\leq 1}$ is uniformly integrable. Since
$$
\E(\P_\W(t_\e^r \geq \e t) \P_{\W}(t_\e^l \geq \e t)) \leq \E(\P_{\W}(t_\e^r \geq \e t))=
\P(t_\e^r \geq \e t),
$$
we have that
\begin{equation}
\lim_{\e\downarrow 0} \E(t_\e^r\wedge t_\e^l)/\e  \ dt=  
\int_{0}^\infty \lim_{\e\downarrow 0} \E( \P_{\W}(t_\e^r \geq \e t)\cdot \P_{\W}(t_\e^r \geq \e t) ) dt = \int_{0}^\infty \exp(-2t \tau/\sqrt{2})=\frac{\sqrt{2}}{2\tau}.
\end{equation}
\end{proof}

\blem
\label{1:1}
$\lim_{\e\daw0}\P\left(t_\e^l \wedge t_\e^r\neq t_\e \right)=0.$
\elem

\begin{proof}
By symmetry, it suffices to prove that
\beq
\lim_{\e\daw0}\P\left(t_\e^l=t_\e^l \wedge t_\e^r, \ t_\e^l \neq t_\e \right)=0.
\eeq
Assume that 
$t_\e^l=t_\e^r \wedge t_\e^l$ and $t_\e \neq t_\e^l$. 
Then there exists a left marked excursion 
$e_{l,\e}$  from $B_0$ starting at $T(e_{l,\e})$
which is at a 
distance $\sqrt{2}\e$ from $B_0$ at time $t_\e^l$. 
Since $t_\e \neq t_\e^l$, $\Bo$  avoids this
excursion implying that 
$\Bo$ follows a right marked excursion $e_r$ during 
the time interval $[T(e_r),T(e_r)+D(e_r)]$
such that $T(e_r)< T(e_{l,\e})< T(e_r)+D(e_r)$. In other words,
$T(e_{l,\e})$ is straddled by a marked right excursion.
The lemma follows from Proposition
\ref{rihght-left-e}.
\end{proof}

%Combining (\ref{wq1}) and Lemmas \ref{1:0}, \ref{1:2} and \ref{1:1}, we have that
By (\ref{wq1}) and Lemmas \ref{1:0}, \ref{1:2} and \ref{1:1}, we have
$\limsup_{\e\daw0} \E(t_\e)/\e \leq \lim_{\e\daw0} \E(t_\e^r \wedge t_\e^l)/\e=\sqrt{2}/(2\tau)$. By (\ref{ref::3333}), Property (3)-(ii) and hence Lemma \ref{intermed} follow. We conclude that $(\W,\W(\tau))$ has the required distribution.

\subsubsection{Markov Property and Stationarity}
We continue with the second and third properties of Theorem \ref{real::stick}.
$\W(\tau_2)$ is constructed by modifying 
$\W=\W(\tau=0)$ according to the marking $\M(\tau_2)$.
In order to prove the Markov property and stationarity, it suffices to prove that this is 
distributionally
equivalent to the 
following procedure: (1) construct $\W(\tau_1)$ from $(\W,\M(\tau_1))$; then (2) construct $\W(\tau_2)$ from $(\W(\tau_1),\M^{\tau_1}(\D \tau))$
where $\M^{\tau_1}(\D \tau)$ is a marking of $\W(\tau_1)$ with intensity $\D \tau \equiv\tau_2-\tau_1$ which, given
the past $(\W,\{\M(\tau)\}_{\tau\leq\tau_1})$, only depends on $\W(\tau_1)$ 
.

\bigskip
 
Recall that given $\W$,
\begin{enumerate}
\item[(i)] for any measurable subset $O\subset\R^2$ with $\Lo(O)<\infty$ (where $\Lo$ is the local time outer measure---see Definition \ref{localtimemeasure}), $[\Mde\setminus\Mun] \cap O$ is a Poisson Point Process on $\R^2$ with intensity measure $(\tau_2-\tau_1)\Lo(\cdot \cap O)$, and 
\item[(ii)] $\{\M(\tau)\}_{\tau\leq\tau_1}$ and $\tilde\M(\D\tau)\equiv\Mde\setminus\Mun$ are independent.
\end{enumerate}
$\tilde \M(\D\tau)$ induces a natural marking on $\W(\tau_1)$.
Indeed, for every $n\geq0$, we can define
$\tilde \M_{n,n}^{\tau_1}(\D \tau)$ as
$\tilde \M(\D\tau) \cap E_n$ where 
$E_n=\{B^{\tau_1}_i\}_{i=0}^{n-1} \cap \{\hat B_j^{\tau_1}\}_{j=0}^{n-1}$ and $\M^{\tau_1}(\D\tau)\equiv\lim_{n\uaw\infty}\M^{\tau_1}_{n,n}(\D\tau)$. We will denote by $\W_{n,n}^{\tau_1}(\D\tau)$ the web obtained from $\W(\tau_1)$ by
switching the direction of all the points in $\M_{n,n}^{\tau_1}(\D \tau)$.

We have already proved that $\W(\tau_1)$ is a Brownian web. 
Hence, $\mathcal{\Lo}(E_n)<\infty$ and, by item (i) above,
conditioned on $\W$, $\M_{n,n}^{\tau_1}(\D \tau)$ is a Poisson Point Process with intensity measure $(\tau_2-\tau_1)\Lo(\cdot \cap E_n)$. 

\blem
\label{LO::tau1}
Let $\Lo^{\tau_1}_{n,n}$ be the local time measure on $\R^2$
induced by $\{B^{\tau_1}_i\}_{i=0}^{n-1} \cup \{\hat B_j^{\tau_1}\}_{j=0}^{n-1}$, i.e.,
$$\mathcal{L}^{\tau_1}_{n,n}(O)=m_\phi\left(\mathcal{P}(\{B_i^{\tau_1}\}_{i=0}^{n-1} \bigcap \{\hat B^{\tau_1}_j\}_{j=0}^{n-1} \cap O)\right)$$
(where $\mathcal{P}$ is the projection on the $t$-axis).
Then $\Lo(O\cap E_n)=\mathcal{L}_{n,n}^{\tau_1}(O)$, where $\Lo$ is the usual local time
measure of (\ref{lo-def}). 
\elem
\begin{proof}
For a web $\W'$, let $\W'\cap\hat\W'$
denote the set of $(1,2)$ points of $\W'$. 
By definition,
\begin{eqnarray*}
\mathcal{L}^{\tau_1}_{n,n}(O) & = &  m_\phi\left(\mathcal{P}(E_n\cap O)\right)
\\
\Lo(O\cap E_n)& = & m_\phi\left(\mathcal{P}([\W\cap\hat \W]\cap E_n\cap O)\right)
\end{eqnarray*}
Hence, in order to prove our lemma it is
sufficient to prove that for every Borel $O$ 
$$
m_\phi\left(\mathcal{P}([E_n\cap O]\setminus[\W\cap\hat \W])\right)=0
$$
which will follow if we can prove that
\begin{equation}\label{e-mphi}
m_\phi\left(\mathcal{P}([\W(\tau_1)\cap\hat\W(\tau_1)]\setminus[\W\cap\hat \W])\right)=0.
\end{equation}
In order to prove (\ref{e-mphi}) we prove
\begin{equation}\label{e-mphi2}
m_\phi\left(\mathcal{P}([\W\cap\hat \W]\setminus[\W(\tau_1)\cap\hat\W(\tau_1)])\right)=0
\end{equation}
instead. The lemma will follow from
the equidistribution of $(\W,\W(\tau_1))$ and $(\W(\tau_1),\W)$. (Recall
that in Subsections \ref{proof::stickiness::webs}-\ref{Lemma::esdd} we 
already proved that $(\W,\W(\tau_1))$ is a sticky pair of webs whose
distribution is invariant under permutation of the two webs.)

We now prove (\ref{e-mphi2}). 
For a given realization of $(\W,\W(\tau_1))$, let us assume that 
$$m_\phi\left({\cal P}([\W\cap\hat\W]\setminus[\W(\tau_1)\cap \hat\W(\tau_1)])\right)>0$$ and find a contradiction. By construction of $\M(\tau_2)$, there would be strictly
positive probability that $\Mde\setminus[\W(\tau_1)\cap\hat\W(\tau_1)]\neq\emptyset$. 

Let $z$ be 
any point in $\Mde$. Then $z$ is a separation point of $\Net(\tau_2)$. Proposition \ref{SP::new}(3) directly implies that for some $i$ the path $B_i^{\tau_1}\in\W(\tau_1)$, from $z_i$, enters $z$. 
Since up to a reversal of the $t$-axis $\Net(\tau_2)$ and $\hat \Net(\tau_2)$ are equidistributed, there is a path $\hat B^{\tau_1}_j\in\hat \W(\tau_1)$ meeting $B^{\tau_1}$ at $z$ and hence that $z$ is in $\W(\tau_1)\cap\hat\W(\tau_1)$. It would follow that $\Mde\subset\W(\tau_1)\cap\hat\W(\tau_1)$, yielding a contradiction.
This ends the proof of the lemma.
\end{proof}

Lemma \ref{LO::tau1}  implies that
$\M_{n,n}^{\tau_1}(\D \tau)$ only depends on $\W(\tau_1)$. Moreover, $\W(\tau_1)$ being a Brownian web, we also have the distributional identities, 
\begin{eqnarray}
(\W(\tau_1),\M_{n,n}^{\tau_1}({\D\tau}))=_{d}(\W,\tilde \M_{n,n}({\D \tau})), \label{w1}
\\
(\W(\tau_1),\W_{n,n}^{\tau_1}({\tau_2}))=_{d}(\W,\W_{n,n}({\D \tau})), \label{w2}
\end{eqnarray}
where $\W_{n,n}(\D\tau)$ and $\tilde \M_{n,n}(\D \tau)$ are defined as in Section \ref{marking_process}.
It remains to prove that $\W_{n,n}^{\tau_1}(\D\tau)$ converges (in $(\mathcal{H},d_{\mathcal{H}})$) to $\W(\tau_2)$.

\blem
$\M^{\tau_1}(\D\tau)$ and $\tilde \M(\D\tau)$ coincide.
\elem
\begin{proof}
By construction, $\M^{\tau_1}(\D\tau)\subset\tilde\M(\D\tau)$ since  $\M^{\tau_1}(\D \tau)$ is
the marking induced by $\tilde \M(\D\tau)$ on $\W(\tau_1)$.
Analogously, we define $\M'(\D \tau)$
as the marking induced by $\M^{\tau_1}(\D \tau) \ (=\lim_{n\uaw\infty}\M_{n,n}^{\tau_1}(\D\tau))$ on $\W$. 
We already proved that $(\W,\W(\tau_1))$ is a  $(1/2\tau)$-sticky pair of webs. 
Therefore, $(\W,\W(\tau_1))$ is equidistributed with $(\W(\tau_1),\W)$ 
and, by (\ref{w1}), $(\W,\W(\tau_1),\M^{\tau_1}(\Delta\tau))$ is
equidistributed with $(\W(\tau_1),\W,\tilde \M({\D\tau}))$.
Thus, $$(\W,\M_{n,n}'(\D\tau))=_d (\W(\tau_1),\M_{n,n}^{\tau_1}(\D\tau)).$$
By (\ref{w1}), we conclude
that
$\M'_{n,n}(\D \tau)$ 
is distributed like $\tilde \M_{n,n}(\D\tau)$. Since by construction, $\M'_{n,n}(\D \tau)\subset\tilde \M_{n,n}(\D\tau)$, it follows that $\M'_{n,n}(\D \tau)= \tilde \M_{n,n}(\D\tau)$
and $\M'(\D\tau)=\tilde \M(\D\tau)$. Since $\M'(\tau)\subset\M^{\tau_1}(\D \tau)$, we deduce that $\tilde \M(\D\tau)\subset \M^{\tau_1}(\D\tau)$ and hence $\M^{\tau_1}(\D\tau)=\tilde \M(\D\tau)$.
\end{proof}

Let $\W'$ be any subsequential limit
of $\{\W_{n,n}^{\tau_1}(\D\tau)\}$. We next prove that $\W'=\W(\tau_2)$
via two inclusions, which completes this subsection.

{\bf(i) $\W'\subseteq \W(\tau_2)$}.

Let $z=(x,t)\in\M(\tau_2)$, and let $\pi\in\W'$ start at $t-2\e$ with $\e>0$, and pass through $z$.
By Proposition \ref{existence_of_Br}, what we need to 
show is that $\pi\sim^z B_{switch}$. 
By construction, there exists a sequence $\{\pi_N\}_{N\geq0}$ so that $\pi_{N}$ belongs to $\cup_{n,m>M}\W_{n,m}^{\tau_1}(\D\tau)$ and $\{\pi_N\}$ converges to $\pi$.
Taking $N$ large enough, we can assume w.l.o.g. that $\pi_N$ belongs to $\Net_{\leq t-\e}$ and $(x,t)\in\mathcal{M}_{n,m}(\tau_2)$  
for $n,m>N$. Moreover, by Proposition \ref{SP::new}(2) we can also assume that 
$\pi_N$ enters the point $z$. We distinguish between two cases.
\begin{enumerate}
\item $z\in\mathcal{M}(\tau_1)$. Here, $z\notin\M^{\tau_1}(\D\tau)$ and by construction, $\pi_N \sim^z B^{\tau_1}$. Since $B^{\tau_1} \sim^z B_{switch}$, Proposition \ref{SP::new}(1) implies that $\pi_N \sim^z B_{switch}$. By Proposition \ref{SP::new}(2), $\pi\sim^z B_{switch}$. 
\item $z\in \tilde \M(\D\tau)$. Since $\M^{\tau_1}(\D\tau)=\tilde\M(\D\tau)$,
 we get that
$\pi_N\sim^z B_{switch}^{\tau_1}$. We claim that $B_{switch}^{\tau_1} \sim^z B_{switch}$, implying that
$\pi \sim^z B_{switch}$ as desired. The claim can be verified as follows. Let us assume that $B_{switch}^{\tau_1}\sim^z B$ (and show that this leads to a contradiction). Then $B^{\tau_1}\sim^z B_{switch}$, implying that $B$ and $B^{\tau_1}$ separate at $z$, or equivalently that $z\in\M(\tau_1)$. Since $\M(\tau_1)$ and $\tilde \M(\D\tau)=\M(\tau_2)\setminus\M(\tau_1)$ are disjoint, the claim follows.
\end{enumerate}

{\bf(ii) $\W'\supseteq \W(\tau_2)$}.

There is at least one path $B'$ in $\W'$ starting from $z_i$. By Proposition \ref{existence_of_Br}(3)(o), there is a unique path $B^{\tau_2}_i\in\W(\tau_2)$ starting from there. Since  $\W'\subseteq \W(\tau_2)$ we get $B'=B^{\tau_2}_i$. Hence, $\W'\supseteq \overline{\{B_i^{\tau_2}\}}=\W(\tau_2)$ (see Proposition
\ref{existence_of_Br}(3)(ii)).

\subsubsection{ $\tau\raw B_0^\tau(t)$ is Piecewise Constant}

For any $\tau\leq\tau_0$, the path $B_0^\tau$ belongs to $\Net(\tau_0)$. 
Given $\Net(\tau_0)$,
$B_0^\tau(t)$ only depends on the direction of the $(1,2)$ points of $\W(\tau)$ which are located at the $(0,t)$-separation points of  $\Net(\tau_0)$. Since the set of $(0,t)$-separation points 
is locally finite (see Proposition \ref{(0,1)S.P_locally_finite}),
$\tau\raw B_0^\tau(t)$ is piecewise constant.

\bigskip

{\bf {\em Acknowledgements.}}
The research reported in this paper was supported in part
by N.S.F. grants DMS-01-04278 and DMS-06-06696. The authors
thank Rongfeng Sun and Jan Swart for very useful
communications and discussions concerning their work on the Brownian net.
They also thank the anonymous referee for a number of very
helpful suggestions which have led to a clearer presentation
of some of our arguments.

\end{document}